\documentclass[a4paper,11pt,onecolumn]{article}

\usepackage{mathbbold}
\usepackage{bbm}
\usepackage{amsfonts}

\usepackage{geometry}
\usepackage{fancyhdr}
\usepackage{fancyvrb}
\usepackage{fancybox}
\usepackage{amsmath,amsfonts,amssymb,amscd,graphicx}
\usepackage{subfigure}
\usepackage{indentfirst}
\usepackage{bm}
\usepackage{multicol}
\usepackage{abstract}
\usepackage{anysize}
\usepackage{hyperref}
\usepackage{listings}
\usepackage{color}
\usepackage{units}
\usepackage{tabularx}
\usepackage{mathrsfs}
\usepackage[all]{xy}

       \font\tenmsb=msbm10
       \font\sevenmsb=msbm7
       \font\fivemsb=msbm5
       \catcode`\@=11
       \ifx\amstexloaded@\relax\catcode`\@=\active
       \endinput\else\let\amstexloaded@\relax\fi
       \def\spaces@{\space\space\space\space\space}
       \def\spaces@@{\spaces@\spaces@\spaces@\spaces@\spaces@}
       \def\space@.{\futurelet\space@\relax}
       \space@. %
       \def\Err@#1{\errhelp\defaulthelp@\errmessage{AmS-TeX error: #1}}
       \def\relaxnext@{\let\next\relax}
       \def\accentfam@{7}
       \def\noaccents@{\def\accentfam@{0}}

       \def\mathcal{\relaxnext@\ifmmode\let\next\mathcal@\else
       \def\next{\Err@{Use \string\mathcal\space only in math mode}}\fi\next}
       \def\mathcal@#1{{\mathcal@@{#1}}}
       \def\mathcal@@#1{\noaccents@\fam\tw@#1}

       \def\Bbb{\relaxnext@\ifmmode\let\next\Bbb@\else
       \def\next{\Err@{Use \string\Bbb\space only in math mode}}\fi\next}
       \def\Bbb@#1{{\Bbb@@{#1}}}
       \def\Bbb@@#1{\noaccents@\fam\msbfam#1}

       \newfam\msbfam
       \textfont\msbfam=\tenmsb
       \scriptfont\msbfam=\sevenmsb
       \scriptscriptfont\msbfam=\fivemsb

\def\notin{ \in \! \!\!\!\!\!  / }

\def\Z{\mathbb{Z}}
\def\R{\mathbb{R}}
\def\C{\mathbb{C}}

\def\H{\mathbb{H}}

\def\V{\mathbb{V}}

\def\qed{$\rlap{$\sqcap$}\sqcup$}

\newtheorem{Theorem}{Theorem}[section]
\newtheorem{Lemma}{Lemma}[section]
\newtheorem{Proposition}{Proposition}[section]
\newtheorem{Corollary}{Corollary}

\newtheorem{Definition}{Definition}[section]

\renewcommand{\theequation}{\thesection.\arabic{equation}}
\@addtoreset{equation}{section}

\newcommand{\bq}{\begin{equation}}
\newcommand{\eq}{\end{equation}}

\begin{document}

\title{Relatively Open Gromov-Witten Invariants for Symplectic Manifolds
of Lower Dimensions}

\author{Hai-Long Her}

\date{}

\maketitle

\begin{abstract}
Let $(X,\omega)$ be a compact symplectic
manifold, $L$ be a Lagrangian submanifold and $V$ be a codimension 2
symplectic submanifold of $X$, we consider  pseudoholomorphic
maps from a Riemann surface with a fixed conformal structure and with boundary $(\Sigma,\partial\Sigma)$
to the pair $(X,L)$ satisfying Lagrangian boundary conditions and
intersecting $V$. In some special  cases, for instance, under the
semi-positivity condition, we
study the stable moduli space of such open pseudoholomorphic maps
involving the intersection data. If $L\cap V=\emptyset$, we
study the problem of orientability of the moduli space. Moreover, assume that there exists an
anti-symplectic involution $\phi$ on $X$ such that $L$ is the fixed
point set of $\phi$ and $V$ is $\phi$-anti-invariant, then we define
the so-called ``relatively open" invariants for the tuple
$(X,\omega,V,\phi)$ if $L$ is orientable and dim$X\le 6$. If $L$ is
nonorientable, we define such invariants under the condition that
dim$X\le4$ and some additional restrictions on the number of marked
points on each boundary component of the domain.
\end{abstract}




\bigskip

\tableofcontents

\section{Introduction}\label{SEC-1}

Gromov-Witten invariants $GW_{g,n}(X,A)$ for a symplectic
manifold $(X,\omega)$, roughly speaking, count with sign the $(J,\nu)$-stable maps $u$
from a compact Riemann surface $\Sigma_{g,n}$ without boundary of genus $g$ with
$n$ marked points to $(X,\omega)$,
representing the homology class $A\in H_2(X)$ and satisfying some
constraints, where $J$ is a $\omega$-compatible or $\omega$-tamed
almost complex structure on $X$ and $\nu$ is an inhomogeneous
perturbation satisfying the $J$-holomorphic map equation
$\bar{\partial}_Ju=\nu$. Such invariants originate from the seminal works of Gromov and Witten as well as Ruan's
introducing Donaldson type invariants in symplectic category (see
\cite{G}\cite{W}\cite{R1}). The rigorous mathematical foundation of GW invariants
was first established by Ruan-Tian \cite{RT1} for semi-positive
symplectic manifolds in 1993. Then many mathematicians contribute to
the development and completion of such theme (see
\cite{FO}\cite{KM}\cite{LiT}\cite{RT2}\cite{Si}). In many cases
these numbers coincide with the enumerative invariants in algebraic
geometry. In order to find effective ways of computing these
invariants, Tian \cite{Ti} first showed a rough description of
studying the degeneration of rational curves under a symplectic
degeneration, an operation of `splitting the target'. Then
Li-Ruan \cite{LR} and Ionel-Parker \cite{IP1}\cite{IP2}
independently defined the so-called relative Gromov-Witten
invariants with respect to a codimension 2 symplectic divisor
$V$ and also established the symplectic connect sum formula for GW
invariants. In the current paper, we construct certain relative GW invariants of open version (see Theorem \ref{thm-1.1} and \ref{thm-1.2} below) which in some sense count pseudo-holomorphic maps from Riemann surfaces with boundary  and relative to a divisor.

For recent years, theoretic physicists have predicted the existence
of enumerative invariants about pseudoholomorphic maps with
Lagrangian boundary conditions by studying dualities for open
strings (see \cite{AKV}\cite{OV}). Roughly speaking, let $L \subset
X$ be a Lagrangian submanifold, such invariants would count
(perturbed) pseudoholomorphic maps $u$ (or say $(J,\nu)$-maps
satisfying equation $\bar{\partial}_Ju=\nu$) from a Riemann surface
with boundary $(\Sigma,\partial\Sigma)$ to $(X,L)$, representing a
relative class $d\in H_2(X,L)$, such that the boundary is mapped in
the Lagrangian submanifold.
We denote by ${\mathcal M}_{k,l}(X,L,d)$ the moduli
space of such maps $u$ with $k$ distinct marked points in
$\partial\Sigma$ and $l$ distinct marked points in $\Sigma$, and its
compactfication by $\overline{\mathcal M}_{k,l}(X,L,d)$. If we consider the special case that the domain surface with fixed conformal structure, to prescribe the domain surface, we will denote  by ${\mathcal M}_{k,l}(\Sigma, X,L,d)$ and $\overline{\mathcal M}_{k,l}(\Sigma,X,L,d)$ the corresponding moduli spaces, respectively.
Note that we have canonical evaluation maps at marked points
$$evb_i:\overline{\mathcal M}_{k,l}(X,L,d)\rightarrow L,\ \ \ \ i=1,\cdots,k,$$
$$evi_j:\overline{\mathcal M}_{k,l}(X,L,d)\rightarrow X,\ \ \ \ j=1,\cdots,l.$$

There are two main difficulties in defining such invariants:
orientability of moduli space and codimension 1 bubbling-off
phenomenon. Note, by \cite{FO3}\cite{Sil}, that ${\mathcal M}_{k,l}(X,L,d)$
might be non-orientable. If $L$ is relatively spin, then Fukaya {\it
et al.} proved that the moduli space is orientable.
Katz-Liu \cite{KL}\cite{Liu} defined an open invariant under the
much restricted assumption that there exists an additional
$S^1$-action on the pair $(X,L)$. To define an orientation on the
moduli space, they also assumed the Lagrangian submanifold $L$ is
orientable and (relatively) spin. The orientability of moduli space
seems important for Katz-Liu's definition of enumerative invariants.
However, Solomon \cite{So} showed that even if $L$ is not
orientable, under the weaker assumption that $L$ is ``relatively
Pin$^\pm$'' and some constraints on the number of boundary marked
points, we still can obtain a canonical isomorphism from the
orientation bundle of ${\mathcal M}_{k,l}(\Sigma,X,L,d)$ to the pull-backed
line bundle from the orientation bundle of $L^k$ by the evaluation
map $\prod_ievb_i$. Although the moduli space ${\mathcal
M}_{k,l}(\Sigma,X,L,d)$ may be non-orientable, the integral of $\det
(TL)$-valued forms over ${\mathcal M}_{k,l}(\Sigma,X,L,d)$ would make
sense.

On the other hand, we at present have no effective method to
completely deal with the codimension 1 boundary of moduli space
coming from the bubbling off discs. According to the discussion in
\cite{AKV}, in order to define an invariant independent of other
choices, it seems necessary to introduce some additional parameter
on $(X,L)$. For instance, the assumption in \cite{Liu} that there
exists an $S^1$-action is used to cancel the effect of codimension 1 boundary
of moduli space. In the paper \cite{So}, under the assumption that
there exists an anti-symplectic involution
$$\phi:X\rightarrow X,\ \ \ \ \phi^*\omega=-\omega,$$ such that
$L={\rm Fix}(\phi)$, and among other technical conditions, Solomon
constructed the open invariants for ${\rm dim}X\le 6$ if $L$ is
orientable and, for ${\rm dim}X\le 4$ even if $L$ might be
non-orientable. Note that, from the viewpoint of real algebraic
geometry, the assumption that symplectic manifold admits an
anti-symplectic involution is more natural. Actually, in that paper,
Solomon also showed that, for the genus zero domain surface and
strongly semi-positive real symplectic manifold of dimension no more
than 6, his open invariants exactly coincide with twice of the
Welschinger's invariants which can be regarded as a lower bound of
the number of real rational curves in a real symplectic manifold(see
\cite{W1}\cite{W2}\cite{W3}). Independently, Cho \cite{C}, with
different way of counting, also defined similar enumerative
invariants for strongly semi-positive symplectic manifolds with
dim$X\le 6$, genus $g=0$, and $L$ being relatively spin.

The aim of this article is to show a definition of ``relatively open
Gromov-Witten invariants"  which intuitively would count stable maps
$u:(\Sigma,\partial\Sigma)\rightarrow(X,L)$ from a compact Riemann
surface with boundary and {\bf with fixed conformal structure} whose images intersect with a codimension 2
symplectic submanifold $V\subset X$, satisfying Lagrangian boundary
conditions. In the current paper, to avoid technical complexity of
constructing virtual cycle or Kuranishi structure, we do not deal
with the most general case. Instead, we only consider the special
case that all moduli spaces are of expected dimensions. That is we
need the semi-positivity assumptions\footnote{See the Definition \ref{semiposi}} on both  $(X,L)$
and $V$ (when $L_V=L\cap V\neq\emptyset$, we have to further assume that the pair $(V,L_V)$ is also semi-positive), or the assumption that each considered moduli space
generically contains no multiply covered maps. These invariants
are designed for a preparation of establishing open symplectic sum
formulas, which will be used to compute open Gromov-Witten
invariants of a symplectic connect sum of two symplectic manifolds
$(X,L_X,V)$ and $(Y,L_Y,\overline{V})$. More concretely, the
symplectic sum is an operation that first removes $V$ and
$\overline{V}$ from $X$ and $Y$, respectively, and then combines
them into a new symplectic manifold $X\# Y$ with symplectic
structures matching on the overlap region. A stable map into the sum
is expected to be a pair of stable maps into the two sides which
match in the middle. So the first step is to count stable maps in
one side $X$ which record the intersection points with $V$ and
multiplicities. At present, we only construct such invariants under
the additional assumptions that $L$ is the fixed point set of an
anti-symplectic involution, dim$X\le 6$ and $L\cap V=\emptyset$.
However, for possible applications in the future, we will define the
related moduli space and study the problem of orientability for more
general situation. An axiomized formulation of the symplectic sum of
open GW invariants has appeared in \cite{LiuY}.

Similar to the absolute case in \cite{IP1}, before considering such
stable maps we have to extend $J$ and $\nu$ to the connect sum. The
$V$-$compatibility$ conditions imposed on the pair $(J,\nu)$ defined
in Section \ref{SEC-4} ensure that such an extension exists.
However, these conditions do not always hold for generic $(J,\nu)$.
Our relatively open invariants count stable maps for these special
$V$-compatible pairs, which are different from the way of counting
of the absolutely open GW invariants. Given such a special
$V$-compatible pair $(J,\nu)$, $V$ is a $J$-holomorphic submanifold,
and a $(J,\nu)$-holomorphic map $u:(\Sigma,\partial\Sigma)\rightarrow
(V,L_V)$ (if $L_V=L\cap V\neq \emptyset$) into $(V,L_V)$ is
automatically a $(J,\nu)$-holomorphic map into $(X,L)$. Therefore, some
domain components (closed or open) of stable maps maybe mapped
entirely into $V$. Such maps are not transverse to $V$ and the
moduli space of such maps may be of dimension larger than the
expected dimension of ${\mathcal M}_{k,l}(\Sigma, X,L,d)$.

To avoid such non-transversal intersections, we restrict  consideration to the so-called $V$-regular open stable maps which
have no components mapped entirely into $V$. In the current paper, we only study these maps for the special case that $L\cap V=\emptyset$
which avoids the appearance of boundary intersections. The author plans to study the more difficult case $L\cap V\neq\emptyset$ in a separate paper.
We note that such maps may intersect $V$ at finite many interior points with multiplicity. According to the ordering of
these intersection points, the space of $V$-regular stable maps separates into components labeled by vector
$\mathbbm{s}=(s_1,\cdots,s_\mathbbm{l})$. The subscript $\mathbbm{l}$ denotes the number of  interior intersection points and $s_j$ denotes
the multiplicity of the $j^{th}$ interior intersection point. Then in Section \ref{SEC-5}, we will
study the moduli space ${\mathcal M}^V_{k,l}(\Sigma, X,L,d)$ of such
$V$-regular stable maps. We will prove that each irreducible (see Definition \ref{multi-cover})
component ${\mathcal M}^{V,\mathbbm{s}}_{k,l,\mathbbm{l}}(\Sigma, X,L,d)$ of $V$-regular maps is an
orbifold with boundary whose dimension is expressed in (\ref{dim}). This applies
to the semi-positive case since then those moduli spaces are generically irreducible.

In Sections \ref{SEC-7} and \ref{SEC-8}, we consider the
compactification of the moduli space of $V$-regular open maps---the
space of $V$-stable open maps, of which each component is a subset
of the closure of the $V$-regular space ${\mathcal
M}^{V,\mathbbm{s}}_{k,l,\mathbbm{l}}(\Sigma, X,L,d)$ in the stable moduli space
$\overline{\mathcal M}_{k,l+\mathbbm{l}}(\Sigma, X,L,d)$. We will show that the irreducible part of this compactification is also an
orbifold with boundary (corners). The key observation is that each sequence
$\{u_n\}$ of $V$-regular open maps limits to a stable map $u$ with
additional restriction. For instance, if the image of a component of
such a limit map lies entirely in $V$, then along this component, we can find a section $\xi$ of the normal bundle of
$V$, such that the elliptic equation $D^N\xi=0$ holds, where $D^N$ is the restriction of the linearization operator
of $\bar{\partial}_{J,\nu}$ to the normal bundle of $V$. Stable maps with this additional restriction are called
$V$-stable maps. Each component of the space of $V$-stable maps is denoted by $\overline{\mathcal M}^{V,\mathbbm{s}}_{k,l,\mathbbm{l}}(\Sigma,X,L,d)$, which is the
compactification of the space of $V$-regular open maps with frontier
strata of codimension at least one. To analyze the convergence of
sequences of $V$-regular maps, we refer to  the renormalization technique
used in \cite{IP1} and \cite{T} for closed curves.
For our open case, the parallel arguments can go through and we will omit the detail.
Since the $V$-stable moduli space contains boundary or corners,  in section \ref{SEC-gl} we will describe a gluing construction which gives the moduli space a local chart.

Then we come to define relative invariants for the case $L\cap
V=\emptyset$. Let us introduce some more delicate notations. Recall
that the maps we consider are $(j,J,\nu)$-holomorphic maps $u$ from
a bordered Riemann surface $(\Sigma,\partial\Sigma)$ with fixed
conformal structure $j$ to a symplectic-Lagrangian pair $(X,L)$
satisfying Lagrangian boundary conditions and representing a fixed
relative homology class $d\in H_2(X,L)$. We suppose the boundary of
$\Sigma$ has $m$ components, $i.e.$
$\partial\Sigma=\bigcup_{a=1}^m(\partial\Sigma)_a$ with each
$(\partial\Sigma)_a\simeq S^1$. Moreover, we also require that the
image of each boundary component represents a fixed homology class,
$i.e.$ $u|_{(\partial\Sigma)_a*}([(\partial\Sigma)_a])=d_a\in
H_1(L)$. Suppose that there are $k_a$ marked points
$z_{a1},\cdots,z_{ak_a}$  on each boundary component
$(\partial\Sigma)_a$ and $l$ marked points $w_1,\cdots,w_l$ and
additional $\mathbbm{l}$ intersecting marked points
$q_1,\cdots,q_\mathbbm{l}$ on $\Sigma$.\footnote{To avoid some kind
of bubbling, in this article, we require that if $\Sigma\simeq D^2$,
$k>0$.} We reset
$${\bf d}=(d,d_1,\cdots,d_m)\in H_2(X,L)\oplus H_1(L)^{\oplus m},$$
$${\bf k}=(k_1,\cdots,k_m),\ \ \ |{\bf k}|=\sum_{a=1}^mk_a,$$
$${\bf
u}=(u,\vec{z},\vec{w},\vec{q}),\ \ \vec{z}=(z_{ai}),\ \
\vec{w}=(w_j),\ \ \vec{q}=(q_\jmath).$$ Then we rewrite the moduli
space of $V$-regular maps ${\bf u}$ by ${\mathcal M}^{V,\
\mathbbm{s}}_{{\bf k},l,\mathbbm{l}}(\Sigma, X,L,{\bf d})$. The
compactification $\overline{\mathcal M}^{V,\mathbbm{s}}_{{\bf
k},l,\mathbbm{l}}(\Sigma, X,L,{\bf d})$ is the space of $V$-stable open
maps.

There are some natural evaluation maps at those marked and
intersection points
$$evb_{ai}:\overline{\mathcal
M}^{V,\mathbbm{s}}_{{\bf k},l,\mathbbm{l}}(\Sigma, X,L,{\bf d})\rightarrow
L,\ \ \ i=1,\cdots,k_a,\ a=1,\cdots,m,$$
$$evi_j:\overline{\mathcal
M}^{V,\mathbbm{s}}_{{\bf k},l,\mathbbm{l}}(\Sigma, X,L,{\bf d})\rightarrow
X,\ \ \ j=1,\cdots,l,$$
$$evi^I_{\jmath}: \overline{\mathcal
M}^{V,\mathbbm{s}}_{{\bf k},l,\mathbbm{l}}(\Sigma, X,L,{\bf d})\rightarrow
V,\ \ \ \jmath=1,\cdots,\mathbbm{l}.$$

The first problem is the orientability of the moduli space. To get a
more general result, we do not expect $\overline{\mathcal M}^{V,\
\mathbbm{s}}_{{\bf k},l,\mathbbm{l}}(\Sigma, X,L,{\bf d})$ to have a
canonical orientation, which means that there is a canonical
orientation on the orientation line bundle $\det(\overline{\mathcal
M}^{V,\mathbbm{s}}_{{\bf k},l,\mathbbm{l}}(\Sigma, X,L,{\bf d}))$.
Instead, we only want to construct an orientation on a modified line
bundle $${\rm det}(T\overline{\mathcal M}^{V,\mathbbm{s}}_{{\bf
k},l,\mathbbm{l}}(\Sigma, X,L,{\bf d}))
\otimes\bigotimes_{a,i}evb_{ai}^*\det(TL).$$ The conditions that
ensure there exists such an orientation are that the Lagrangian
submanifold $L$ is ``relatively $Pin^\pm$" (see Definition
\ref{def-9.1}) and some restrictions imposed on the  boundary marked
points. To state the following theorems, we recall

\begin{Definition}\label{semiposi}
(1) A $2n$-dimensional symplectic manifold $(M,\omega)$ is called
semi-positive if $\omega(\beta)\le 0$ for any $\beta\in\pi_2(M)$
with $3-n\le c_1(M)[\beta]<0$, where $c_1(M)$ is the first Chern
class of $M$.\ \ \ (2) Let $L$ be a $n$-dimesional Lagrangian
submanifold of $(M,\omega)$. A symplectic-Lagrangian pair $(M,L)$ is
called semi-positive if $\omega(\beta)\le 0$ for any
$\beta\in\pi_2(M,L)$ with $3-n\le \mu_L(\beta)<0$, where $\mu_L$ is
the Maslov index.
\end{Definition}

In Section \ref{SEC-9} we obtain the following result.

\begin{Theorem}\label{thm-1.1}
Assume that both $(X,L)$ and $V$ are semi-positive or the moduli
space $\overline{\mathcal M}^{V,\mathbbm{s}}_{{\bf
k},l,\mathbbm{l}}(\Sigma, X,L,{\bf d})$ generically contains no multiply
covered maps\footnote{See Definition \ref{multi-cover}. }, $L$ is relatively $Pin^\pm$ and fix a relatively
$Pin^\pm$ structure on $(X,L)$. If $L$ is not orientable, we assume
$k_a\cong w_1(d_a)+1 \mod 2$. If $L$ is orientable, fix an
orientation. Then the relatively $Pin^\pm$ structure on $(X,L)$ and
the orientations of $L$ if it is orientable determines a canonical
isomorphism
\begin{equation}\label{thm-1.1-isom}
{\rm det}(T\overline{\mathcal M}^{V,\mathbbm{s}}_{{\bf
k},l,\mathbbm{l}}(\Sigma, X,L,{\bf d}))
\widetilde{\longrightarrow}\bigotimes_{a,i}evb_{ai}^*\det(TL).
\end{equation}
\end{Theorem}
{\it Remark}. By the Wu relation \cite{Mi-St}, if $n\le 3$ then $L$
is always $Pin^-$.

\bigskip

Given the isomorphism (\ref{thm-1.1-isom}), we can define the
relatively open invariants as follows. Denote the images of marked
points under evaluation maps by
$$x_{ai}=u(z_{ai}),\ y_j=u(w_j),\ \mathbbm{q}_\jmath=u(q_\jmath)\ .$$
Let $\Omega^*(L,{\rm det}(TL))$  denote differential forms on $L$
with values in ${\rm det}(TL)$, and let $\Omega^*(X)$ (resp.
$\Omega^*(V)$) denote ordinary differential forms on $X$ (resp.
$V$). Let $\alpha_{ai}\in\Omega^n(L,{\rm det}(TL))$, $a=1,\cdots,m;
i=1,\cdots,k_a$, represent the Poincar\'{e} dual of the point
$x_{ai}$ in $H^n(L,{\rm det}(TL))$, which is the cohomology of $L$
with coefficients in the flat line bundle ${\rm det}(TL)$. Let
$\gamma_j\in\Omega^{2n}(X)$ represent the Poincar\'{e} dual of $y_j$
for $j=1,\cdots,l$. And let $\eta_\jmath\in\Omega^{2n-2}(V)$
represent the Poincar\'{e} dual of $\mathbbm{q}_\jmath$ for
$\jmath=1,\cdots,\mathbbm{l}$. Then we define
\begin{eqnarray}\label{ROGW-d-inv}
&&\mathscr{R}\mathscr{N}:=\mathscr{R}\mathscr{N}(V,{\bf d},{\bf k},l,\mathbbm{l})\nonumber\\
&&\ \ \ \ \ \ \ \ =\int_{\overline{\mathcal M}^{V,\mathbbm{s}}_{{\bf
k},l,\mathbbm{l}}(\Sigma, X,L,{\bf d})
}\bigwedge_{a,i}evb_{ai}^*(\alpha_{ai}) \bigwedge_{j}
evi^*_j(\gamma_j)\bigwedge_{\jmath} evi^{I*}_\jmath(\eta_\jmath)\ .\
\ \ \ \ \ \ \
\end{eqnarray}

The integral (\ref{ROGW-d-inv}) makes sense because by Theorem
\ref{thm-1.1}, the integrand is a differential form taking values in
the orientation line bundle of the $V$-stable moduli space. Denote
by $\mu:H_2(X,L)\rightarrow\Z$ the Maslov index, and by $g$ the
genus of the closed Riemann surface
$\Sigma\cup_{\partial\Sigma}\overline{\Sigma}$ which is the complex
double of $\Sigma$. From the dimensional calculation (\ref{dim})  we
know that the integral above vanishes unless
\begin{eqnarray}\label{indexthm-condition}
(|{\bf k}|+2l)n+2\mathbbm{l}(n-1)&=&\mu(d)+n(1-g)+k \nonumber\\
   & &+2(l+\mathbbm{l}-{\rm
deg}\mathbbm{s})-{\rm dim} Aut(\Sigma).
\end{eqnarray}
where ${\rm deg}\mathbbm{s}=\sum_{j=1}^\mathbbm{l}s_j$.

In general, the integral (\ref{ROGW-d-inv}) might depend on the
choice of the forms $\alpha_{ai}$, $etc.$, because the codimension 1
boundary strata would contribute to the integral. To prove the
invariance of the integral (\ref{ROGW-d-inv}), we suppose that there
exists an anti-symplectic involution $\phi$ on $X$, $i.e.$
$\phi^*\omega=-\omega$, such that $L={\rm Fix}(\phi)$ and the
restriction of $\phi$ to $V$ is also an involution on the
submanifold. It is not difficult to find an example for our setting.
For example, let $(X,L)=(\C P^1,\R P^1)$, the symplectic form is
induced by the Fubini-Study metric, the anti-symplectic involution
$\phi$ is the complex conjugation, and $V=\{x^+,x^-\}$ such that
$x^-\neq x^+$ and $x^-=\phi(x^+)$.

Furthermore, suppose $\Sigma$ is biholomorphic to its conjugation
$\bar{\Sigma}$, $i.e.$ there exists an anti-holomorphic involution
$c:\Sigma\rightarrow\Sigma$. Denote by ${\mathcal J_\omega}$ the set of
$\omega$-tamed almost complex structures on $X$. Let $\mathcal P$ denote
the set of $J$-anti-linear inhomogeneous perturbation terms. Define
$${\mathcal J}_{\omega,\phi}:=\{J\in{\mathcal J}_{\omega}|\ \phi^*J=-J\}.$$
Now for each $J\in{\mathcal J}_{\omega,\phi}$, we denote by ${\mathcal
P}^J_{\phi,c}$ the set of $\nu\in{\mathcal P}$ satisfying
$d\phi\circ\nu\circ dc=\nu$. Denote
$$\mathbb{J}:=\{(J,\nu)|J\in{\mathcal J}_{\omega},\nu\in{\mathcal P}\},$$
$$\mathbb{J}_\phi:=\{(J,\nu)|J\in{\mathcal J}_{\omega,\phi},\nu\in{\mathcal P}^J_{\phi,c}\}\subset \mathbb{J}.$$

In Section 3, we will impose the $V$-compatibility conditions  on the pair
$(J,\nu)$, and denote by  $\mathbb{J}^V$ (resp. $\mathbb{J}^V_\phi$) the set of all $V$-compatible pairs $(J,\nu)\in\mathbb{J}$ (resp. $(J,\nu)\in\mathbb{J}_\phi$).
Fix $(J,\nu)\in\mathbb{J}^V_\phi$. Thus, from  the $V$-stable
$(J,\nu)$-holomorphic map $u:(\Sigma,\partial\Sigma)\mapsto(X,L)$ we
can define its conjugate $V$-stable $(J,\nu)$-holomorphic  map
$\tilde{u}=\phi\circ u\circ c$ representing the homology class
$[\tilde{u}]=-\phi_*d$, simply denoted by $\tilde{d}$. Denote
$\tilde{\bf d}=(\tilde{d},d_1,\cdots,d_m)$. So we have an induced
map
$$\phi':\overline{\mathcal
M}^{V,\mathbbm{s}}_{{\bf k},l,\mathbbm{l}}(\Sigma,X,L,{\bf d})\rightarrow
\overline{\mathcal M}^{V,\mathbbm{s}}_{{\bf
k},l,\mathbbm{l}}(X,L,\tilde{\bf d})$$ given by
$${\bf u}=(u,\vec{z},\vec{w},\vec{q})\mapsto\tilde{\bf u}=
(\tilde{u},\ (c|_{\partial\Sigma})^{|{\bf k}|}(\vec{z}),\
c^l(\vec{w}),\ c^\mathbbm{l}(\vec{q})).$$

Define
$$\Omega_\phi^{*}(X):=\{\gamma\in\Omega^{*}(X)|\ \phi^*\gamma=\gamma\},$$
$$\Omega_\phi^{*}(V):=\{\eta\in\Omega^{*}(V)|\ \phi^*\eta=\eta\}.$$
Now in the integral (\ref{ROGW-d-inv}) we take the forms
$\gamma_j\in\Omega_\phi^{2n}(X)$ and
$\eta_\jmath\in\Omega_\phi^{2n-2}(V)$.

Then we study the lower dimensional cases, $i.e.$ $n\le 3$. If
$d=\tilde{d}$, then using the method in this paper (see Sections
\ref{SEC-9}--\ref{SEC-12}), we can show that the integrals
(\ref{ROGW-d-inv}) are invariants of the tuple $(X,\omega, V,\phi)$.
Roughly speaking,  we can show that the map $\phi'$ is an induced
involution on the $V$-stable moduli space $\overline{\mathcal M}^{V,\
\mathbbm{s}}_{{\bf k},l,\mathbbm{l}}(\Sigma,X,L,{\bf d})$. And using the
isomorphism (\ref{thm-1.1-isom}), we can assign a sign to each
$V$-stable map (see Section \ref{SEC-10}). Then we can prove that on
each codimension 1 boundary stratum, the induced involution $\phi'$
is fixed point free and orientation reversing, $i.e.$ changing the
sign of each disc-bubbling stable map. Therefore, the contributions
from the codimension 1 boundary of $\overline{\mathcal M}^{V,\
\mathbbm{s}}_{{\bf k},l,\mathbbm{l}}(\Sigma,X,L,{\bf d})$ to the integrals
(\ref{ROGW-d-inv}) are eliminated. That implies the invariance.

However, for general situation, the relative homology class $d$
might not be $\phi$-anti-invariant, then $\phi'$ might not be a map
from the space of $V$-stable maps to itself. Thus the numbers
$\mathscr{R}\mathscr{N}$ are not always well-defined invariants. So
we have to modify the definition above.

The idea is to take all related moduli space together to eliminate
the contributions from the codimension 1 boundaries to the integral.
In fact, such boundary-canceling method has been used in \cite{C}.

For a homology class $\alpha\in H_2(X)$ and a relative class
$\beta\in H_2(X,L)$, we say $\beta_\C=\mathbbm{d}$ if there is a
holomorphic map $u$ of class $\beta$ such that its complex double
$u_\C$ represents the homology class $\alpha$. In this sense, we may
say $\alpha=\beta_\C$ is the doubling of $\beta$. Note that
$\beta_\C=(\tilde{\beta})_\C$. Denote by $\mathbbm{d}=d_\C$ the
homology class in $H_2(X)$ which is the doubling of $d$, and by
$\bar{\alpha}=(\alpha_{ai})$, $\bar{\gamma}=(\gamma_{j})$,
$\bar{\eta}=(\eta_{\jmath})$. Then we define
\begin{eqnarray}\label{ROGW-int}
&&\mathcal {I}:=\mathcal {I}^{\ V,\
\mathbbm{s}}_{X,\phi,g,\mathbbm{d},{\bf
k},l}(\bar{\alpha},\bar{\gamma},\bar{\eta})\nonumber\\
&&\ \ \ \ =\sum_{\begin{array}{c}\forall
\beta:\beta_\C=\mathbbm{d}\end{array}}\ \ \sum_{\begin{array}{c}
\gamma_j:[\gamma_j]={\rm
PD}(\xi_j),\eta_\jmath:[\eta_\jmath]={\rm PD}(\lambda_\jmath), \\
\forall(\vec{x},\vec{\xi},\vec{\lambda})\in {\mathcal R}\end{array}
}\nonumber\\
&&\ \ \ \ \ \ \ \ \int_{\overline{\mathcal M}^{V,\mathbbm{s}}_{{\bf
k},l,\mathbbm{l}}(X,L,\bar{\beta})
}\bigwedge_{a,i}evb_{ai}^*(\alpha_{ai}) \bigwedge_{j}
evi^*_j(\gamma_j)\bigwedge_{\jmath} evi^{I*}_\jmath(\eta_\jmath)\ ,\
\ \ \ \ \ \ \
\end{eqnarray}
where $\mathcal R$ is the space of real configurations (see (\ref{11-real-confg}) for detailed definition),
$$\bar{\beta}=(\beta,d_1,\cdots,d_m)\in H_2(X,L)\oplus H_1(L)^{\oplus
m},$$ PD$(\cdot)$ is the Poincar\'{e} dual of the point class. Now,
we state  our main result.

\begin{Theorem}\label{thm-1.2}
Fix a conformal structure $j$ on a bordered Riemann surface $\Sigma$.
Assume that $L\cap V=\emptyset$, and that ${\rm dim}X\le 6$, if $L$
is not orientable, we assume that ${\rm dim}X\le 4$ and $k_a\cong
w_1(d_a)+1 \mod 2$. If $(X,L)$ is semi-positive or
each moduli space $\overline{\mathcal M}^{V,\mathbbm{s}}_{{\bf
k},l,\mathbbm{l}}(\Sigma,X,L,\bar{\beta})$ generically contains no
$\phi$-multiply covered pseudoholomorphic map\footnote{A nonconstant
map $u:(\Sigma,\partial\Sigma)\rightarrow(X,L)$ is called
$\phi$-{\it multiply covered} or {\it not $\phi$-somewhere
injective} if there exists a regular point $z\in\Sigma$ such that
either $u(z)\in u(\Sigma\setminus\{z\})$ or $u(z)\in{\rm
Im}(\phi\circ u)$.}, then the integers $\mathcal {I}=\mathcal {I}^{\
V,\mathbbm{s}}_{X,\phi,g,\mathbbm{d},{\bf k},l}$ are independent
of the generic choice of pair $(J,\nu)\in\mathbb{J}^V_\phi$, the
choice of conformal structure $j$, or the choice of forms
$\alpha_{ai}\in\Omega^n(L,{\rm det}(TL))$,
$\gamma_j\in\Omega_\phi^{2n}(X)$,
$\eta_\jmath\in\Omega_\phi^{2n-2}(V)$. Therefore, $\mathcal {I}$ are
invariants of the tuple $(X,\omega,V,\phi)$.
\end{Theorem}
\noindent {\it Remark}. In particular, if $d=\tilde{d}$,
$\mathscr{R}\mathscr{N}:=\mathscr{R}\mathscr{N}(V,{\bf d},{\bf
k},l,\mathbbm{l})$ are invariants of the tuple $(X,\omega,V,\phi)$.
$\mathscr{R}\mathscr{N}$ can be regarded as an extension of
Solomon's definition of open GW invariants $N_{\Sigma,{\bf d},{\bf
k},l}$ to the relative case. After the article was posted on arxiv, the author was informed by Welschinger that a definition of relatively open invariants for some special cases, $i.e.$ cotangent
bundles of 2-sphere and real projective plane, has appeared in \cite{W4}.

\smallskip

\noindent {\it Remark}. In the proof of Theorem \ref{thm-1.1} and \ref{thm-1.2}, we use the assumption that (generically) there is no
($\phi$-)multiply covered map. Actually, we
believe, by using the virtual cycle techniques to the case of open
maps or  the equivariant Kuranishi structure (see, for example, section 7 of \cite{So}) or the expected polyfold techniques being developed by Hofer {\it et al}, that neither such assumption nor the semi-positivity assumption is
necessary and the Theorem  \ref{thm-1.1} and \ref{thm-1.2} holds for general
case. \\

In Section \ref{SEC-9}, we deal with the problem of orientation and
derive the conclusion of Theorem \ref{thm-1.1}. Then in Section
\ref{SEC-10}, we assign a sign to each $V$-stable map and study how
the induced action by involution changes the sign. The last two
sections devote to the proof of Theorem \ref{thm-1.2}. The next
Sections \ref{SEC-2} is preparation for
succeeding discussion,  we show the
definition of {\rm stable open $(J,\nu)$-holomorphic maps} and
describe the moduli space of such maps.  For the convenience of the reader, in \ref{SEC-APP} we
review some definitions and important conclusions in \cite{So} about
the orientation of determinant of real linear Cauchy-Riemann
operator.\\

\noindent{\bf Acknowledgments.}\ 
The author thanks the Institute of Mathematical Science at Nanjing University, the Department of
Mathematics at Princeton University, the IH\'{E}S and the BICMR for hospitality and support, part
of this work was done when he was working at or visiting these institutions. He also thanks Gang Tian for inspiring conversations on
relative invariants,  An-Min Li, Thomas Parker, Jack Solomon  and Jean-Yves Welschinger for explaining to him some points in
the definitions of relative GW invariants, open GW invariants and real enumerative invariants, respectively.
Thanks also go to Bohui Chen, Xiaojun Chen, Huijun Fan, Jianxun Hu, Guangcun Lu and Shanzhong Sun
for helpful conversations, and Scott Kenney for providing direct help for his working and living in Princeton.
He is very grateful to Yiming Long, Jiangong You for their long-term encouragement,
and Haihong Lin for special financial support in the spring of 2008 when he had no income.
He wants to thank  referees for their constructive comments and giving him many nice suggestions.

\section {Stable open $(J,\nu)$-holomorphic maps}\label{SEC-2}

In this section, we recall the definition of {\rm stable open
$(J,\nu)$-holomorphic maps}, and describe the moduli space of such
maps. The domains of such stable open maps are bordered Riemann
surface $i.e.$ compact Riemann surface with boundary, smooth or
allowing nodal singularities. The boundary is mapped in a Lagrangian
submanifold.

In the following, both marked points and double points (or singular
points or nodes) are called {\it special points}. And we always fix
a conformal structure on the open domain curve $\Sigma$, which is a
bordered compact Riemann surface. We say such curve is of genus $g$
if the closed Riemann surface
$\Sigma\cup_{\partial\Sigma}\bar{\Sigma}$, which is the complex
double of $\Sigma$ and also denoted by $\Sigma_\C$ (see \cite{AG}),
is of genus $g$. Denote the genus of the closed surface
$\Sigma/\partial\Sigma$ by $g_0$.

Assume that there are altogether $m$ boundary components, $i.e.$
$\partial\Sigma=\bigcup_{a=1}^m(\partial\Sigma)_a$. We say $\Sigma$
is of topological type $(g_0,m)$.
\begin{Definition}\label{aut}
An automorphism of a bordered Riemann surface $\Sigma$ is a
diffeomorphism $\varphi:\Sigma\rightarrow\Sigma$ preserving the
conformal structure and the ordering of the boundary components. The
set of all automorphism of $\Sigma$ is denoted by ${\rm
Aut}(\Sigma)$. We say $\Sigma$ is stable if ${\rm Aut}(\Sigma)$ is
finite.
\end{Definition}
Denote by $\Sigma_\C$ the complex double of $\Sigma$, which is a
closed Riemann surface with an antiholomorphic involution $\sigma$.
If $\varphi:\Sigma\rightarrow\Sigma$ is an automorphisim, then its
complex double $\varphi_\C:\Sigma_\C\rightarrow\Sigma_\C$ is an
automorphism of $(\Sigma_\C,\sigma)$. This provides a natural
inclusion ${\rm Aut}(\Sigma)\subset{\rm Aut}(\Sigma_\C,\sigma)$. The
following statements are equivalent:\\

$\bullet$ $\Sigma$ is stable, $i.e.$, ${\rm Aut}(\Sigma)$ is finite.

$\bullet$ $\Sigma_\C$ is stable.

$\bullet$ The genus $g=2g_0+m-1$ of $\Sigma_\C$ is bigger than one.

$\bullet$ The Euler characteristic $\chi(\Sigma)=2-2g_0-m$ of
$\Sigma$ is negative.

\begin{Definition}\label{prestable}
A prestable bordered Riemann surface with fixed conformal structure
is either a smooth bordered Riemann surface or the union of a smooth
bordered Riemann surface with finite sphere and disc bubbles.
\end{Definition}
In this paper, we only consider the prestable bordered Riemann
surfaces with fixed conformal structure, and simply call them
prestable bordered Riemann surfaces.

\begin{Definition}\label{aut-mark}
An automorphism of a prestable bordered Riemann surface $\Sigma$ of
topological type $(g_0,m)$ with $({\bf k},l)$ marked points
$$(\Sigma,\partial\Sigma,\vec{z},\vec{w})=
(\Sigma,\{(\partial\Sigma)_a\}_{a=1}^m,\{z_{a1},\cdots,z_{ak_a}\}_{a=1}^m,w_1,\cdots,w_l)$$
is an automorphism of $\Sigma$ preserving all marked points. The set
of all automorphism is denoted by ${\rm Aut}(\Sigma_{{\bf k},l})$ (or simply by ${\rm Aut}(\Sigma)$). A
prestable $({\bf k},l)$-marked bordered Riemann surface $\Sigma$ is
stable if ${\rm Aut}(\Sigma_{{\bf k},l})$ is finite.
\end{Definition}

\begin{Definition}\label{2-type}
A $({\bf k},l)$-marked bordered bubble domain curve $\Sigma$ of type
$(g_0,m)$, which is also called a prestable $({\bf k},l)$-marked
bordered Riemann surface and where the vector ${\bf
k}=(k_1,\cdots,k_m)$, is a finite connected union of smooth oriented
compact surfaces $\Sigma_i$, at least one surface with boundary,
joined at interior or boundary double points together with $k_{a}$
distinct marked points in the boundary $(\partial\Sigma)_a$,
$a=1,\cdots,m$, and $l$ distinct interior marked points, none of
which are double points. The $\Sigma_i$, with their special points,
are of two types:

{\rm (1)} stable components, and

{\rm (2)} unstable components, which are unstable sphere bubbles or
unstable disc bubbles.

\noindent And there must be at least one stable component.
\end{Definition}

There exists a natural stablization map
\begin{equation}\label{stable}
{\rm st}: \Sigma\rightarrow\widehat{\Sigma}
\end{equation}
that collapses the unstable components to points, thus we get a
connected domain $\widehat{\Sigma}={\rm st}(\Sigma)$ which is a
stable genus $g$ open curve.

Bordered bubble domains can be constructed from a stable bordered
Riemann surface ${\Sigma}_0$ by replacing points by finite chains of
2-spheres or 2-discs or their combination. Alternatively, they can
be obtained from a smooth Riemann surface $\Sigma_0$ by pinching a
set of nonintersecting embedded circles in the interior of
${\Sigma}_0$ and (or) a set of half-circles in ${\Sigma}_0$ with
centers in the boundary $\partial\Sigma_0$. The latter viewpoint can
be formalized as follows. Assume that there are $b_a$ double points
on each  boundary component $(\partial\Sigma)_a$, $a=1,\cdots,m$,
and there are $d$ interior double points. Denote by ${\bf
b}=(b_1,\cdots,b_m)$.

\begin{Definition}\label{resl}
A resolution of $(g_0,m)$-type $({\bf k},l)$-marked bordered bubble
domain $\Sigma$ with $({\bf b},d)$-double points is a smooth
bordered Riemann surface with genus $g$, $d$  embedded circles
$\gamma_c$ in the interior part and $b_a$  embedded half-circles
$\gamma_{ah}$ with centers in each $(\partial\Sigma_0)_a$,
$a=1,\cdots,m$, (any two of distinct circles or half-circles are
disjoint), and the $({\bf k},l)$-marked points are apart from
$\gamma_{ah}$ and $\gamma_{c}$, together with a resolution map
$$\mathscr{R}:\Sigma_0\rightarrow\Sigma$$
which respects orientation and marked points, takes each $\gamma_c$
(resp. $\gamma_{ah}$) to an interior (resp. boundary) double point
of $\Sigma$, and restricts to a diffeomorphism from the complement
of $\bigcup\bar{\gamma}_{ah}\cup\bar{\gamma}_c$ in $\Sigma$ to the
complement of the double points, where $\bar{\gamma}_{ah}$ (resp.
$\bar{\gamma}_c$) denotes the closure of half-disc (resp. disc)
contained by $\gamma_{ah}$ (resp. $\gamma_c$).
\end{Definition}

We next define $(J,\nu)$-holomorphic maps from bubble domains. Such
maps depend on the choice of an $\omega$-compatible (tamed) almost
complex structure $J\in{\mathcal J}_\omega$ and an inhomogeneous perturbation $\nu$ to the
Cauchy-Riemann equation. Let $\mathcal B$ be a parameter space which will be specified later. Let $\pi_i$, $i=1,2$, denote the projection
from $\Sigma\times X\times\mathcal{B}$ to the $i^{th}$ factor and let $\pi'_i$ denote
the restriction of $\pi_i$ to $\partial\Sigma\times L\times\mathcal{B}$. We define
the inhomogeneous term to be the section
\begin{equation}\label{2-def-nu}
\nu\in\Gamma(\Sigma\times X\times\mathcal{B},{\rm Hom}(\pi_1^*T\Sigma,\pi_2^*TX))
\end{equation}
such that

(1) $\nu$ is $(j_\Sigma,J)$-anti-linear:\ \  $\nu\circ
j_\Sigma=-J\circ\nu$;

(2) $\nu|_{\partial\Sigma\times L\times\mathcal{B}}$ carries a sub-bundle
${\pi'}_1^*T\partial\Sigma\subset \pi_1^*T\Sigma$ to the sub-bundle
${\pi'}_2^*(TL)\subset\pi_2^*TX$.\footnote{Note here our setting, i.e. the following equation (\ref{def-Jnu}), is a little different from the one in \cite{So}.}\\

Denote by $\mathcal P$ the set of all such inhomogeneous terms. Let
$\mathbb{J}$ denote the space of such pairs $(J,\nu)$. The parameter space $\mathcal{B}$ is
often taken to be the ambient space of relevant pseudo-holomorphic maps. For instance, in this section, let
$\mathcal{B}=B^{1,p}(\Sigma,L,{\bf d})$ be the Banach manifold of $W^{1,p}$ maps
$u:(\Sigma,\partial\Sigma)\rightarrow(X,L)$ such that
$u_*([\Sigma,\partial\Sigma])=d$ and
$u|_{(\partial\Sigma)_a*}([(\partial\Sigma)_a])=d_a$. In section 3, we will take $\mathcal{B}$ as the space (\ref{B-kl-2}).

\begin{Definition}\label{def-3}
A $(J,\nu)$-holomorphic open map from a bordered bubble domain curve
$(\Sigma,\partial\Sigma)$ with complex structure $j_\Sigma$ is a map
$$u:(\Sigma,\partial\Sigma)\rightarrow (X,L)$$ such that, on each
component $\Sigma_i$ of $\Sigma$, $u$ is a solution of the
inhomogeneous Cauchy-Riemann equation
\begin{equation}\label{def-Jnu}
\bar{\partial}_J u=\frac{1}{2}(du+J\circ du\circ
j_\Sigma)=\nu(\cdot,u(\cdot),\bf{u}),
\end{equation}
or equivalently, $$\bar{\partial}_{(J,\nu)} u=0$$ where
$\bar{\partial}_{(J,\nu)}$ denotes the perturbed nonlinear elliptic
operator $\frac{1}{2}(d+J\circ d\circ j_\Sigma)-\nu$. In particular,
$\bar{\partial}_J u=0$ on each unstable bubble.
\end{Definition}

\noindent{\it Remark}. It is not difficult to show that the operator
$\bar{\partial}_{(J,\nu)}$ gives rise to an elliptic boundary value
problem. We refer to the Lemma 4.1 in \cite{So}. \\

The symplectic area  of the image is the number
\begin{equation}\label{area}
A(u)=\int_{u(\Sigma)}\omega=\int_\Sigma u^*\omega
\end{equation}
which depends only on the homology class of the curve modulo its
boundary. And the energy of $u$ is
\begin{equation}\label{en}
E(u)=\frac{1}{2}\int_{\Sigma}|du|^2_{J,\mu}d\mu
\end{equation}
where $|\cdot|_{J,\mu}$ is the norm defined by the metric on $X$
determined by $J$ and the metric $\mu$ on the domain. For
$(J,0)$-holomorphic maps, $E(u)=A(u)$.

We also define a {\it modified energy} $\mathbb{E}$ componentwise
\begin{equation}\label{E-new}
\mathbb{E}(u_i)= \left\{
\begin{array}{l}
1+\frac{1}{2}\int_{\Sigma_i}|du_i|^2_{J,\mu}d\mu,\ \ \ \ \ \ \ \ \Sigma_i\ {\rm is\ stable\ and\ } u_i\  {\rm is\ not\ a\ ghost  ;}\\
\frac{1}{2}\int_{\Sigma_i}|du_i|^2_{J,\mu}d\mu,\ \ \ \ \ \ \ \ \ \ \
\ \ \Sigma_i\ {\rm is\ unstable.}
\end{array}
\right.
\end{equation}And
\begin{equation}\label{E-mod}
\mathbb{E}(u)=\sum_i\mathbb{E}(u_i).
\end{equation}

Now we come to the definition of {\it stable} map
\begin{Definition}\label{def-stable}
A $(J,\nu)$-holomorphic map $u$ is stable if each of its component
maps $u_i=u|_{\Sigma_i}$ has positive modified energy $i.e.$
$\mathbb{E}(u_i)>0$ for each $i$.
\end{Definition}
That is to say, either each component $\Sigma_i$ of the domain is   stable, or else the image of $(\Sigma_i,\partial\Sigma_i)$
carries a nontrivial homology class. We have
\begin{Lemma}\label{L-abc}
Let $(X,\omega)$ be a compact symplectic manifold and  $L$ be a
compact Lagrangian submanifold. Then

{\rm (1)} every $(J,\nu)$-holomorphic map has $\mathbb{E}(u)\geq 1$.

{\rm (2)} There exists a constant $\hbar>0$ such that  for every
component $u_i$ of every stable $(J,\nu)$-holomorphic map
$u:(\Sigma,\partial\Sigma)\rightarrow(X,L)$ with Lagrangian boundary
condition, we have
\bq\label{Eh} \mathbb{E}(u_i)>\hbar \eq 
\end{Lemma}
Proof.\  (1) The conclusion is direct since at least one component
is stable.

(2) On the stable components, we have $\mathbb{E}(u_i)\geq 1$. On
the unstable components, since $u_i$ is a $J$-holomorphic map, the
Proposition 4.1.4 of \cite{MS} implies there exists $\hbar>0$
depending only on $(X,J)$ such that
$$\mathbb{E}(u_i)=E(u_i)=\frac{1}{2}\int_{\Sigma_i}|du_i|^2_{J,\mu}d\mu>\hbar,$$
provided the images of $(\Sigma_i,\partial\Sigma_i)$ carries a
nontrivial homology class. Otherwise, if $u_i$ represents a trivial
homology class, then $\mathbb{E}(u_i)=E(u_i)=A(u_i)=\omega\cap
[u_i(\Sigma_i)]=0$, contrary to the definition of stable map. \qed

The following Gromov Convergence Theorem, which is also called
Compactness Theorem, is the important fact about
$(J,\nu)$-holomorphic maps. It means that every sequence of
$(J,\nu)$-holomorphic maps from a smooth (bordered) domain has a
subsequence which converges modulo automorphism to a stable map.
Various forms of such convergence theorem have been proved for
closed curves (c.f. \cite{G}, \cite{PW}, \cite{RT1}, \cite{Y},
$etc.$). For open curves, we refer to  \cite{FO3},\cite{Liu}, \cite{MS}. Moreover, we assume
$|\bf{k}|>0$ to exclude some exceptional case\footnote{When $|\bf{k}|=0$, it turns out that the standard
stable map moduli space might be not compact. See sections 3.8 and 7.4 of \cite{FO3} for discussion of such case.}.

\begin{Theorem}\label{thm-bubble}{\rm (Bubble Convergence)}.\ \
Let $(X,\omega)$ be a compact symplectic manifold, $L\subset X$ be a
Lagrangian submanifold. Given any sequence $\{u_j\}$ of $({\bf
k},l)$-marked $(J_j,\nu_j)$-holomorphic maps from a bordered Riemann
surface $\Sigma_0$ satisfying Lagrangian boundary conditions, with
$\mathbb{E}(u_j)<E_0$ and $(J_j,\nu_j)\rightarrow(J,\nu)$ in $C^k$,
$k\geq 0$, then we can obtain a subsequence and

(1) a $({\bf k},l)$-marked bordered bubble domain $\Sigma$ with
resolution $\mathscr{R}:\Sigma_0\rightarrow\Sigma$, and

(2) automorphisms $\varphi_j$ of $\Sigma_0$ preserving the
orientation and the marked points,\\
such that the modified subsequence $\{u_j\circ\varphi_j\}$ converges
to a limit
$$\CD
  \Sigma_0 @>\mathscr{R}>> \Sigma @>u>> X
\endCD$$
where $u$ is a stable $(J,\nu)$-holomorphic map. This convergence is
in $C^0$, and in $C^k$ on compact sets not intersecting the
collapsing curves $\gamma_{ah}$ and $\gamma_c$ of the resolution
$\mathscr{R}$, and the area (\ref{area}) and energy (\ref{E-mod})
are preserved in the limit.
\end{Theorem}

Given a bordered bubble domain curve with fixed complex structure
$(\Sigma,\partial\Sigma)$. Denote the space of equivalence classes
of stable $({\bf k},l)$-marked open pseudoholomorphic maps with
Lagrangian boundary conditions representing the homology class ${\bf
d}$ by
$\overline{\mathcal M}_{{\bf k},l}(\Sigma,X,L,{\bf d})$ which is just the
Gromov compactification of ${\mathcal M}_{{\bf k},l}(\Sigma,X,L,{\bf d})$.

\begin{Definition}\label{multi-cover}
We say a map $u:(\Sigma,\partial\Sigma)\rightarrow(X,L)$ is {\rm
multiply covered} if there does not exist a point $z\in\Sigma$ such
that
$$du(z)\ne 0,\ \ \ {\rm and}\ \ \ u(z)\notin\ u(\Sigma\setminus\{z\}).$$
A  map is called {\rm irreducible} or {\rm somewhere injective} if
it is not multiply covered. The moduli space ${\mathcal M}_{{\bf
k},l}(\Sigma,X,L,{\bf d})$ or $\overline{\mathcal M}_{{\bf k},l}(\Sigma,X,L,{\bf d})$
is called irreducible if it contains no multiply covered
pseudoholomorphic map.
\end{Definition}
\noindent{\it Remark}. As mentioned in the introduction, if there exists an
anti-symplectic involution $\phi$ on the symplectic manifold $X$
such that $L={\rm Fix}(\phi)$, there is a modified definition of the
so-called $\phi$-{\it multiply covered map}, we refer to the
footnote of the Theorem \ref{thm-1.2}. The multiply covered maps are often singular points in the moduli
space ${\mathcal M}_{{\bf k},l}(\Sigma,X,L,{\bf d})$. For keeping the paper in
suitable length, we will avoid dealing with these maps
in the present paper.

Let $\overline{\mathcal M}_{{\bf k},l}(\Sigma,X,L,{\bf d})^*$ be the space of
irreducible stable open $(J,\nu)$-maps. The following theorem is direct from
Theorem \ref{thm-bubble} , the index theorem and the standard arguments of transversality theorem
\begin{Theorem}\label{stablemap}
Fix a conformal structure $j$ on $\Sigma$ (forgetting $({\bf k},l)$ marked points), denote the group of automorphisms of $\Sigma$ by $Aut(\Sigma)$. We assume $|\bf{k}|>0$, then

(1) $\overline{\mathcal M}_{{\bf k},l}(\Sigma,X,L,{\bf d})$ is compact and
Hausdorff, and there exists a continuous evaluation map
\begin{equation}\label{EV1}
{\rm ev}: \overline{\mathcal M}_{{\bf k},l}(\Sigma,X,L,{\bf d})\longrightarrow
X^l\times L^{|{\bf k}|}
\end{equation}

(2) For generic $(J,\nu)\in\mathbb{J}$, $\overline{\mathcal M}_{{\bf
k},l}(\Sigma,X,L,{\bf d})^*$ is a manifold with boundary (corners) of
dimension
\begin{equation}\label{1-dim}
{\rm dim}\ \overline{\mathcal M}_{{\bf k},l}(\Sigma,X,L,{\bf d})^*=\mu(d)+n(1-g)+k+2l-\dim Aut(\Sigma).
\end{equation}
\end{Theorem}
\noindent {\it Remark}. The conclusions of the theorem above did not appear exactly in the literature, while there exist  related
results in some different settings. In
\cite{FO3} and \cite{Liu}, they show that, for general symplectic manifold, the stable moduli space of open unperturbed $J$-maps is compact and Hausdorff and carries the smooth Kuranishi structure, from which they respectively constructed the virtual fundamental chain (see Theorem 2.1.29.-32. of \cite{FO3}, or sections 5.3 and 7.1 of \cite{Liu}). In  \cite{MS}, they only deal with the genus 0 unperturbed $J$-maps (in fact, they claim that their arguments go through to higher genus case with fixed complex structure), and conclude that
 the irreducible part of the moduli space is generically a manifold. In \cite{RT1}\cite{RT2} and \cite{LiT}, they study the closed $(J,\nu)$-maps
 in  semi-positive and general symplectic manifolds, respectively. Their arguments may apply to our open case with minor modifications, since here we just describe the stable moduli space instead of defining invariants. We will not give a proof of the theorem above here, and refer to literature mentioned above and a similar arguments, corresponding to Theorem \ref{stablemap} (2), in the proof of Lemma \ref{dim-M}.

\section{$V$-compatible pair $(J,\nu)$}\label{SEC-4}

We now come to extend Solomon's open symplectic invariants
$N_{\Sigma,{\bf d},{\bf k},l}$\footnote{Also one can define the Solomon-type open invariant for more general case, for one such generalization see \cite{H}.}  to open invariants of $(X,\omega,\phi)$ relative to a
codimension 2 symplectic submanifold $V$ such that $\phi|_V$ is also
an anti-symplectic involution on $(V,\omega|_V)$. Recall the
assumption that $L=$ Fix$(\phi)\neq\emptyset$ is a Lagrangian
submanifold. Note that $V\cap L$ might be empty set. If $V\cap
L={\rm Fix}(\phi|_V)\neq\emptyset$, then we denote by $L_V=V\cap L$
which is a Lagrangian submanifold of $(V,\omega|_V)$. We consider
open curves generically intersect $V$ in a finite collection of
points. Such relative open invariants will count these open curves
satisfying some constraints. In particular, if $L\cap V=\emptyset$,
we will not encounter extra codimension 1 boundary of moduli space
except the moduli space of pseudoholomorphic maps with a bubble
disc. If $\dim L\le 3$, under some assumptions, we will define
relatively open invariants for domain curves of any genus with fixed
conformal structures.

Let us reset our notations for new discussion. Suppose ${\rm
dim}X=2n$. We denote still by $(\Sigma,\partial\Sigma)$ the domain
Riemann surface with $m$ boundary components and with fixed
conformal structure $j_\Sigma$. That means we also don't deal with
the case that the degenerations of $\Sigma$ may occur. Denote by $g$
the genus of the closed Riemann surface
$\Sigma\cup_{\partial\Sigma}\bar{\Sigma}$ obtained by doubling
$\Sigma$. If $\Sigma$ is a closed Riemann surface then $g$ is just
the genus of itself. Assume that there are $k_a$ distinct marked
points $\{z_{a1},\cdots,z_{ak_a}\}$ on the boundary component
$(\partial\Sigma)_a$, $a=1,\cdots,m$, and $l$ distinct interior
marked points $\{w_1,\cdots,w_l\}$. Denote by ${\bf
k}=({k_1,\cdots,k_m})$. Additionally, we assume that there are
$\mathbbm{k}_a$ marked points $\{p_{a1},\cdots,p_{a\mathbbm{k}_a}\}$
on each boundary component $(\partial\Sigma)_a$, $a=1,\cdots,m$,
which are different from $\{z_{a1},\cdots,z_{ak_a}\}$ and are mapped
to the intersection points of $V$ and the image of our open curve
$u:(\Sigma,\partial\Sigma)\rightarrow (X,L)$ if $V\cap
L\neq\emptyset$. Denote by
$\mathbbm{k}=(\mathbbm{k}_{1},\cdots,\mathbbm{k}_m)$. Also, we
assume there are $\mathbbm{l}$ interior marked points
$\{q_{1},\cdots,q_{\mathbbm{l}}\}$  which are different from
$\{w_1,\cdots,w_l\}$ and are mapped to the intersection points of
$V$ and the image of $u$. For given homology class ${\bf
d}=(d,d_1,\cdots,d_m)\in H_2(X,L)\oplus H_1(L)^{\oplus m}$, and
given a pair $(J,\nu)$, denoted by ${\mathcal M}_{{\bf
k},l,\mathbbm{l}}(\Sigma,X,L,{\bf d})$ the moduli space of
$({\bf k},l,\mathbbm{l})$-marked $(J,\nu)$-pseudoholomorphic maps such that $u_*([\Sigma,\partial\Sigma])=d$
and $u|_{(\partial\Sigma)_a*}([(\partial\Sigma)_a])=d_a$, $i.e.$
representing ${\bf d}$, which is the zero set of
\begin{equation}\label{J-eqn}
\Phi(u,J,\nu)=\bar{\partial}_{(J,\nu)}\
u=\bar{\partial}_Ju-\nu=\frac{1}{2}(du+J\circ du\circ j_\Sigma)-\nu.
\end{equation}
We denote the ambient space, which is the Sobolev completion  of
${\mathcal M}_{{\bf k},l,\mathbbm{l}}(\Sigma,X,L,{\bf d})$, by
\begin{equation}\label{B-kl-2}
\mathcal{B}=B_{{\bf k},l,\mathbbm{l}}^{1,p}(\Sigma,X,L,{\bf d}):=B^{1,p}(\Sigma,X,L,{\bf
d})\times\prod_a(\partial\Sigma)_a^{(k_a+\mathbbm{k}_a)}\times\Sigma^{(l+\mathbbm{l})}\backslash
\triangle.
\end{equation}  If we use vectors
$\vec{z}=(z_{ai}),\vec{w}=(w_j),\vec{p}=(p_{a\imath}),\vec{q}=(q_{\jmath})$
to denote marked points respectively, then we denote elements of
$B_{{\bf k},l,\mathbbm{l}}^{1,p}(\Sigma,X,L,{\bf d})$ by ${\bf
u}=(u,\vec{z},\vec{w},\vec{p},\vec{q})$. We simply denote the
restricted pullback bundle $(u|_{\partial\Sigma})^*TL$ by $u^*TL$.
By straightforward computations we have the following proposition.

\begin{Proposition}\label{Du}
The linearization of (\ref{J-eqn}) at ${\bf u}\in{\mathcal M}_{{\bf
k},l,\mathbbm{l}}(\Sigma,X,L,{\bf d})$ is
$$D_{\bf u}:=D_{\bf u}\bar{\partial}_{(J,\nu)}:\Gamma[(\Sigma,\partial\Sigma),(u^*TX,u^*TL)]
\rightarrow\Omega^{0,1}(\Sigma,u^*TX)$$
\begin{eqnarray}\label{Du1}
D_{\bf u}(\xi)&=&\frac{1}{2}[\nabla\xi+J\circ\nabla\xi\circ j+\nabla_\xi J\circ du\circ j]\nonumber\\
& &+\frac{1}{2}[T(\xi,du)+JT(\xi,du\circ j)]-\nabla_\xi\nu
\end{eqnarray}
where $\nabla$ is the pullbacked connection on $u^*TX$ and
$T(\zeta,\eta)=\nabla_\zeta\eta-\nabla_\eta\zeta-[\zeta,\eta]$ is
the torsion of $\nabla$.
\end{Proposition}
Denote by
\begin{equation}\label{Du2}
{\mathcal D}_{u}(\xi)=\frac{1}{2}[\nabla\xi+J\circ\nabla\xi\circ
j+\nabla_\xi J\circ du\circ j]-\nabla_\xi\nu.
\end{equation}

We can similarly define the Banach space bundle ${\mathcal E}\rightarrow
B_{{\bf k},l,\mathbbm{l}}^{1,p}(\Sigma,X,L,{\bf d})$ fiberwise
by
$${\mathcal E}_{\bf u}:=L^p(\Sigma,\Omega^{0,1}(u^*TX)).$$
So we still think $$\bar{\partial}_{(J,\nu)}:B_{{\bf
k},l,\mathbbm{l}}^{1,p}(\Sigma,X,L,{\bf d})\rightarrow{\mathcal
E}$$ as a section of $\mathcal E$. And we also denote by
$$D:=D\bar{\partial}_{(J,\nu)}:TB_{{\bf k},l,\mathbbm{l}}^{1,p}(X,L,{\bf
d})\rightarrow{\mathcal E}$$ the vertical component of the linearization
of $\bar{\partial}_{(J,\nu)}$. From the Proposition 6.14 in
\cite{Liu} we know that for each ${\bf u}\in B_{{\bf
k},l,\mathbbm{l}}^{1,p}(\Sigma, X,L,{\bf d})$, both $D_{\bf u}$
and ${\mathcal D}_u$ are Fredholm operators with the same index
ind$(D_{\bf u})={\rm ind}({\mathcal D}_u)=\mu+n(1-g)$, where
$\mu=\mu(u^*TX,u^*TL)$ is the total boundary Maslov index associated
with the vector bundle pair $(u^*TX,u^*TL)$ defined in the Appendix (\ref{eq-A1}).
In particular, if $u$ is a $(J,\nu)$-holomorphic map with closed
genus $g$ domain curve, then ind$(D_{\bf u})={\rm ind}({\mathcal
D}_u)=2c_1+2n(1-g)$, where $c_1$ is the first Chern number of the
bundle $u^*TX\rightarrow\Sigma$.

Also we suppose $\Sigma$ is biholomorphic to its conjugation
$\bar{\Sigma}$, $i.e.$ there exists an anti-holomorphic involution
$c:\Sigma\rightarrow\Sigma$. Now for each $J\in{\mathcal
J}_{\omega,\phi}$, we denote by ${\mathcal P}^J_{\phi,c}$ the set of
$\nu\in{\mathcal P}$ satisfying $d\phi\circ\nu\circ dc=\nu$. Recall
$$\mathbb{J}:=\{(J,\nu)|J\in{\mathcal J}_{\omega},\nu\in{\mathcal P}\},$$
and
$$\mathbb{J}_\phi:=\{(J,\nu)|J\in{\mathcal J}_{\omega,\phi},\nu\in{\mathcal P}^J_{\phi,c}\}\subset \mathbb{J}.$$

We still assume $L$ is relatively $Pin^\pm$ and fix a relative
$Pin^\pm$ structure $\mathfrak{P}$ on L. If $L$ is orientable, we
fix an orientation on $L$.

Denote the orthogonal projection onto the normal bundle $N_V$ by
$\xi\mapsto\xi^N$. Since $L_V$, if nonempty, is a submanifold of
$L$, we denote the normal bundle of $L_V$ in $L$ by $N_{L_V}$. Then
for each $(J,\nu)$-holomorphic open map $u$ whose image lies in
$(V,L_V)$, the operator
$$D_u^N: \ \Gamma[(\Sigma,\partial\Sigma),(u^*N_V,u^*N_{L_V})]\rightarrow\Omega^{0,1}(\Sigma,u^*N_V),$$
\begin{equation}\label{D-N}
D_{\bf u}^N(\xi)=[D_{\bf u}(\xi)]^N,
\end{equation}
restricting the vertical linearization of
$\bar{\partial}_{(J,\nu)}$ at $u$ to the normal bundle, is also a
real linear Fredholm operator.

Similar to the construction of relative GW invariants, we will
restrict attention to a subspace of $\mathbb{J}$
(or $\mathbb{J}_\phi$ involving the involution
) consisting of pairs $(J,\nu)$ that are compatible with $V$ in the
following sense. Recall  we  let $\pi_i$, $i=1,2$, denote the projection
from $\Sigma\times X$ to the $i^{th}$ factor and let $\pi_i'$ denote
the restriction of $\pi_i$ to $\partial\Sigma\times L$.

\begin{Definition}\label{3.1} We say the pair $(J,\nu)\in\mathbb{J}
$ or $\mathbb{J}_\phi$ is $V$-compatible if the
following three conditions hold:

(1) $J$ preserves $TV$,  $\nu^N|_V=0$, and if $L\cap V\neq\emptyset$, then $\nu|_{\partial\Sigma\times L_V\times\mathcal{B}}$ carries a sub-bundle
${\pi'}_1^*T\partial\Sigma\subset \pi_1^*T\Sigma$ to the sub-bundle
${\pi'}_2^*(TL_V)\subset\pi_2^*TX$; \\


and for all $\xi\in (N_V,N_{L_V})$, $v\in TV$ and $\vartheta\in T\Sigma$

(2) $[(\nabla_\xi J+J\nabla_{J\xi}J)(v)]^N=[(J\nabla_{Jv}J)\xi+(\nabla_{v}J)\xi]^N$;\\

(3) $[(J\nabla_{\nu(\vartheta)}
J)\xi]^N=[(\nabla_\xi\nu+J\nabla_{J\xi}\nu)(\vartheta)]^N$.
\end{Definition}

We denote by $\mathbb{J}^V$ (resp. $\mathbb{J}_\phi^V$) the set of
all $V$-compatible pairs $(J,\nu)\in\mathbb{J}$ (resp.
$(J,\nu)\in\mathbb{J}_\phi$).\\

\noindent{\it Remark}. The definitions above are  modification of the closed version of the Definition 3.2 of \cite{IP1}. However, in our open case, the condition (1), compared with the condition (a) of the Definition 3.2 of \cite{IP1}, involves the Lagrangian boundary condition (see (\ref{2-def-nu}) for the definition of inhomogeneous  perturbation $\nu$). It implies that $V$ is a $J$-holomorphic submanifold, and that $(J,\nu)$-holomorphic (open) curves in $V$ or $(V,L_V)$ are also $(J,\nu)$-holomorphic in $X$ or $(X,L)$, moreover, $\bar{\partial}_{(J,\nu)}$ gives rise to elliptic boundary value problems for both two kinds of open $(J,\nu)$-maps into $(X,L)$ and $(V,L_V)$. When $L\cap V=\emptyset$, conditions (2) and (3) are the same as the conditions (b) and (c) of the Definition 3.2 of \cite{IP1} (note $\xi\in(N_V,N_{L_V})$ means the section of normal bundle $N_V$ takes value in $N_{L_V}$ over $L_V$, if nonempty). In fact, they are mainly used to assure the following Lemma \ref{comp-op} holds  and the same conclusion of Lemma \ref{comp-op} is important for the discussion in \cite{IP1}, since they desire the boundary strata of codimension at least 2. \qed

Conditions (2) and (3) ensure that for each $(J,\nu)$-holomorphic
map with closed domain curve whose image lies in $V$
$$u:\Sigma\rightarrow V,$$ the operator
$${\mathcal D}_u^N: \ \Gamma[(\Sigma,u^*N_V)\rightarrow\Omega^{0,1}(\Sigma,u^*N_V),$$
\begin{equation}\label{D^N}
{\mathcal D}_{u}^N(\xi)=[{\mathcal D}_{u}(\xi)]^N
\end{equation}
by restricting the linearization of $\bar{\partial}_{(J,\nu)}$ at
$u$ to the normal bundle, is a complex linear operator, $i.e.$ we
have the following lemma which is taken from \cite{IP1}.

\begin{Lemma}\label{comp-op}
Choose $(J,\nu)\in \mathbb{J}^V$ (resp. $\mathbb{J}_\phi^V$). Then
for each $(J,\nu)$-holomorphic map $u$ from closed domain curve
whose image lies in $V$, the operator ${\mathcal D}_{u}^N$ is a complex linear operator.
\end{Lemma}

\bigskip

In the following we will not distinguish $D$ and $\mathcal {D}$ and
denote always by $D$ the linearization operator. In the end of the
section, we show a local normal form for open holomorphic maps near
the points where they intersect $V$. Here the argument is similar to
the Lemma 3.4 in \cite{IP1}.

Take a pair $(J,\nu)$ satisfying the condition (1) in the Definition
\ref{3.1}. Let $V$ be a codimension two $J$-holomorphic submanifold
of $X$, and $u:(\Sigma,\partial\Sigma)\rightarrow(X,L)$ be a
$(J,\nu)$-holomorphic map which intersects $V$. Suppose $L_V=V\cap
L\neq\emptyset$. Take a boundary marked point
$p_{a\imath}\in(\partial\Sigma)_a$ and an interior marked point
$q_\jmath$ satisfying $\mathbbm{p}=u(p_{a\imath})\in L_V$ and
$\mathbbm{q}=u(q_{\jmath})\in V$. Let $\H=\{z=x+{\rm i}y |\ y\ge
0\}$ be the upper half complex plane. Fix a local holomorphic
coordinate $z\in \H$ on a half open disc $\mathfrak{D}$ in $\Sigma$
containing $p_{a\imath}$ or a local holomorphic coordinate $z'\in\mathbb{C}$ on an
open disc $\mathcal {O}$ in $\Sigma$ containing $q_\jmath$. Also fix
local coordinates $\{v^i\}=\{v^{i_1},v^{i_2}\}$ ($1\le i_1 \le n-1$,
$n\le i_2 \le 2n-2$) in an open set $\mathscr{O}_V$ in $V$ such that
$\{v^{i_1}\}$ are the local coordinates in the open set
$\mathscr{O}_V\cap L_V$ in $L_V$. And we extend $\{v^i\}$ to local
coordinates $\{v^i,x\}$ for $X$ with $x\equiv 0$ along $V$ and so
that $x=x^1+{\rm i}x^2$ with $J(\frac{\partial}{\partial
x^1})=\frac{\partial}{\partial x^2}$ and $J(\frac{\partial}{\partial
x^2})=-\frac{\partial}{\partial x^1}$, moreover, we require that
$\{v^{i_1},x^1\}$ are just local coordinates for $L$.

\begin{Lemma}\label{L-normal}{\rm (normal form)}.
Suppose that $\Sigma$ is a smooth connected Riemann surface with
boundary and $u$ is a $(J,\nu)$-holomorphic map which intersects $V$
at $\mathbbm{p}$ or $\mathbbm{q}$. Then either (1) $u(\Sigma)\subset
V$, or (2) there exist an integer $K>0$ (depending on $\mathbbm{p}$ or $\mathbbm{q}$) and a nonzero $a_0\in\C$
such that  near interior intersection point
\begin{equation}\label{normal}
u(z',\bar{z'})=(\mathbbm{q}^i+O(|z'|),\ a_0z'^K+O(|z'|^{K+1})),
\end{equation} or near boundary intersection point
\begin{equation}\label{normal-1}
u(z,\bar{z})=(\mathbbm{p}^i+O(|z|),\ c_0z^K+O(|z|^{K+1})),
\end{equation}
in the local coordinates $z'$ or $z$ and $\{v^i,x\}$, respectively, where
$O(|z|^k)$ denotes a function of $z$ and $\bar{z}$ that vanished to
order $k$ at $z=0$. In particular, writing $z=r+\sqrt{-1}\ s$ and
restricting (\ref{normal-1}) to the real part, we have
\begin{equation}\label{normal-2}
u(z,\bar{z})|_\R=u(r)=(\mathbbm{p}^i+O(r),\ c_0r^K+O(r^{K+1})),
\end{equation}
in the local coordinates $\{v^{i_1},x^1\}$ in $L$, where $c_0$ is a
nonzero real number.
\end{Lemma}
{\it Proof}. When we consider an interior intersection point $\mathbbm{q}$, the argument is the same as the one in Lemma 3.4 of \cite{IP1}. When we consider a boundary intersection point $\mathbbm{p}$, the argument can not apply directly, since the map $u$ is only continuous at the boundary. Fortunately, this is just a local problem, the above similar conclusion still holds. Indeed, under the local coordinates $(v^i,x)$ around $\mathbbm{p}$ and $z\in A\subset\mathbb{H}$ around $p_{a\imath}$, every $v^i$ or $x$ is a function $:A\to B\subset\mathbb{C}$ which is only continuous on $A\cap\{{\rm Im} z=0\}$. Since the map $u$ is locally holomorphic, by definition, each $v^i$ (or $x$) can be extended to a holomorphic function $v^i_{\mathbb{C}}$ (or $x_{\mathbb{C}}$) $:U\to\mathbb{C}$, where $U$ is an open neighborhood of $A$ in $\mathbb{C}$. Then the argument in \cite{IP1} can go through for $u_\mathbb{C}$ and obtain the formula (\ref{normal-1}) by restricting to $\mathbb{H}$. \qed

\section{Moduli space of  $V$-regular and  $V$-stable maps }\label{SEC-5}

In the rest of the paper, we only consider the special case that
$L\cap V=\emptyset$. Also we {\bf assume that all moduli spaces are
generically irreducible}, that is they generically contain no
multiply covered maps which, for instance, can be verified by the
imposed semi-positivity conditions on $(X,L)$ and $V$.

\subsection{$V$-regular open stable maps}
\begin{Definition}\label{$V$-regular}
Given a codimension two symplectic submanifold $V$ of $(X,\omega)$.
A stable $(J,\nu)$-holomorphic map
$u:(\Sigma,\partial\Sigma)\rightarrow(X,L)$ is called $V$-{\rm
regular} if no component of its domain is mapped entirely into $V$
and if neither any marked point nor any double point is mapped into
$V$.
\end{Definition}
The set of $V$-regular open maps forms an open subset of the space of
open stable maps, denote it by ${\mathcal M}^V(\Sigma,X,L,{\bf d})$. Denote by
${\mathcal M}^V_{{\bf k},l,\mathbbm{l}}(\Sigma,X,L,{\bf d})$ the
space of $V$-regular maps with marked points.

We denote by
$\mathbbm{s}=(s_1,\cdots,s_\mathbbm{l})$ the list of multiplicities
of interior intersection points, where
$s_j$ are positive integer numbers. And we define the
{\it degree, length}, and {\it order} of
$\mathbbm{s}$ by
$${\rm deg}\ \mathbbm{s}:=\sum_{j=1}^\mathbbm{l}s_j,\ \ \ \ \ \ \ \ {\rm
length}(\mathbbm{s}):=\mathbbm{l},\ \ \ \ \ \ \ \ \ {\rm
ord}(\mathbbm{s}):=|\mathbbm{s}|=\prod_{j=1}^{\mathbbm{l}}s_j.
$$
The vector
$\mathbbm{s}$ labels the component of ${\mathcal M}^V_{{\bf k},l,\mathbbm{l}}(\Sigma,X,L,{\bf d})$.
We denote each
$\mathbbm{s}$-labeled component of ${\mathcal M}^V_{{\bf
k},l,\mathbbm{l}}(\Sigma,X,L,{\bf d})$ by
$${\mathcal M}^{V,\mathbbm{s}}_{{\bf k},l,\mathbbm{l}}(\Sigma,X,L,{\bf d})\subset\overline{\mathcal M}_{{\bf k},l+\mathbbm{l}}(\Sigma,X,L,{\bf d}).$$
Forgetting the additional $\mathbbm{l}$ marked points defines a projection
$${\mathcal M}^{V,\mathbbm{s}}_{{\bf k},l,\mathbbm{l}}(\Sigma,X,L,{\bf d})$$
\bq\label{down}\downarrow\eq
$${\mathcal M}^{V}_{{\bf k},l}(\Sigma,X,L,{\bf d})$$ onto one component of
${\mathcal M}^{V}_{{\bf k},l}(\Sigma,X,L,{\bf d})$, which is the disjoint union
of such components. So for each fixed
$\mathbbm{s}$, (\ref{down}) is a covering map to its
image whose deck transformation group is the the group of renumberings of the
$\mathbbm{l}$ interior marked points.

For any integer $m\ge 1$, let $\mathbb{J}^{V,m}$ be the completion
of the space $\mathbb{J}^V$ of all smooth pairs $(J,\nu)$ in the
$C^m$-topology. Then $\mathbb{J}^{V,m}$ is a $C^m$-smooth Banach
manifold. Denote  by $B_{{\bf k},l}^{*,V,1,p}(\Sigma,X,L,{\bf d})$ the irreducible part of ambient Sobolev space of
$V$-regular maps, which is also a Banach manifold. Let
$\mathfrak{M}^V_{{\bf k},l}=B_{{\bf k},l}^{*,V,1,p}(\Sigma,X,L,{\bf
d})\times\mathbb{J}^{V,m}$, then we define the following universal
moduli space
$$\mathcal{UM}^V_{{\bf k},l}=\mathcal{UM}^V_{{\bf k},l}(\Sigma,X,L,{\bf d}):=\{({\bf u},(J,\nu))\ | \
(J,\nu)\in\mathbb{J}^{V,m},{\bf u}\in{\mathcal M}^V_{{\bf k},l}(\Sigma,X,L,{\bf
d})\}$$ and its irreducible part
$$\mathcal{UM}^{* V}_{{\bf k},l}=\{(u,\vec{z},\vec{w},J,\nu)\in \mathfrak{M}^V_
{{\bf k},l}\  |\  \bar{\partial}_{J,\nu}u=0 \}.$$ The irreducible
part of the universal moduli space involving the intersection data
is denoted by $\mathcal{UM}^{*\ V, \mathbbm{s}}_{{\bf
k},l,\mathbbm{l}}$. When $p>2$, by elliptic regularity, each element $u$ of $\mathcal{UM}^{* V}_{{\bf k},l}$ is smooth.

\begin{Proposition}\label{prop5.1}
The irreducible universal moduli space $\mathcal{UM}^{* V}_{{\bf k},l}$ is a
smooth Banach submanifold of $\mathfrak{M}^V_{{\bf k},l}$.
\end{Proposition}
Proof.  There is an infinite dimensional vector bundle $${\mathcal
E}\rightarrow\mathfrak{M}^V_{{\bf k},l},$$ with each fiber ${\mathcal
E}_{({\bf u},J,\nu)}:=L^p(\Sigma,\Omega^{0,1}(u^*TX)).$ The $(J,
\nu)$-holomorphic equation defines a section of this bundle by
\begin{equation}\label{univer-section}
{\mathcal S}:\mathfrak{M}^V_{{\bf k},l} \rightarrow{\mathcal E}, \ \  {\mathcal
S}({\bf u},J,\nu)(x)=\bar{\partial}_Ju(x)-\nu(u(x)).
\end{equation}
Since $\mathcal{UM}^{* V}_{{\bf k},l}={\mathcal S}^{-1}(0)$, one has to show that ${\mathcal S}$ is transverse to the zero
section. Let ${\mathcal S}({\bf u},J,\nu)=0$. We have
$$T_{\bf u}B_{{\bf k},l}^{*,V,1,p}(\Sigma,X,L,{\bf d})=W^{V,1,p}(\Omega^0(u^*TX,u^*TL))$$
$$T_{(J,\nu)}\mathbb{J}^{V,m}=C^m({\rm End}(TX,J))\oplus C^m({\rm Hom}_J(T\Sigma,TX))$$
where ${\rm End}(TX,J)=\{ A|\  A: TX\rightarrow TX, AJ+JA=0\}$.
Recall that ${\rm Hom}_J(T\Sigma,TX)$ is the space of
anti-$J$-linear homomorphism with respect to the complex structure.
There is an identification
\begin{equation}\label{5-ident}
{\rm Hom}_J(T\Sigma,TX)|_{\Gamma_u}=\Omega^{0,1}_J(u^*TX),
\end{equation}
where $\Gamma_u\subset\Sigma\times X$ is the graph of $u$.

Now we consider the vertical differential
\begin{eqnarray}\label{5-diff}
D{\mathcal S}({\bf u},J,\nu):W^{V,1,p}(\Omega^0(u^*TX,u^*TL))\oplus C^m({\rm End}(TX,J))\nonumber\\
\oplus C^m({\rm Hom}_J(T\Sigma,TX)) \longrightarrow
L^p(\Sigma,\Omega^{0,1}(u^*TX)).
\end{eqnarray}
Then we have
\begin{equation}\label{5-diff-formula}
D{\mathcal S}({\bf u},J,\nu)(\xi, A, \mu)=D_u\xi+\frac{1}{2}u^{*}A\circ du\circ
j_{\Sigma}-\mu|_{\Gamma_u},
\end{equation}
where
$$D_{\bf u}:=D_{\bf u}\bar{\partial}_{(J,\nu)}:W^{V,1,p}(\Omega^0(u^*TX,u^*TL))
\rightarrow L^p(\Omega^{0,1}(u^*TX)).$$

By elliptic regularity theory,  for any $p>2$,  $u$ is in the space
$C^m$ if $(J,\nu)$ is in $C^m$. Moreover, the cokernel of $D_{\bf
u}$ is contained in $C^m(\Omega^{0,1}(u^*TX))$. Since $D_{\bf u}$ is
Fredholm, its cokernel is of finite dimension.  On the other hand,
considering the definition of $\nu$ (\ref{2-def-nu}) and the
identification (\ref{5-ident}), the map
\begin{equation}\label{5-surjective}
\mu|_{\Gamma_u}: C^m({\rm Hom}_J(T\Sigma,TX)) \longrightarrow
C^m(\Omega^{0,1}(u^*TX))
\end{equation}
is surjective. Therfore, by (\ref{5-diff-formula}), $D{\mathcal S}({\bf
u},J,\nu)$ is surjective. Thus, the implicit function theorem
implies that the universal moduli space $\mathcal{UM}^{* V}_{{\bf
k},l}$ is a smooth Banach submanifold of $\mathfrak{M}^V_{{\bf
k},l}$.\qed

\begin{Lemma}\label{dim-M}
Fix a conformal structure $j$ on $\Sigma$. For generic  $V$-compatible pair $(J,\nu)\in\mathbb{J}^V$, the irreducible part of ${\mathcal
M}^{V,\mathbbm{s}}_{{\bf k},l,\mathbbm{l}}(\Sigma,X,L,{\bf d})$ is an orbifold with boundary of dimension
\begin{equation}\label{dim}
{\rm dim}\ {\mathcal M}^{*\ V,\mathbbm{s}}_{{\bf k},l,\mathbbm{l}}(\Sigma,X,L,{\bf
d})=\mu(d)+n(1-g)+k+2(l+\mathbbm{l}-{\rm deg}\mathbbm{s})-{\rm dim} Aut(\Sigma).
\end{equation}
\end{Lemma}

\noindent Proof. The argument will be a slight modification of the ones used in the Lemma 4.2
and Lemma 4.3 of \cite{IP1}, since no boundary points contact $V$.  The main points are as follows.

Similar to the Proposition 17.1 in \cite{FO3} for special case $g=0$
and $l=0$, and the Theorem C.10 in \cite{MS}, it is easy to
calculate the dimension of ${\mathcal M}^{*\ V}_{{\bf k},l+\mathbbm{l}}$ which is
$${\rm dim}{\mathcal M}^{*\ V}_{{\bf k},l+\mathbbm{l}}(\Sigma,X,L,{\bf d})=\mu(d)+n(1-g)+k+2(l+\mathbbm{l})-{\rm dim}Aut(\Sigma).$$
If we show that the irreducible part of universal moduli space $\mathcal{UM}^{*\ V, \mathbbm{s}}_{{\bf k},l, \mathbbm{l}}$ is an orbifold,
then the Sard-Smale transversality theorem implies that for generic $(J,\nu)\in\mathbb{J}^{V,m}$ the irreducible part of the
$V$-regular moduli space is an orbifold.

The Proposition \ref{prop5.1} shows that $\mathcal{UM}^{* V}_{{\bf
k},l}$ 
is a smooth Banach orbifold. Then the point is to show that the contact condition corresponding to each sequence
$\mathbbm{s}$ is transverse. That will imply that the irreducible universal moduli space $\mathcal{UM}^{*\ V, \mathbbm{s}}_{{\bf k},l,\mathbbm{l}}$
is an orbifold and the fomula (\ref{dim}) holds. From the proof of Lemma 4.2 and Lemma 4.3 of \cite{IP1},
we know that the transverality of contact condition arises from the local analysis around those intersecting points with $V$.
Geometrically, that means the $V$-compatibility conditions assure that local variations of the inhomogeneous perturbation and the $V$-regular maps would realize other intersecting multiplicities with given degree.
Since we assume that $L\cap V=\emptyset$ and use the similar $V$-compatibility conditions as  ones in \cite{IP1}, the similar argument can go through.\qed

\subsection{Limits of $V$-regular open maps}\label{SEC-7}

We come to construct a compactification of each component of the
moduli space of $V$-regular open maps.
To compactify ${\mathcal M}^{V,\mathbbm{s}}_{{\bf
k},l,\mathbbm{l}}(\Sigma,X,L,{\bf d})$, we take its closure
\begin{equation}\label{clo}
\mathcal{CM}^{V,\mathbbm{s}}_{{\bf
k},l,\mathbbm{l}}(\Sigma,X,L,{\bf d})
\end{equation}
in the stable moduli space $\overline{\mathcal M}_{{\bf
k},l+\mathbbm{l}}(\Sigma,X,L,{\bf d})$. In fact, the closure
lies in the subset of $\overline{\mathcal M}_{{\bf
k},l+\mathbbm{l}}(\Sigma,X,L,{\bf d})$ consisting of open
stable maps whose last $\mathbbm{l}$ marked points are mapped into
$V$, still with associated multiplicities $\mathbbm{s}$, although the actual order of contact might be infinite.

We will show that the closure is an orbifold with boundary. That is
to prove that the frontier $\mathcal{CM}^{V}\setminus{\mathcal M}^V$ is a
subset of codimension at least 1. Since such frontier is a subset of
the space of stable maps, it is stratified according to the type of
bubble structure of the domain. The following proposition is the
main result in this section describing the structure of the closure
$\mathcal{CM}^V$.
\begin{Proposition}\label{prop-1}
For generic pair $(J,\nu)\in\mathbb{J}^V$, each stratum of the
irreducible part of
$$\mathcal{CM}^{V,\mathbbm{s}}_{{\bf
k},l,\mathbbm{l}}(\Sigma,X,L,{\bf d})\setminus{\mathcal
M}^{V,\mathbbm{s}}_{{\bf
k},l,\mathbbm{l}}(\Sigma,X,L,{\bf d})$$ is an orbifold of
dimension at least one less than the dimension (\ref{dim}) of ${\mathcal M}^{*\ V,\mathbbm{s}}_{{\bf k},l,\mathbbm{l}}(\Sigma,X,L,{\bf d})$.
\end{Proposition}

We follow the way of \cite{IP1}, which studies moduli spaces of closed $(J,\nu)$-maps according to the
different types of limits, to prove this proposition for two basic cases:

Case 1. stable maps with no components or special (({\bf k},l)-marked and double) points lying entirely in $V$;

Case 2. stable maps with some components in $V$ and some off $V$.\\

The full proof of this proposition will finish in the Theorem \ref{Thm-stable}.\\

\noindent {\it Case 1}.\ \ Since we do not admit
the degeneration of domain, each stratum of this type is labeled by
the pair of positive integer numbers $(\mathbbm{b},\mathbbm{i})$ of
boundary and interior double points of their nodal domain curve
$\widehat{\Sigma}$. For such fixed $\widehat{\Sigma}$, the
corresponding stratum is  ${\mathcal M}^{V,\mathbbm{s}}_{{\bf k},l,\mathbbm{l}}(\widehat{\Sigma},X,L,{\bf d})$.

\begin{Lemma}\label{case-1}
In the `Case 1', for generic pair $(J,\nu)\in\mathbb{J}^V$, the
irreducible part of the stratum ${\mathcal
M}^{V,\mathbbm{s}}_{{\bf
k},l,\mathbbm{l}}(\widehat{\Sigma},X,L,{\bf d})$ of
$\mathcal{CM}^V$ is an  orbifold of dimension $\mathbbm{b}+2\mathbbm{i}$
less than  (\ref{dim}).
\end{Lemma}
Proof. The argument is standard and is open version parallel to the ones in \cite{IP1}.  Note that the dimension argument is direct since the degree of freedom of $z_i\in\partial\Sigma$  is one.\qed

\bigskip

\noindent {\it Case 2}. Now we consider the case that there is a sequence of $V$-regular open maps converges to a limit map
$u\in\mathcal{CM}^{V,\mathbbm{s}}_{{\bf k},l,\mathbbm{l}}(\Sigma, X,L,{\bf d})$ whose domain is the
union $\Sigma=\Sigma_1\cup\Sigma_2$, where $\Sigma_1$  is either a bubble domain or a union of several unconnected  bubble domains of (total) genus $g_1$ (for simplicity and without loss of generality, we only consider the one bubble domain case),  $\Sigma_2$
is a bubble domains of genu $g_2$, $\partial\Sigma_2=\emptyset$, and $u$ restricts to a
$V$-regular map $u_1:\Sigma_1\rightarrow X$ and a holomorphic map
$u_2:\Sigma_2\rightarrow V$ into $V$. Limit maps of this type arise
from sequences of $V$-regular maps in which either

(1) two of the
$\mathbbm{l}$-contact intersection points collide in the domain, or

(2) one of the original $l$-marked points whose image sinks into $V$, collides with a contact point.

In either case the collision produces a ghost bubble map
$u_2:\Sigma_2\rightarrow V$ whose energy is at least $\hbar_V$ by Lemma \ref{L-abc}.

Then note that $u_1^{-1}(V)$ consists of the nodal points
$\Sigma_1\cap\Sigma_2$ and some of the last $\mathbbm{l}$ marked points $q_\jmath\in\Sigma$. The
nodes are defined by identifying points $x_j\in\Sigma_1$ and
$y_j\in\Sigma_2$. Since $u_1$ is $V$-regular and $u_1(x_j)\in V$, then Lemma \ref{L-normal}
associates a multiplicity $s'_j$  to each $x_j$, the multiplicity vector is denoted by
$$\mathbbm{s}'=(s_1',s_2',\cdots).$$
For convenience, images of these nodes $x_j$ and $y_j$ is denoted by
vectors $\mathbbm{l}'$. Similarly, since $u$ arises as a limit of $V$-regular maps, the $q_\jmath$, which are limits of the contact points with $V$,
also have associated multiplicities. We split the set
of $q_\jmath$ into the points $\{q^1_\jmath\}$ on $\Sigma_1$  and $\{q^2_\jmath\}$  on $\Sigma_2$, denote by
$$\mathbbm{s}^1=(s^1_1,s^1_2,\cdots)\ \ \ \ \ {\rm and}\ \ \ \ \ \mathbbm{s}^2=(s^2_1,s^2_2,\cdots)$$
the associated multiplicity vectors. Hence we write $u$ as a pair
\begin{equation}\label{pair}
(u_1,u_2)\in {\mathcal M}^{V,\mathbbm{s}^1+\mathbbm{s}'}_{g_1,{\bf k}_1,l_1,\mathbbm{l}_1+\mathbbm{l}'}(\Sigma_1, X,L,{\bf
d}_1)\times \overline{\mathcal M}_{g_2,{\bf k}_2,l+\mathbbm{l}_2+\mathbbm{l}'}(\Sigma_2, V,L_V,{\bf d}_2)
\end{equation}
where ${\bf d}_i=[u_i]$ with ${\bf d}_1+{\bf d}_2={\bf d}$
satisfying the matching conditions $u_1(x_j)=u_2(y_j)$, and $({\bf
k}_1,l_1)+({\bf k}_2,l_2)=({\bf k},l)$.

\begin{Lemma}\label{case-3}
In this case 2, the only elements (\ref{pair}) lying in ${\mathcal
CM}^{V, \mathbbm{s}}$ are those for which there exists
a (singular) section
$\xi\in\Gamma((\Sigma,\partial\Sigma),(u^*_2N_V,u^*_2N_{L_V}))$
nontrivial on at least one component of $\Sigma_2$ with zeros of
order $s_{\jmath}^2$ at $q_{\jmath}^2\in\Sigma$, poles of order $s'_j$ at $y_j$,
and $D^N_{u_2}\xi=0$ where $D^N_{u_2}\xi$ is as in (\ref{D-N}).
\end{Lemma}
Proof.\ It is an open version parallel to the delicate renormalization argument used in proposition 6.6 of
\cite{IP1}. The renormalization will be done in a projective compactification $\mathbb{P}((N_V,N_{L_V})\oplus(\mathbb{C},\mathbb{R}))$ of the normal bundle.
By applying the renormalization argument we see that it can still work well in our open case. We will not repeat the discussion. \qed

\subsection{The space of $V$-stable open maps}\label{SEC-8}

From last subsection we see that the limit of a sequence of $V$-regular
open maps is a stable map whose components are of the basic types
described in Case 1 and 2. Actually, the components of the limit map are
also partially ordered according to the rate at which they sink into
$V$. In this section we will make this precise by extending the concept {\it layer structure},
which is originally introduced in \cite{IP1} for closed domain curves, to our domain of open Riemann surfaces. And then we construct a
compactification of the space of $V$-regular maps. Let $\Sigma$ be a stable curve with boundary.
\begin{Definition}\label{def-layer}
A layer structure on $\Sigma$ is the assignment of a non-negative
integer $\lambda_j$ to each irreducible component $\Sigma_j$ to
$\Sigma$, such that at least one component must have $\lambda_j=0$
or 1.
\end{Definition}
The union of all the components with $\lambda_j=K$ is called the
{\it layer $K$ part of $\Sigma$}, denoted by $\Lambda_K$. Note that
$\Lambda_K$ might not be a connected.

\begin{Definition}\label{mark-layer}
A marked layer structure on a $({\bf k},l+\mathbbm{l})$-marked bordered bubble surface
$\Sigma_{{\bf k},l+\mathbbm{l}}$ is a layer structure on $\Sigma$ together with

{\rm (1)} a  vector $\mathbbm{s}$ recording the multiplicities of the last $\mathbbm{l}$  interior
intersection marked points, and

{\rm (2)} a pair of vectors $(\bbalpha,\bbbeta)$ which assigns multiplicities respectively to each boundary double point of
$\partial\Lambda_K\cap\partial\Lambda_L$, and to each interior double point of $\Lambda_K\cap\Lambda_L$, $K\neq L$.
\end{Definition}Note that double points within a layer are not assigned a multiplicity.

Each layer $\Lambda_K$ then has interior points $q_{K,\jmath}$ of the  type (1) with multiplicity
vector $\mathbbm{s}_K=(s_{K,\jmath})$,  and has double pints with multiplicities. The double points in $\Lambda_K$ that are assigned
multiplicities are divided into two types.

Let $\bbalpha_K^+$ (resp.
$\bbbeta_K^+$) be the vector derived from $\bbalpha$ (resp.
$\bbbeta$) that gives the multiplicities of the boundary (resp.
interior) double points $y^+_{K,i}$ where $\Lambda_K$ meets the
higher layers at boundary (resp. interior part), $i.e.$ the points
$\partial\Lambda_K\cap\partial\Sigma_j$ (resp.
$\Lambda_K\cap\Sigma_j$) with $\lambda_j>K$. Let $\bbalpha_K^-$
(resp. $\bbbeta_K^-$) be the similar vector of multiplicities of
boundary (resp. interior) double points $y^-_{K,i}$ where
$\Lambda_K$ meets the lower layers.

Then we associate operators $D_K^N$ similar to (\ref{D-N})
defined on the layers $\Lambda_K$, $K\ge 1$, as follows. For the
interior double points, the argument is the same as the one in
\cite{IP1}, so we only consider the special case that we assume
there are only boundary double points. For each choice of
$\bbalpha=\bbalpha_K^-=\{\alpha_{K,i}^-\}$ and $\rho$, fix smooth
weighting function $W_{\bbalpha,\rho}$ which has the form
$|z_i|^{\rho+\alpha^-_{K,i}}$ in some local coordinates $z_i$ in a
small half-disc centered at $y^-_{K,i}$ and has no other zeros. Then
given a stable map $u:\Lambda_K\rightarrow V$, denote by
$L^m_{\bbalpha^-_K,\delta}(u^*N_V,u^*N_{L_V})$ the Hilbert space of
all $L^m_{loc}$ sections of $(u^*N_V,u^*N_{L_V})$ over
$(\Lambda_K,\partial\Lambda_K)\setminus\{y^-_{K,i}\}$ which are
finite in the norm
$$\|\xi\|^2_{m,\bbalpha,\delta}=\sum_{l=1}^m\int_{\Lambda_K}|W_{\bbalpha,l+\delta}\cdot \nabla^l\xi|^2.$$
For large $m$ the elements $\xi$ in this space have poles with
$|\xi|\le c |z_i|^{-\alpha_{K,i}^--\delta}$ at each $y_{K,i}^-$ and
have $m-1$ continuous derivatives elsewhere on $\Lambda_K$. For such
$m$, denote by $L^m_{K,\delta}(u^*N_V,u^*N_{L_V})$ the closed
subspace of $L^m_{\bbalpha^-_K,\delta}(u^*N_V,u^*N_{L_V})$
consisting of all sections that vanish to order $r_{K,j}$ and
$s_{K,j}$ at $p_{K,j}$ and $q_{K,j}$ and order $\alpha^+_{K,i}$ at
$y_{K,i}^+$. Hence by elliptic theory for weighted norms, the
operator $D^N$ defines a bounded operator
\begin{equation}\label{DKN}
D^N_K:
L^m_{K,\delta}((\Lambda_K,\partial\Lambda_K),(u^*N_V,u^*N_{L_V}))\rightarrow
L^{m-1}_{K,\delta+1}(T^*\Lambda_K\otimes u^*N_V).
\end{equation}
For generic $0<\delta<1$, $D^N_K$ is Fredholm with
\begin{eqnarray}\label{ind-DKN}
{\rm index}\ D^N_K &=&
\mu_{N_V}([u(\Lambda_K)])+\chi(\Lambda_K)+{\rm deg}\ \bbalpha_K^- -{\rm deg}\ \bbalpha_K^+\nonumber\\&
&+2({\rm deg}\ \bbbeta_K^- -{\rm deg}\mathbbm{s}_K-{\rm deg}\
\bbbeta_K^+)\nonumber\\&=&\chi(\Lambda_K)
\end{eqnarray}
where $\chi(\Lambda_K)$is the Euler characteristic of $\Lambda_K$,
and the index formula would be more general for cases admitting
interior double points. The formula is derived from the fact that
the Euler class of the pair of complex line bundle and totally real
sub-bundle $(u^*N_V,u^*N_{L_V})$ can be computed from the zeros and
poles of a section.

Then under the assumption $L\cap V=\emptyset$, we give the following definition of the stable maps we want.
\begin{Definition}\label{Vstable}
A $V$-stable map is a stable map $u\in\overline{\mathcal M}_{{\bf
k},l+\mathbbm{l}}(\Sigma,X,L,{\bf d})$ together with\\

{\rm (1)} a marked layer structure on its domain $\Sigma$ with
$u|_{\Lambda_0}$ being $V$-regular, and\\

{\rm (2)} for each $K\geq 1$, an element $\xi_K\in\ker D^N_K$
defined on the layer $\Lambda_K$ that is a section nontrivial on
every irreducible component of $\Lambda_K$.
\end{Definition}

Denote by $\overline{\mathcal M}^{V,\mathbbm{s}}_{{\bf k},l,\mathbbm{l}}(\Sigma,X,L,{\bf d})$ the
$\mathbbm{s}$-labeled component of the set of all $V$-stable open maps. This contains the set ${\mathcal
M}^{V,\mathbbm{s}}_{{\bf k},l,\mathbbm{l}}(\Sigma,X,L,{\bf d})$ of $V$-regular maps as
the open subset, $i.e.$ the $V$-stable maps whose entire domain lies in layer 0. Forgetting the data $\xi_K$ defines a map
\begin{equation}\label{forget}
\mathscr{F}:\overline{\mathcal M}^{V,\mathbbm{s}}_{{\bf k},l,\mathbbm{l}}(\Sigma,X,L,{\bf d})\rightarrow\overline{\mathcal M}_{{\bf k},l+\mathbbm{l}}(\Sigma,X,L,{\bf d}).
\end{equation}

Note that the number of layers of a $V$-stable map must be finite.
Assume there are altogether $r+1$ layers
$\Lambda_0,\Lambda_1,\cdots,\Lambda_r$. Each $V$-stable map $({\bf
u},\xi_1,\cdots,\xi_r)$ determines an element of the space\footnote{Here, the space ${\mathcal
H}^V_{(X,L)}$ is not really a parallel generalization of the space ${\mathcal
H}^V_{X}$ appeared in (5.8) of \cite{IP1}. Our space ${\mathcal
H}^V_{(X,L)}$ is less delicate than the space used in \cite{IP1}. However, if $L\cap V=\emptyset$, we can also define the space of homology-intersection data as a covering map of
$\displaystyle H_2(X,L)\times\bigsqcup_{\mathbbm{s}}V_{\mathbbm{s}}$, that will just be the case we consider to define relatively open invariants in this article.}
$${\mathcal
H}^V_{(X,L)}=H_2(X,L)\times\bigsqcup_{\mathbbm{s}}V_{\mathbbm{s}}$$
as follows, where $V_{\mathbbm{s}}$ is the space, diffeomorphic to $V^{\mathbbm{l}(\mathbbm{s})}$, of all sets of contact intersection points with multiplicities.

For sufficiently small
$\varepsilon$, we can push the components in $V$ off $V$ by
composing $u$ with $\exp(\varepsilon^K\xi_K)$ and for each $K$,
smoothing the domain at the nodes
$\displaystyle\Lambda_K\cap(\cup_{L>K}\Lambda_L)$ and smoothly joining the images
where the zeros of $\varepsilon^K\xi_K$ on $\Lambda_K$ approximate
the poles of $\varepsilon^{K+1}\xi_{K+1}$. The resulting map
$$u_\xi=u|_{\Lambda_0}\#\exp(\varepsilon\xi_1)\#\cdots\#\exp(\varepsilon^r\xi_r)$$
is $V$-regular and represents an element or a class in ${\mathcal H}^V_{(X,L)}$. Note that this
class is independent of the choice of the small $\varepsilon$, and
depends on each $\xi_K$ up to a nonzero multiplier.

Thus, we can associate a well-defined map
\begin{equation}\label{proj-int}
\mathscr{H}:\overline{\mathcal M}^{V,\mathbbm{s}}_{{\bf k},l,\mathbbm{l}}(\Sigma,X,L,{\bf d})\rightarrow{\mathcal H}^{V,\mathbbm{s}}_{(X,L),{\bf d}}.
\end{equation}

The following proposition is an open version similar to the
Proposition 7.3 in \cite{IP1}.
\begin{Proposition}\label{compacthm}
There exists a topology on $\overline{\mathcal M}^{V,\mathbbm{s}}_{{\bf k},l,\mathbbm{l}}(\Sigma,X,L,{\bf d})$ such that this space of
$V$-stable maps is compact and the maps $\mathscr{F}$ of
(\ref{forget}) and $\mathscr{H}$ of (\ref{proj-int}) are continuous and differential on each stratum.
\end{Proposition}
Proof.\  The argument is parallel to the one in Proposition 7.3 of \cite{IP1}. The idea is to analyse a sequence of $V$-regular maps and a more
general sequence of $V$-stable maps to define the topology on $\overline{\mathcal M}^{V,\mathbbm{s}}_{{\bf k},l,\mathbbm{l}}(\Sigma,X,L,{\bf d})$.
The modifications are that we would use the bubble convergence Theorem \ref{thm-bubble} to our open case,
and that the argument of open version renormalization which is simply described in the proof of lemma \ref{case-3} will define the structure of a $V$-stable map on the limit map.
\qed

The next theorem is the key result needed to define the relative
invariants, which implies the Proposition \ref{prop-1}.
\begin{Theorem}\label{Thm-stable}
The space of $V$-stable maps is compact and there exists a continuous map
\begin{equation}\label{continuous}
\CD \varepsilon_V  : \overline{\mathcal M}^{V,\mathbbm{s}}_{{\bf
k},l,\mathbbm{l}}(\Sigma,X,L,{\bf d})
 @>{\bf ev}\times \mathscr{H}>>L^{|{\bf k}|}\times X^l \times{\mathcal H}^{V,\mathbbm{s}}_{(X,L),{\bf
 d}}.
\endCD
\end{equation}
Moreover, for generic pair $(J,\nu)\in\mathbb{J}^V$, the complement of ${\mathcal
M}^{V,\mathbbm{s}}_{{\bf k},l,\mathbbm{l}}(\Sigma,X,L,{\bf d})$ in the irreducible part
of $\overline{\mathcal M}^{V,\mathbbm{s}}_{{\bf k},l,\mathbbm{l}}(\Sigma,X,L,{\bf d})$ has codimension at least two.
\end{Theorem}

\noindent Proof. Since we only consider the case without boundary intersection point, it is parallel to the proof of Theorem 7.4 of \cite{IP1}.\qed

In fact, the following lemma is useful for the proof of the theorem above and the discussion below. Since it is a version parallel to the lemma 7.6 of \cite{IP1}, we will not repeat the proof.
We use ${\mathcal S}=\cup \Lambda_K $ (with a marked layer structure but no specified complex structure) to label the strata of the space of $V$-stable maps, and denote these strata by
$\overline{\mathcal M}^{V,\mathbbm{s}}_{{\bf k},l,,\mathbbm{l}}(\mathcal S)$
\begin{Lemma}\label{strata-dim}
Each irreducible stratum $\overline{\mathcal M}^{V,\mathbbm{s}}_{{\bf k},l,,\mathbbm{l}}(\mathcal S)$ is an  orbifold whose dimension is
$$2r+ {\mathcal L}^b_0+2\sum_K {\mathcal L}^i_K $$ less than that in (\ref{dim}), where $r$ is the total number of nontrivial layers without boundary and ${\mathcal L}^b_K$ (resp.
${\mathcal L}^i_K$) is the number of boundary (resp. interior) double points, whether or not intersecting other layers, in layer
$\Lambda_K$, $K\ge 0$.
\end{Lemma}


\bigskip

\subsection{Gluing construction}\label{SEC-gl}

We can refer to Chapter 7 of \cite{FO3} for the construction of Kuranishi structures and especially for gluing construction of local charts of interior and codimension part of moduli space of open stable maps. Since we essentially only consider semi-positive $(X,L)$ and $V$, we do not need to construct the   Kuranishi structure, the situation is much simpler, we just give a rough description.

Note from Theorem \ref{Thm-stable} and Lemma \ref{strata-dim} that, for our special case $L\cap V=\emptyset$,  the $V$-stable maps with some components sinking in $V$ belong to strata of codimension at least 4. The codimension 1 stratum consists of $V$-regular maps with only one boundary node,
denoted by  ${\mathcal M}^{V,\mathbbm{s}}_{{\bf k},l,\mathbbm{l}}(\widehat{\Sigma},X,L,{\bf d})^{1}$, where $\widehat{\Sigma}=\Sigma\cup D^2/z_0'\sim z_0''$,
$z_0'$ is the extra marked point on the $b^{th}$ boundary component $(\partial\Sigma)_b$, $z_0''$ is the extra marked point on $\partial D^2$.
The top stratum of it is the (irreducible) moduli space of $V$-regular maps with smooth domain surface.
Higher codimension strata consist of both $V$-regular maps with more than 1 disc or sphere bubbles on domain and $V$-stable maps with some components sinking in $V$.

We just consider the gluing process for boundary 1-nodal strata. That is, we will glue an open $V$-regular map with  smooth domain with a $(J,\nu)$-holomorphic disc to obtain an open $V$-regular map with a smooth domain.

For simplicity, we only consider the special case that ${\bf k}=\bf 0$, $l=0$, $i.e.$ without marked points.  Denote the moduli space by ${\mathcal M}^{V,\mathbbm{s}}_{\mathbbm{l}}(\Sigma,X,L,{\bf d})$. Let $S$ be a stratum data of the codimension 1 stratum, which records the graph structure of the domain. We admit different complex structures on $\widehat{\Sigma}$, denote the space of complex structures on $\widehat{\Sigma}$ by $M_{\widehat{\Sigma}}$. As for the closed case, one can define an orbifold half-line bundle $$L_{\widehat{\Sigma}}\to M_{\widehat{\Sigma}}$$ which is the normal bundle of $M_{\widehat{\Sigma}}$ in the stable moduli  space of open Riemann surfaces, denoted by  $\overline{M}_\Sigma$. The corresponding codimension 1 stratum is denoted by ${\mathcal M}^{V,\mathbbm{s}}_{\mathbbm{l}}(\widehat{\Sigma}_S,X,L,{\bf d})^{1}$, which is the space of pairs $(u,j)$, where $\widehat{\Sigma}_S$ denotes the 1-nodal surfaces admitting various complex structures.

Then the forgetting map $${\mathcal F}:  {\mathcal M}^{V,\mathbbm{s}}_\mathbbm{l}(\widehat{\Sigma}_S,X,L,{\bf d})^{1}\to M_{\widehat{\Sigma}},$$
$${\mathcal F}(u,j)=j$$ induces an orbifold half-line bundle $${\mathcal F}^*L_{\widehat{\Sigma}}\to {\mathcal M}^{V,\mathbbm{s}}_\mathbbm{l}(\widehat{\Sigma}_S,X,L,{\bf d})^{1}.$$
Given a point $\mu=(u,j,\rho)\in{\mathcal F}^*L_{\widehat{\Sigma}}$, we want to construct a holomorphic map $Gl(\mu)\in {\mathcal M}^{V,\mathbbm{s}}_{\mathbbm{l}}(\Sigma,X,L,{\bf d}).$
We first do the process of  pre-gluing that gives an approximation map $pgl(\mu)$.

\noindent (1) {\it Pre-gluing}. Firstly, the gluing domain surface $\Sigma_{z_0,\rho}$ can be obtained as follows. Let $2T=\log(1+|\rho|)$. We use the holomorphic strip-like coordinates on $\Sigma$ and $D^2$ near $z_0'$ and $z_0''$, that is
$$\Sigma\setminus\{z_0'\}=\Sigma_o\cup\{[-2T,+\infty)\times[0,1]\},$$
$$D^2\setminus\{z_0'' \}=D^{2}_o\cup\{(-\infty,2T]\times[0,1]\},$$ where $\Sigma_o$ and $D^2_o$ are the remaining parts of $\Sigma$ and $D^2$ taking away a small neighborhood of boundary nodes $z_0'$ and $z_0''$, respectively.

Then we cut off the part $[0,\infty)\times[0,1]$ of $\Sigma$ and $(-\infty,0]\times[0,1]$ of $D^2$ with strip-like coordinates and glue the remainders by identifying $\{0\}\times[0,1]$. The derived surface is also denoted by $\Sigma_{z_0,T}$.

For $u=(u_1,u_2)\in{\mathcal M}^{V,\mathbbm{s}}_\mathbbm{l}(\widehat{\Sigma}_S,X,L,{\bf d})^{1}$, the pre-gluing map $pgl(u,j_o,T)$ should be a map on surface $\Sigma_{z_0,T}$ with the same tangents with $V$. We simply denote this map by $\psi=pgl(u,j_o,T): \Sigma_{z_0,T}\to (X,L)$, and define it partially by
$$
\psi(x)=\left\{\begin{array}{ll}
u_1(x) & x\in\Sigma_o,\\
p=u_1(z_0')=u_2(z_0'') & x\in [-T,T]\times[0,1],\\
u_2(x) & x\in D^2_o.
\end{array}
\right.
$$
On the rest part the map can be defined by using a fixed  smooth cutoff function. We do not explicitly write the form of the map on this part and just claim
that it is a minor modification of the construction for the closed domain case, and parallel estimates can show that the $L^p$-norm of
the differential $\bar{\partial}_{j_{o T},J,\nu}\psi$ can be bounded from above by $CT^{1/2p}$.\\

\noindent (2) {\it Right inverses}. With the pre-gluing map $\psi$, one can consider its linearized map $D_{\psi,j_{o\rho}}$ and then construct the right inverse $R_{\psi,j_{oT}}$ to  $D_{\psi,j_{oT}}$. The construction is standard, we will refer to \cite{FO3} and not repeat the argument. \\

\noindent (3) {\it Gluing maps}. With the pre-gluing map $\psi$ and the right inverse $R_{\psi, j_{oT}}$, we can construct a holomorphic curve by perturbing the pre-gluing map.
The construction is still standard. Roughly, when the $L^p$-norm of the differential $\bar{\partial}_{j_{oT},J,\nu}\psi$ is sufficiently small, there exists a unique map $g(u,j_o,T)$ in $L^p(\Lambda^{0,1}_{j_{oT}}\psi^*TM)$ such that the map $\exp_\psi R_{\psi,j_{oT}}g(u,j_o,T)$ is holomorphic. Thus the gluing map is defined to be
$$Gl: {\mathcal F}^*L_{\widehat{\Sigma}}|_U\longrightarrow {\mathcal M}^{V,\mathbbm{s}}_{\mathbbm{l}}(\Sigma,X,L,{\bf d})$$
$$Gl(u,j,T)=pgl(u,j,T)+g(u,j,T).$$
where $U\subset {\mathcal M}^{V,\mathbbm{s}}_\mathbbm{l}(\widehat{\Sigma}_S,X,L,{\bf d})^{1}$ is any proper open subset.

In general, we can glue maps with marked points. We consider the case that there is only one  node $z_0$ which is different  from any marked point, one domain surface is $\Sigma$ with an extra point $z_0'\in(\partial\Sigma)_b$, the other  is a disc $D^2$ with an extra point $z_0''\in\partial D^2$. Then by the same method, the gluing map can be constructed as
\begin{eqnarray}\label{gluemap}
Gl:& {\mathcal M}^{V,\mathbbm{s}^1}_{{\bf k_1}+e_b,l_1,\mathbbm{l}_1}(\Sigma,X,L,{\bf d}_1)_{evb_{z_0'}}\times_{evb_{z_0''}}{\mathcal M}^{V,\mathbbm{s}^2}_{k_2+1,l_2,\mathbbm{l}_2}(D^2,X,L,d_2)\times(0,+\infty)\nonumber\\
     &\longrightarrow{\mathcal M}^{V,\mathbbm{s}}_{{\bf k},l,\mathbbm{l}}({\Sigma},X,L,{\bf d}),
\end{eqnarray}
where $|{\bf k}|=|{\bf k_1}|+k_2$, $l=l_1+l_2$, $\mathbbm{l}=\mathbbm{l}_1+\mathbbm{l}_2$, ${\bf d}={\bf d}_1+d_2$,
$e_b=(0,\cdots,1,\cdots,0)$ is a vector with only the $b^{th}$ component is $1$, others are  zero. This induces  local coordinate charts for the moduli space.


\section{Orietation}\label{SEC-9}

\subsection{Orienting the determinant line bundle over moduli space}

Let ${\rm V}\rightarrow B$ be a vector bundle. Denote by $w_i({\rm V})$ the $i^{th}$
Stiefel-Whitney class of ${\rm V}$. Define two
characteristic classes $p^\pm({\rm V})\in H^2(B,\Z/2\Z)$ by
\begin{equation}\label{9-A5}
p^+({\rm V})=w_2({\rm V}),\ \ \ \ \ p^-({\rm V})=w_2({\rm
V})+w_1({\rm V})^2.
\end{equation}
From \cite{KT} we know that $p^\pm({\rm V})$ is the obstruction to
the existence of a $Pin^\pm$ structure on ${\rm V}$.
\begin{Definition}\label{def-9.1}
Given a symplectic Lagrangian pair $(X,L)$, we say $L$ is relatively
$Pin^\pm$ if
$$p^\pm(TL)\in{\rm Im}(i^*: H^2(X,\Z/2\Z)\rightarrow H^2(L,\Z/2\Z)).$$
and is $Pin^\pm$ if $p^\pm(TL)=0$. If $L$ is  $Pin^\pm$, we define a
$Pin^\pm$ structure for L to be a $Pin^\pm$ structure for $TL$. If
$L$ is  relatively $Pin^\pm$, we define a relative $Pin^\pm$
structure for $(X,L)$ consists of the choices of

1. a triangulation for the pair $(X,L)$,

2. an oriented vector bundle $\V$ over the three skeleton of $X$
such that $w_2(\V)=p^\pm(TL)$,

3. a $Pin^\pm$ structure on $TL|_{L^{(3)}}\oplus \V|_{L^{(3)}}.$
\end{Definition}

Let us first introduce some notations. Recall that
$${\mathcal M}^{V,\mathbbm{s}}_{{\bf
k},l,\mathbbm{l}}(\Sigma,X,L,{\bf d})$$ is the moduli space of
$(J,\nu)$-holomorphic $V$-regular open maps $u :
(\Sigma,\partial\Sigma)\rightarrow(X, L)$ with $k_a$ marked points
$z_{a1},\cdots,z_{ak_a}$ on each boundary component
$(\partial\Sigma)_a$ and $l$ marked points $w_1,\cdots,w_l$ and
additional $\mathbbm{l}$ interior intersecting marked points
$q_1,\cdots,q_\mathbbm{l}$ on $\Sigma$ such that
$u_*([\Sigma,\partial\Sigma])=d$ and
$u|_{(\partial\Sigma)_a*}([(\partial\Sigma)_a])=d_a$, and
$\mathbbm{l}>0$. The compactification $\overline{\mathcal M}^{V,\
\mathbbm{s}}_{{\bf k},l,\mathbbm{l}}(\Sigma,X,L,{\bf d})$ is the space of
$V$-stable open maps. We can define the evaluation maps at the
$({\bf k},l)$-marked points and intersection points
$$evb_{ai}:\overline{\mathcal
M}^{V,\mathbbm{s}}_{{\bf k},l, \mathbbm{l}}(\Sigma,X,L,{\bf d})\rightarrow L,\ \ \ i=1,\cdots,k_a,\ a=1,\cdots,m,$$
$$evi_j:\overline{\mathcal M}^{V,\mathbbm{s}}_{{\bf k},l, \mathbbm{l}}(\Sigma,X,L,{\bf d})\rightarrow X,\ \ \ j=1,\cdots,l,$$
$$evi^I_{\jmath}: \overline{\mathcal M}^{V,\mathbbm{s}}_{{\bf k},l, \mathbbm{l}}(\Sigma,X,L,{\bf d})\rightarrow V,\ \ \ \jmath=1,\cdots,\mathbbm{l}.$$ In fact, the
moduli space above can be considered as the zero locus of a Fredholm
section of a Banach bundle. We denote by $B^{1,p}(\Sigma,X,L,{\bf d})$ the
Banach manifold of $W^{1,p}$ maps
$u:(\Sigma,\partial\Sigma)\rightarrow(X,L)$ such that
$u_*([\Sigma,\partial\Sigma])=d$ and
$u|_{(\partial\Sigma)_a*}([(\partial\Sigma)_a])=d_a$. And define
$$B^{1,p}_{{\bf k},l+\mathbbm{l}}(\Sigma,X,L,{\bf d}):=B^{1,p}(\Sigma,X,L,{\bf d})\times
\prod_a(\partial\Sigma)_a^{k_a}
\times\Sigma^{l+\mathbbm{l}}\setminus\bigtriangleup,$$ where
$\bigtriangleup$ denotes the subset of the product in which two
marked points coincide. Elements of $B^{1,p}_{{\bf
k},l+\mathbbm{l}}(\Sigma,X,L,{\bf d})$ are denoted by ${\bf
u}=(u,\vec{z},\vec{w},\vec{q})$, where $\vec{z}=(z_{ai})$,
$\vec{w}=(w_j)$, $\vec{q}=(q_\jmath)$.

Then the Banach space bundle ${\mathcal E}\rightarrow B^{1,p}_{{\bf
k},l+\mathbbm{l}}(\Sigma,X,L,{\bf d})$ is defined fiberwise with $${\mathcal
E}_{\bf u}:=L^p(\Sigma,\Omega^{0,1}(u^*TX)).$$ We define the section
of this bundle as
$$\bar{\partial}_{(J,\nu)}:B^{1,p}_{{\bf
k},l+\mathbbm{l}}(\Sigma,X,L,{\bf d})\rightarrow{\mathcal E}$$
$$\bar{\partial}_{(J,\nu)}u=du\circ j_{\Sigma}+J\circ du-\nu(\cdot,u(\cdot),{\bf u}),$$
which is the $\nu$-perturbed Cauchy-Riemann operator. We require the
inhomogeneous term
$$\nu\in\Gamma(\Sigma\times X\times B^{1,p}_{{\bf
k},l+\mathbbm{l}}(\Sigma,X,L,{\bf d}), {\rm
Hom}(\pi^*_1T\Sigma,\pi^*_2TX)),$$ such that

(1) $\nu$ is $(j_\Sigma,J)$-anti-linear, $i.e.$ $\nu\circ
j_\Sigma=-J\circ\nu$;

(2) $\nu|_{\partial\Sigma\times X\times B^{1,p}_{{\bf
k},l+\mathbbm{l}}(\Sigma,X,L,{\bf d})}$ carries the sub-bundle
${\pi'}_1^*T\partial\Sigma\subset{\pi'}_1^*T\Sigma$ to the
sub-bundle ${\pi'}_2^*(JTL)\subset{\pi'}_2^*TX$, where $\pi_i$,
$i=1,2$, is the projection from $\Sigma\times X\times B^{1,p}_{{\bf
k},l+\mathbbm{l}}(\Sigma,X,L,{\bf d})$ to the $i^{th}$ factor.

The vertical component of the linearization of this section is
denoted by
$$D:=D\bar{\partial}_{(J,\nu)}:TB^{1,p}_{{\bf
k},l+\mathbbm{l}}(\Sigma,X,L,{\bf d})\rightarrow{\mathcal E}.$$

\begin{Definition}\label{B-regular}
A $W^{1,p}$-map $u\in B^{1,p}_{{\bf k},l+\mathbbm{l}}(\Sigma,X,L,{\bf d})$
is called $V$-regular if no component of its domain is mapped
entirely into $V$ and if neither any of the $|k|+l$ marked points
nor any double points are mapped into $V$.
\end{Definition}
Thus the parameterized $V$-regular moduli space can be regarded as
the zero locus of the section $\bar{\partial}_{(J,\nu)}$ which is
denoted by
\begin{equation}\label{9-para-mod-1}
\widetilde{\mathcal M}^{V,\mathbbm{s}}_{{\bf
k},l,\mathbbm{l}}(\Sigma,X,L,{\bf d}):=\bar{\partial}^{-1}_{J,\nu}(0)\cap
B^{1,p,\ V}_{{\bf k},l+\mathbbm{l}}(\Sigma,X,L,{\bf d}).
\end{equation}
Alternatively, similar to \cite{So}, we can also define our moduli space ${\mathcal
M}^{V,\mathbbm{s}}_{{\bf k},l,\mathbbm{l}}(\Sigma,X,L,{\bf d})$ as an
appropriate section (or say slice) of the reparameterization group
action, we refer the reader to the section 4 in \cite{So} for
details of construction.

We simply denote $B^V:=B^{1,p,\ V}_{{\bf
k},l+\mathbbm{l}}(\Sigma,X,L,{\bf d})$. The evaluation maps can also be
similarly defined on this larger space
$$evb_{ai}:B^V\rightarrow L,\ \ \
i=1,\cdots,k_a,\ a=1,\cdots,m,$$
$$evi_j:B^V\rightarrow X,\ \ \
j=1,\cdots,l,$$
$$evi^I_{\jmath}: B^V\rightarrow V,\ \ \
\jmath=1,\cdots,\mathbbm{l},\ \ \ \ \ \ \ \ \ \ {\rm if}\ \
\mathbbm{l}\ge 1.$$ such that
$$evb_{ai}({\bf u})=u(z_{ai}),\ \ \ evi_j({\bf u})=u(w_j), \ \ \ evi^I_\jmath({\bf u})=u(q_\jmath).$$
The total evaluation map is denoted by
$${\bf ev}:=\prod_{a,i}evb_{ai}\times \prod_jevi_j\times \prod_{\mathbbm{l}\ge 1,\jmath} evi^I_\jmath$$
\begin{equation}\label{9-1-ev}
{\bf ev}:B^{1,p,\ V}_{{\bf k},l+\mathbbm{l}}(\Sigma,X,L,{\bf d})\rightarrow
L^{|{\bf k}|}\times X^l\times V^{\mathbbm{l}}.
\end{equation}

Then let
\begin{equation}\label{9-det-line}
{\mathcal L}:={\rm det}(D)\rightarrow B^{1,p,\ V}_{{\bf
k},l+\mathbbm{l}}(\Sigma,X,L,{\bf d})
\end{equation}
be the determinant line bundle of the family of Fredholm operators
$D$ restricted to $B^V$. And let
\begin{equation}\label{det-LL}
{\mathcal L}'={\mathcal L}\otimes\bigotimes_{a,i}evb^*_{ai}{\rm det}(TL)^*.
\end{equation}

The following two propositions are the fundamental results of
orienting the moduli space. We distinguish the following two
situations\\

(1) $L$ is orientable and provided with an orientation. In this
situation no restriction is required;

(2) $L$ is nonorientable, then we suppose the number of marked
points on each boundary component satisfies
\begin{equation}\label{9-bound}
k_a\cong w_1(d_a)+1,\ \ \ \ ({\rm mod}\ 2);
\end{equation}

Denote by $M^{(i)}$ the $i$-skeleton of a manifold $M$. Then we
define a subspace
$$\mathfrak{B}:=\{(u,\vec{z},\vec{w},\vec{q})\in B^{1,p,\ V}_{{\bf
k},l+\mathbbm{l}}(\Sigma,X,L,{\bf d})\ |\
u:(\Sigma,\partial\Sigma)\rightarrow(X^{(3)},L^{(3)}) \}.$$ We first
prove a lemma

\begin{Lemma}\label{9-A-orient}
Let $(X,L)$ be a symplectic Lagrangian pair and $L$ be relatively $Pin^\pm$. The combination of the relative $Pin^\pm$ structure of $(X,L)$ and
the orientations of $L$ if it is orientable determines a canonical orientation of ${\mathcal L}'|_{\mathfrak{B}}$.
\end{Lemma}
Proof. We just need canonically orient each individual line ${\mathcal
L}'_{\bf u}$ for each ${\bf u}\in\mathfrak{B}$ in a way that varies
continuously in families. Recall that the relative $Pin^\pm$
structure of $(X,L)$ gives a vector bundle $\V\rightarrow X^{(3)}$
and a $Pin^\pm$ structure on $\V|_{L^{(3)}}\oplus TL|_{L^{(3)}}$.
Denote simply by
\begin{equation}\label{9-V-R}
\V_\R:=\V|_{L^{(3)}},\ \ \ \ \ \ \ \V_\C:=\V\otimes\C.
\end{equation}
Choosing an arbitrary Cauchy-Riemann operator $D_0$ on $u^*\V_\C$,
we consider the operator $D_{\bf u}\oplus D_0$,
$$D_{\bf u}\oplus D_0:TB_{\bf u}\oplus W^{1,p}(u^*\V_\C,u^*\V_\R)\longrightarrow
{\mathcal E}\oplus L^p(\Omega^{0,1}(u^*\V_\C)),$$ That is
$$W^{1,p}_{\mathbbm{s}}(u^*(TX\oplus\V_\C),u|^*_{\partial\Sigma}(TL\oplus\V_\R))\oplus\R^{|{\bf k}|}
\oplus\C^{l+\mathbbm{l}}\rightarrow
L^p(\Omega^{0,1}(u^*(TX\oplus\V_\C))).$$ We remark that the choice
of $D_0$ is irrelevant since the space of real linear Cauchy-Riemann
operators on $u^*\V_\C$ is contractible.

Thus, we have a short exact sequence of Fredholm operators
$$0\rightarrow D_{\bf u}\rightarrow D_{\bf u}\oplus D_0\rightarrow D_0\rightarrow 0.$$
By Lemma \ref{A-iso} there exists an isomorphism $${\rm det}(D_{\bf
u})\simeq{\rm det}(D_{\bf u}\oplus D_0)\otimes{\rm det}(D_0)^*.$$
Tensor by $\bigotimes_{a,i}(evb^*_{ai}{\rm det}(TL)^*)_{\bf u}$ on
both sides,
\begin{eqnarray}
{\mathcal L}'_{\bf u}&=&{\rm det}(D_{\bf
u})\otimes\bigotimes_{a,i}(evb^*_{ai}{\rm det}(TL)^*)_{\bf u}\nonumber\\
 &\simeq& {\rm det}(D_{\bf
u}\oplus D_0)\otimes{\rm
det}(D_0)^*\otimes\bigotimes_{a,i}(evb^*_{ai}{\rm det}(TL)^*)_{\bf
u}.
\end{eqnarray}
Actually, we only need to orient
\begin{equation}\label{9-A-det}
{\mathcal L}'_{\bf u}\otimes{\rm det}(D_0)\simeq{\rm det}(D_{\bf
u}\oplus D_0)\otimes\bigotimes_{a,i}(evb^*_{ai}{\rm det}(TL)^*)_{\bf
u},
\end{equation}
since the Cauchy-Riemann $Pin$ boundary value problem
$$\underline{D}_0=(\Sigma,u^*\V_\C,u^*\V_\R,\mathfrak{P}_0,D_0)$$
determines a canonical orientation on ${\rm det}(D_0)$. Indeed, note
that any real vector bundle over a Riemann surface $\Sigma$ with
nonempty boundary admits a $Pin$ structure because the second
cohomology $H^2(\Sigma)$ is trivial. So we can choose a $Pin$
structure $\widetilde{\mathfrak{P}}_0$ on $u^*\V\rightarrow\Sigma$
and define $\mathfrak{P}_0$ to be its restriction to
$u^*\V_\R\rightarrow\partial\Sigma$. By the Lemma 2.11 in \cite{So},
the canonical orientation on ${\rm det}(D_0)$ does not depend on the
choice of $\widetilde{\mathfrak{P}}_0$.

Recall that the relative $Pin$ structure on $(X,L)$ gives a $Pin$
structure on $TL|_{L^{(3)}}\oplus \V|_{L^{(3)}}$, and so gives a
$Pin$ structure on $u|_{\partial\Sigma}^*(TL\oplus \V_\R)$. The
notation is simplified since $u\in\mathfrak{B}$.

If $L$ is orientable and given an orientation, we have induced
orientations on  $u|_{\partial\Sigma}^*(TL\oplus \V_\R)$. Thus we
consider the restricted $Pin$ boundary value problem
$\underline{D_{\bf u}\oplus D_0}$, the Proposition \ref{A.1.1} and
the Remark after it ensures there exists a canonical orientation on
det$(D_{\bf u}\oplus D_0)$.  Then the orientation of $L$ gives the
orientation of det$(TL)$. Therefore, we provide a canonical
orientation on the right-hand side of (\ref{9-A-det}).

If $L$ is non-orientable, on each boundary component
$(\partial\Sigma)_a$ such that $k_a\ne 0$, choose arbitrarily an
orientation on $(evb^*_{a1}TL)_{\bf u}$. Still note that each
boundary component $(\partial\Sigma)_a$ has an orientation induced
from the natural orientation on $\Sigma$. For each $i\in[2,k_a]$ if
$k_a\ne 0$, we can obtain an orientation on $(evb^*_{ai}TL)_{\bf u}$
by trivializing $u|^*_{\partial\Sigma}TL$ along the oriented line
segment in $(\partial\Sigma)_a$ from $z_{a1}$ to $z_{ai}$. If for
some $a$, $u|^*_{(\partial\Sigma)_a}TL$ is orientable, and if
$k_a\ne 0$, then the orientation on $(evb^*_{a1}TL)_{\bf u}$ induces
an orientation on $u|^*_{(\partial\Sigma)_a}TL$; otherwise, if
$k_a=0$, then we arbitrarily choose an orientation on
$u|^*_{(\partial\Sigma)_a}TL$. Thus for such $a$,
$u|_{(\partial\Sigma)_a}^*(TL\oplus \V_\R)$ also has an orientation
by choosing an orientation on $\V_\R$. Similar to the above argument
$u|_{\partial\Sigma}^*(TL\oplus \V_\R)$ admits a $Pin$ structure
given by the relative $Pin$ structure on $(X,L)$. Then the
Proposition \ref{A.1.1} and the Remark after it still apply to
$\underline{D_{\bf u}\oplus D_0}$, and det$(D_{\bf u}\oplus D_0)$
has a canonical orientation. Hence we also provide a canonical
orientation on the right-hand side of (\ref{9-A-det}). Note that
under the additional assumption (\ref{9-bound}), the choice of
orientation on $(evb^*_{a1}TL)_{\bf u}$ is not important. Because
changing the orientation on $(evb^*_{a1}TL)_{\bf u}$ will change all
the orientations on $u|_{(\partial\Sigma)_a}^*(TL\oplus \V_\R)$ and
on $(evb^*_{ai}TL)_{\bf u}$, $i\in[2,k_a]$ it induces, by
Proposition \ref{A.1.1} and Remark after it, the orientation on
det$(D_{\bf u}\oplus D_0)$ will also change. Then the assumption
(\ref{9-bound}) ensures that the number of orientation changes is
even. So the orientation on the right-hand side of (\ref{9-A-det})
is invariant.

Note that the construction above varies continuously in a
one-parameter family, we thus canonically oriented ${\mathcal
L}'|_{\mathfrak{B}}$.\qed

\begin{Proposition}\label{9-A.2.1}
Under the assumptions in Lemma \ref{9-A-orient}, the combination of
the orientations of $L$  if it is orientable and the choice of
relative $Pin^\pm$ structure $\mathfrak{P}$ on $L$ provides a
canonical orientation on ${\mathcal L}'$, that is to say, there exists a
canonical isomorphism of line bundles
\begin{equation}\label{9-iso-L}
{\mathcal L}\widetilde{\longrightarrow}\bigotimes_{a,i}evb^*_{ai}{\rm
det}(TL).
\end{equation}
\end{Proposition}
Proof. It suffices to provide a canonical orientation for the fiber
${\mathcal L}_{\bf u}'$ over each ${\bf u}\in B^V=B^{1,p,\ V}_{{\bf
k},l+\mathbbm{l}}(\Sigma,X,L,{\bf d})$ such that it varies continuously
with ${\bf u}$.

By definition the relative $Pin^\pm$ structure gives a triangulation
of the pair $(X,L)$. The map
$u:(\Sigma,\partial\Sigma)\rightarrow(X,L)$ is homotopic to a map
$\hat{u}:(\Sigma,\partial\Sigma)\rightarrow(X^{(2)},L^{(2)})$ by
using simplicial approximation. The homotopy map is denoted by
$$\Phi:[0,1]\times(\Sigma,\partial\Sigma)\rightarrow(X,L).$$
We will show that the choice of map $\Phi$ is unique up to homotopy.
Indeed, let $\Phi'$ be another such map. We can get a new map by
concatenating $\Phi$ and $\Phi'$
$$\Phi\#\Phi':[-1,1]\times(\Sigma,\partial\Sigma)\rightarrow(X,L).$$
We can use simplicial approximation again to homotope $\Phi\#\Phi'$
to a map into the three skeleton $(X^{(3)},L^{(3)})$. By retaking
suitable parameters, we get a homotopy map from $\Phi$ to $\Phi'$,
denote it by
$$\Psi:[0,1]^2\times(\Sigma,\partial\Sigma)\rightarrow (X,L),$$
satisfying $$\Psi(0,t)=\Phi(t),\ \ \ \ \ \Psi(1,t)=\Phi'(t),\ \ \ \
\ \Psi(s,0)=u.$$ That is to say $\Phi$ is unique up to homotopy.

On the other hand, the two maps $\Phi$ and $\Psi$ can be considered
as maps from $[0,1]$ and $[0,1]^2$ to $B^V$, respectively. Note that
$\hat{u}\in\mathfrak{B}$, by the Lemma \ref{9-A-orient}, the
relative $Pin^\pm$ structure of $(X,L)$ determines a canonical
orientation of ${\mathcal L}'|_{\mathfrak{B}}$. Then the orientation on
${\mathcal L}'|_{\hat{\bf u}}$ induces an orientation of ${\mathcal
L}'|_{\bf u}$ by trivializing $\Phi^*{\mathcal L}'$. Such orientation is
the same as the one induced by any other homotopy $\Phi'$ since we
can trivialize $\Psi^*{\mathcal L}'$. It is easy to see such induced
orientation on ${\mathcal L}'|_{\bf u}$ varies continuously with ${\bf
u}$. Therefore, ${\mathcal L}'$ admits a canonical orientation. \qed

The isomorphism (\ref{9-iso-L}) doesn't involve the effect of the
ordering of the marked points on boundary components. In fact,
$B^{1,p,\ V}_{{\bf k},l+\mathbbm{l}}(\Sigma,X,L,{\bf d})$ consists many
connected components, each one corresponds to each ordering of
boundary marked points. Denote
$$\varpi=(\varpi_1,\cdots,\varpi_m)$$ where each $\varpi_a$ is a
permutation of the integers $1,\cdots,k_a$. We define the sign of
$\varpi$
\begin{equation}\label{9-varpi}
{\rm sign}(\varpi):=\sum_a{\rm sign}(\varpi_a).
\end{equation}
Denote by $B^V_\varpi$ the component of $B^{1,p,\ V}_{{\bf
k},l+\mathbbm{l}}(\Sigma,X,L,{\bf d})$ in which $\vec{z}$ are ordered in
$\partial\Sigma$ by $\varpi$. Now we modify the canonical
isomorphism in Proposition \ref{9-A.2.1} by

\begin{Definition}\label{9-A-ord-1} When {\rm dim}$L\simeq 0$ ({\rm mod} 2) we
define the canonical isomorphism to be the isomorphism constructed
in the Proposition \ref{9-A.2.1} twisted by $(-1)^{{\rm
sign}(\varpi)}$ over the component $B^{1,p}_{{\bf
k},l,\varpi}(\Sigma,L,{\bf d})$. Otherwise, we define the
isomorphism to be exactly the isomorphism constructed in the
Proposition \ref{9-A.2.1}.
\end{Definition}

Now we consider the orientation of moduli space of $V$-stable maps.
We will restrict attention to the special case: $V$-stable map of
two components, one of which is from the original Riemann surface
$\Sigma$, and the other of which is a disc bubble, no component is
mapped into $V$;


Recall that the domain is equipped with a marked layer structure
(see Definition \ref{mark-layer}). Similar to (\ref{DKN}), we can
define associated operators $D_K^N$ on the layers $\Lambda_K$, $K\ge
1$.

\begin{Definition}\label{B-stable} A $V$-stable $W^{1,p}$-map is a  map
$u\in B^{1,p,\ V}_{{\bf
k},l+\mathbbm{l}}(\Sigma,X,L,{\bf d})$ together with\\

{\rm (i)} a marked layer structure on its (nodal) domain $\Sigma$
with
$u|_{\Lambda_0}$ being $V$-regular, and\\

{\rm (ii)} for each $K\geq 1$, an element $\xi_K\in\ker D^N_K$
defined on the layer $\Lambda_K$ that is a section nontrivial on
every irreducible component of $\Lambda_K$.
\end{Definition}

In particular, denote the space of $V$-stable $W^{1,p}$-maps with
only the layer $K=0$ components by $B^{1,p,\ V}_{{\bf
k},l,\mathbbm{l}}(\Sigma,X,L,{\bf d},0)$, and denote the space of
$V$-stable $W^{1,p}$-maps with only the layer $K=1$ components by
$B^{1,p,\ V}_{{\bf
k},l,\mathbbm{l}}(\Sigma,X,L,{\bf d},1)$.\\

We will give a description in more detail. Suppose that only one
disc bubbles off the boundary component $(\partial\Sigma)_b$ along
with $k''$ of the marked points on $(\partial\Sigma)_b$ and $l''$ of
the interior marked points. The domain is a nodal surface
$$\widehat{\Sigma}=\Sigma\cup D^2/z_0'\sim z_0''$$
with $|{\bf k}'|+1$ boundary marked points and $l'+\mathbbm{l}$ interior marked points on $\Sigma$, and $k''+1$ boundary marked
points and $l''$ interior marked points on $D^2$ such that
$${\bf k}'=(k_1,\cdots,k',\cdots,k_m),\  \ k'=k_b-k'',\ \ l=l'+l''.$$
We denote by $z_0'$ (resp. $z_0''$) the extra marked point on $\Sigma$ (resp. $D^2$)
which is different from any $z_{bi}$. Denote by
\begin{eqnarray}
B^{V\#} &=& B^{1,p,\ V}_{{\bf
k},\sigma,l,\rho,\mathbbm{l},\varrho}(\widehat{\Sigma},X,L,{\bf
d}',d'',0)\nonumber\\
   &:=& B^{1,p,\ V}_{{\bf k}'+e_b,l',\mathbbm{l}'}(\Sigma,X,L,{\bf d}',0)_{evb_{z_0'}}\times_{evb_{z_0''}}B^{1,p, V}_{ k''+1,l'', \mathbbm{l}''}(D^2,X,L,d'',0)\nonumber
\end{eqnarray}
the space of $V$-regular $W^{1,p}$ stable maps with only one disc
bubbling off and no component is mapped into $V$ as above. The disc
bubble represents the class $d''\in H_2(X,L)$, denote
$${\bf d}'=(d',d_1,\cdots,d_b',\cdots,d_m)\in H_2(X,L)\oplus
H_1(L)^{\oplus m},$$ satisfying $d'+d''=d$, $d_b'+\partial d''=d_b$,
and $e_b=(0,\cdots,1,\cdots,0)$ is the vector with only the $b^{th}$
element equal to 1, others are zeros.

In the notation of the $W^{1,p}$ space above, $\sigma\subset[1,k_b]$
and $\rho\subset[1,l]$ denote the subsets of boundary and interior
bubbling-off marked points, respectively. And
$\varrho\subset[1,\mathbbm{l}]$ denote the subsets of $V$-intersecting
interior bubbling-off marked points.

We write the element ${\bf u}=({\bf u}',{\bf u}'')\in B^{1,p,\
V}_{{\bf k},\sigma,l,\rho,\mathbbm{l},\varrho}(\hat{\Sigma},X,L,{\bf
d}',d'',0)$, where $${\bf u}'\in B':= B^{1,p,\ V}_{{\bf k}'+e_b,l',\mathbbm{l}'}(\Sigma,X,L,{\bf d}',0),$$
$${\bf u}''\in B'':=B^{1,p,\ V}_{k''+1,l'',\mathbbm{l}''}(D^2,X,L,d'',0)$$ such that
$$evb_{z_0'}({\bf u}')=evb_{z_0''}({\bf u}''),$$ where $evb_{z_0'}$ (resp. $evb_{z_0''}$) is the evaluation
map at $z_0'$ (resp. $z_0''$). For each such ${\bf u}\in B^{V\#}$,
we denote by
$$\widehat{\Sigma}_{\bf u}:=\Sigma\cup D^2/z'_0\sim z_0''.$$
the associated domain curve. The stable map is
$u:(\widehat{\Sigma},\partial\widehat{\Sigma})\rightarrow(X,L)$, and the
node of $\widehat{\Sigma}$ is denoted by $z_0$.

If we denote the two natural projections by
$$p_1:B^{V\#}\rightarrow B',\ \ \  \ \ \ p_2:B^{V\#}\rightarrow B'',$$then
we can define the Banach bundle ${\mathcal E}^\#\rightarrow B^{V\#}$ by
$${\mathcal E}^\#:=p_1^*{\mathcal E}'\oplus p_2^*{\mathcal E}''.$$ with fiber
$${\mathcal E}^\#_{\bf u}:=L^p(\widehat{\Sigma}_{\bf u},\ \Omega^{0,1}(u'^*TX)\oplus\Omega^{0,1}(u''^*TX) ).$$
For $(J,\nu)\in\mathbb{J}^V$, we denote by
$$\bar{\partial}^\#_{(J,\nu)}:B^{V\#}\rightarrow{\mathcal E}^\#$$
the section given by the $\nu$-perturbed Cauchy-Riemann operator.
The vanishing set of this section is the parameterized
one-disc-bubble moduli space, denoted by
\begin{equation}\label{9-para-mod-2}
(\bar{\partial}^\#_{(J,\nu)})^{-1}(0)=\widetilde{\mathcal M}^{\ V,\mathbbm{s}}_{{\bf k},\sigma,l,\rho,\mathbbm{l},\varrho}(\widehat{\Sigma},X,L,{\bf d}',d'',0).
\end{equation}
Its vertical linearization is
\begin{equation}\label{DD}
D^\#:=D\bar{\partial}^\#_{(J,\nu)}:TB^{V\#}\rightarrow{\mathcal E}^\#.
\end{equation}
Similarly, we have two operators
\begin{equation}\label{D12}
D':TB^{V\#}\rightarrow p_1^*{\mathcal E}' ,\hspace{1cm} D'': TB^{V\#}\rightarrow p_2^*{\mathcal E}''.
\end{equation}
Denote the determinant line bundle of family of Fredholm operators
$D^\#$ by $${\mathcal L}^\#:={\rm det}(D^\#)\rightarrow B^{V\#}.$$
Denote
$${\mathcal L}^{\#'}={\mathcal L}^\#\otimes \bigotimes_{a,i}(evb^*_{ai}{\rm
det}(TL)^*.$$ We define a subspace
$$\mathfrak{B}^\#:=\{(u,\vec{z},\vec{w},\vec{q})\in B^{V\#}\ |\
u:(\widehat{\Sigma},\partial\widehat{\Sigma})\rightarrow(X^{(3)},L^{(3)})
\}.$$ We can obtain a lemma similar to the Lemma \ref{9-A-orient}

\begin{Lemma}\label{9-A-orient-2}
Let $(X,L)$ be a symplectic Lagrangian pair and $L$ be relatively $Pin^\pm$. The combination of the relative $Pin^\pm$ structure of $(X,L)$ and
the orientation of $L$ if it is orientable determines a canonical orientation of ${\mathcal L}^{\#'}|_{\mathfrak{B}^\#}$.
\end{Lemma}
Proof. It is similar to the proof of Lemma \ref{9-A-orient} and the
arguments in Proposition 3.3 of \cite{So}, the modifications take
place when we apply the Proposition \ref{A.1.1} and the Remark after
it to the two restricted $Pin$ boundary value problems
$$\underline{D_{\bf u}'\oplus D_0'}\ \ \ \  {\rm and}\ \ \ \  \underline{D_{\bf
u}''\oplus D_0''}.$$ And we have an isomorphism
\begin{equation}\label{9-det-iso}
{\rm det}(D_{\bf u}^\#\oplus D_0^\#)\simeq {\rm det}(D_{\bf
u}'\oplus D_0')\otimes {\rm det}(D_{\bf u}''\oplus D_0'')\otimes
evb_0^*{\rm det}(TL\oplus V_\R)^*_{\bf u}.
\end{equation}\qed

Using the same argument of homotopy uniqueness in the proof of
Proposition \ref{9-A.2.1}, we can similarly obtain the following
\begin{Proposition}\label{9-A.2.2}
The combination of the orientations of $L$ if it is orientable and
the choice of relative $Pin^\pm$ structure $\mathfrak{P}$ on $L$
provide a canonical orientation on ${\mathcal L}^{\#'}$, that is to say,
there exists a canonical isomorphism of line bundles
\begin{equation}\label{9-iso-LL}
{\mathcal L}^\#\widetilde{\longrightarrow}\bigotimes_{a,i}evb^*_{ai}{\rm
det}(TL).
\end{equation}
\end{Proposition}

In order to involve the effect of the ordering of the marked points,
we also need modify the isomorphism in the Proposition
\ref{9-A.2.2}. Recall from the proof above that the  marked points
on each $(\partial\widehat{\Sigma})_a$ can be canonically ordered. So we
can divide the space $B^{1,p}_{{\bf k},l}(\Sigma,L,{\bf d})$ into
components $B^{1,p}_{{\bf k},l,\varpi}(\Sigma,L,{\bf d})$. To be
consistent with the isomorphism defined in Definition
\ref{9-A-ord-1}, we can modify the canonical isomorphism
(\ref{9-iso-LL}) and have the following definition.
\begin{Definition}\label{9-A-ord-2} When {\rm dim}$L\simeq 0$ ({\rm mod} 2) we
define the canonical isomorphism to be the isomorphism constructed
in the Proposition \ref{9-A.2.2} twisted by $(-1)^{{\rm
sign}(\varpi)}$ over the component $B^{1,p}_{{\bf
k},l,\varpi}(\Sigma,L,{\bf d})$. Otherwise, we define the
isomorphism to be exactly the isomorphism constructed in the
Proposition \ref{9-A.2.2}.
\end{Definition}

\subsection{Compatiblility with gluing}
We study the compatibility of the orientation constructed above with the gluing process described in subsection \ref{SEC-gl}.
The argument is a minor modification  of the one in section 8.3 of \cite{FO3}.
For $u=(u_1,u_2)$, we denote the evaluation maps at the boundary marked points for two components and node by
$$evb^1_{a_1i_1}:\widetilde{\mathcal M}^{\ V, \mathbbm{s}}_{{\bf k},\sigma,l,\rho,\mathbbm{l},\varrho}(\widehat{\Sigma},X,L,{\bf d}',d'',0) \rightarrow L,$$
$$evb^2_{a_2i_2}:\widetilde{\mathcal M}^{\ V, \mathbbm{s}}_{{\bf k},\sigma,l,\rho,\mathbbm{l},\varrho}(\widehat{\Sigma},X,L,{\bf d}',d'',0)\rightarrow L,$$
$$evb_0:\widetilde{\mathcal M}^{\ V, \mathbbm{s}}_{{\bf k},\sigma,l,\rho,\mathbbm{l},\varrho}(\widehat{\Sigma},X,L,{\bf d}',d'',0)\rightarrow L.$$
For $u_1$ and $u_2$, we denote the evaluation maps at the extra boundary marked points $z_0'$ and $z_0''$ by
$$evb_{z_0'}:\widetilde{\mathcal M}^{V,\mathbbm{s}^1}_{{\bf k_1},l_1,\mathbbm{l}_1}(\Sigma,X,L,{\bf d}',0)\rightarrow L,$$
$$evb_{z_0''}:\widetilde{\mathcal M}^{V,\mathbbm{s}^2}_{k_2,l_2,\mathbbm{l}_2}(D^2,X,L,d'',0)\rightarrow L.$$
Recall we have defined operators $D^\#$, $D'$ and $D''$ in (\ref{DD}) and (\ref{D12}). Let
$${\mathcal T}={\mathcal L}^{\#'}|_{\widetilde{\mathcal M}^{\ V, \mathbbm{s}}_{{\bf k},\sigma,l,\rho,\mathbbm{l},\varrho}(\widehat{\Sigma},X,L,{\bf d}',d'',0)}
=\det(D^\#)
\otimes\bigotimes_{a,i}evb_{ai}^*\det(TL),$$
$${\mathcal T}_1=\det(D')\otimes\bigotimes_{a_1,i_1}(evb^1_{a_1i_1})^*\det(TL),$$
$${\mathcal T}_2=\det(D'')\otimes\bigotimes_{a_2,i_2}(evb^2_{a_2i_2})^*\det(TL)$$
be the corresponding modified line bundles over moduli spaces $\widetilde{\mathcal M}^{\ V,\mathbbm{s}}_{{\bf k},\sigma,l,\rho,\mathbbm{l},\varrho}(\widehat{\Sigma},X,L,{\bf d}',d'',0)$.

\begin{Lemma}\label{L-T12}
For any $u=(u_1,u_2)\in\widetilde{\mathcal M}^{\ V,\mathbbm{s}}_{{\bf k},\sigma,l,\rho,\mathbbm{l},\varrho}(\widehat{\Sigma},X,L,{\bf d}',d'',0)$, there exists an isomorphism of modified determinant lines
\begin{equation}\label{iso-T12}
{\mathcal T}|_u\cong (-1)^{(|{\bf k}_1|-1)(k_2-1)}{\mathcal T}_1|_{u}\otimes {\mathcal T}_2|_{u}\otimes evb_0^*\det(TL)^*_u.
\end{equation}
which is orientation preserving.
\end{Lemma}
\noindent Proof.  Recall that we have two natural projections
$$p_1: \widetilde{\mathcal M}^{\ V, \mathbbm{s}}_{{\bf k},\sigma,l,\rho,\mathbbm{l},\varrho}(\widehat{\Sigma},X,L,{\bf d}',d'',0) \rightarrow \widetilde{\mathcal M}^{V,\mathbbm{s}^1}_{{\bf k_1},l_1,\mathbbm{l}_1}(\Sigma,X,L,{\bf d}',0),$$
$$p_2: \widetilde{\mathcal M}^{\ V, \mathbbm{s}}_{{\bf k},\sigma,l,\rho,\mathbbm{l},\varrho}(\widehat{\Sigma},X,L,{\bf d}',d'',0) \rightarrow \widetilde{\mathcal M}^{V,\mathbbm{s}^2}_{k_2,l_2,\mathbbm{l}_2}(D^2,X,L,d'',0),$$
and the  bundle
$${\mathcal E}^\#:=p_1^*{\mathcal E}'\oplus p_2^*{\mathcal E}''.$$
We define an operator $$\delta:p_1^*T\widetilde{\mathcal M}^{V,\mathbbm{s}^1}_{{\bf k_1},l_1,\mathbbm{l}_1}(\Sigma,X,L,{\bf d}',0)\oplus p_2^*T\widetilde{\mathcal M}^{V,\mathbbm{s}^2}_{k_2,l_2,\mathbbm{l}_2}(D^2,X,L,d'',0)\rightarrow evb_0^*(TL)$$ by
$$\xi'\oplus\xi''\mapsto d(evb_0')(\xi')-d(evb_0'')(\xi'').$$
So it implies $$T\widetilde{\mathcal M}^{\ V, \mathbbm{s}}_{{\bf k},\sigma,l,\rho,\mathbbm{l},\varrho}(\widehat{\Sigma},X,L,{\bf d}',d'',0) =\ker(\delta).$$
We thus have a diagram of exact sequence
$$
\xymatrixcolsep{3.5pc}
\xymatrixrowsep{2.5pc}
\xymatrix{
{\mathcal E}^\#_u\ar[r] & p_1^*{\mathcal E}_u'\oplus p_2^*{\mathcal E}_u''\ar[r] & 0\\
T\widetilde{\mathcal M}^{\ V,\mathbbm{s}}_{{\bf k},\sigma,l,\rho,\mathbbm{l},\varrho}|_u\ar[u]^-{D_u^\#} \ar[r]& p_1^*T\widetilde{\mathcal M}^{V,\mathbbm{s}^1}_{{\bf k_1},l_1,\mathbbm{l}_1}|_u\oplus p_2^*T\widetilde{\mathcal M}^{V,\mathbbm{s}^2}_{k_2,l_2,\mathbbm{l}_2}|_u\ar[r]^-{\delta}\ar[u]^-{D'_u\oplus D_u''}  & evb_0^*(TL)_u\ar[u].
}
$$
which implies an isomorphism
$$\det D_u^\#\cong\det D'_u\otimes\det D_u''\otimes evb_0^*\det(TL)^*_u.$$
Then (\ref{iso-T12}) holds and the coefficient $(-1)^{(|{\bf k}_1|-1)(k_2-1)}$ records the change of cyclic order of boundary marked points.\qed

\begin{Lemma}
The gluing map
\begin{eqnarray}\label{glmap-5}
Gl: \widetilde{\mathcal M}^{V,\mathbbm{s}^1}_{{\bf k_1},l_1,\mathbbm{l}_1}(\Sigma,X,L,{\bf d}',0)_{evb_{z_0'}}\times_{evb_{z_0''}}\widetilde{\mathcal M}^{V,\mathbbm{s}^2}_{k_2,l_2,\mathbbm{l}_2}(D^2,X,L,d'',0)\nonumber\\
\longrightarrow\widetilde{\mathcal M}^{\ V,\mathbbm{s}}_{{\bf k}, l,\mathbbm{l}}(\Sigma,X,L,{\bf d}'+d'',0) \hspace{2cm}
\end{eqnarray}
 is orientation preserving in the sense of the following isomorphism,
 \begin{equation}\label{iso-LT12}
Gl^*{\mathcal L}' |_u\cong (-1)^{(|{\bf k}_1|-1)(k_2-1)}{\mathcal T}_1|_{u}\otimes {\mathcal T}_2|_{u},
\end{equation}
where $|{\bf k}|=|{\bf k_1}|+k_2$, $l=l_1+l_2$, ${\mathcal L}'$ is the bundle defined in (\ref{det-LL}).
\end{Lemma}

\noindent Proof.  We claim that the argument is a combination of the conclusion of the lemma above and modification of ones in the lemma 8.3.5 and lemma 8.3.10 of \cite{FO3} which study the case of gluing of discs.  The modification is direct, we will omit the detail.
Note in our case we only assume relative $Pin^{\pm}$ structure on $L$, one has to use the modified determinant lines to study the preserving of orientation. \qed

\begin{Proposition}
The gluing map (\ref{glmap-5}) induces an isomorphism
\begin{eqnarray}\label{P-iso}
\partial{\mathcal M}^{V,\mathbbm{s}}_{{\bf k},l,\mathbbm{l}}({\Sigma},X,L,{\bf d},0)\hspace{8cm}\\
\cong\bigcup(-1)^*{\mathcal M}^{V,\mathbbm{s}^1}_{{\bf k_1}+e_b,l_1,\mathbbm{l}_1}(\Sigma,X,L,{\bf d}',0)_{evb_{z_0'}}
\times_{evb_{z_0''}}{\mathcal M}^{V,\mathbbm{s}^2}_{k_2+1,l_2,\mathbbm{l}_2}(D^2,X,L,d'',0)\nonumber
\end{eqnarray}
compatible with orientations in the  sense of the Lemma above, where the union is taken over all $b=1,\cdots,m$, ${\bf k}$ such that $|{\bf k}|=|{\bf k}_1|+k_2$, $l=l_1+l_2$,
$\mathbbm{l}=\mathbbm{l}_1+\mathbbm{l}_2$ and ${\bf d}={\bf d}'+d''$, $*=(|{\bf k}_1|-1)(k_2-1)+(n+|{\bf k}_1|)$.
\end{Proposition}
\noindent Proof.  The argument is parallel to the one in the Proposition 8.3.3 of \cite{FO3}.
From definition we know that
$${\mathcal M}^{V,\mathbbm{s}}_{{\bf k},l,\mathbbm{l}}({\Sigma},X,L,{\bf d},0)\cong
\widetilde{\mathcal M}^{\ V,\mathbbm{s}}_{{\bf k}, l,\mathbbm{l}}(\Sigma,X,L,{\bf d},0)/{\rm Aut}(\Sigma_{\bf d}),$$
$${\mathcal M}^{V,\mathbbm{s}^1}_{{\bf k}_1+e_b,l_1,\mathbbm{l}_1}({\Sigma},X,L,{\bf d}',0)\cong
\widetilde{\mathcal M}^{\ V,\mathbbm{s}^1}_{{\bf k}_1+e_b, l_1,\mathbbm{l}_1}(\Sigma,X,L,{\bf d}',0)/{\rm Aut}(\Sigma_{{\bf d}'}),$$
$${\mathcal M}^{V,\mathbbm{s}^2}_{k_2+1,l_2,\mathbbm{l}_2}({D^2},X,L,d'',0)\cong\widetilde{\mathcal M}^{V,\mathbbm{s}^2}_{k_2+1,l_2,\mathbbm{l}_2}({D^2},X,L,d'',0)/{\rm Aut}(D^2_{d''}).$$
We can  take   two boundary marked points $x_0$ and $x_1$ on the gluing component  $(\partial\Sigma)_b$, and $y_0, y_1$ on $\partial D^2$.  Denote by
${\rm Aut}(\Sigma; x_0,x_1)$  (resp. ${\rm Aut}(D^2;y_0,y_1)$) the bi-holomorphic automorphisms group of
$\Sigma$ fixing $x_0$ and $x_1$ (resp. $y_0$ and $y_1$). Then we have
\begin{equation}\label{p-Aut}
\R_{p}\times\mathrm{Aut}(\Sigma_{{\bf d}'+d''}; x_0,y_1)\cong\mathrm{Aut}(\Sigma_{{\bf d}'}; x_0,x_1)\times\mathrm{Aut}(D^2_{d''}; y_0,y_1),
\end{equation}
where $\R_p$ is regarded as the space of gluing parameters. Note that
$$\widetilde{\mathcal M}^{\ V,\mathbbm{s}}_{{\bf k}-2e_b, l,\mathbbm{l}}(\Sigma,X,L,{\bf d}'+d'',0)\cong
{\mathcal M}^{V,\mathbbm{s}}_{{\bf k},l,\mathbbm{l}}({\Sigma},X,L,{\bf d}'+d'',0)\times\mathrm{Aut}(\Sigma_{{\bf d}'+d''}; x_0,y_1),$$
$$\widetilde{\mathcal M}^{\ V,\mathbbm{s}^1}_{{\bf k}_1-e_b, l_1,\mathbbm{l}_1}(\Sigma,X,L,{\bf d}',0)\cong
{\mathcal M}^{V,\mathbbm{s}^1}_{{\bf k}_1+e_b,l_1,\mathbbm{l}_1}({\Sigma},X,L,{\bf d}',0)\times\mathrm{Aut}(\Sigma_{{\bf d}'}; x_0,x_1),$$
$$\widetilde{\mathcal M}^{V,\mathbbm{s}^2}_{k_2-1,l_2,\mathbbm{l}_2}(D^2,X,L,d'',0)\cong{\mathcal M}^{V,\mathbbm{s}^2}_{k_2+1,l_2,\mathbbm{l}_2}(D^2,X,L,d'',0)\times\mathrm{Aut}(D^2_{d''}; y_0,y_1).$$
From (\ref{glmap-5}) we have
\begin{eqnarray}
\widetilde{\mathcal M}^{\ V, \mathbbm{s}}_{{\bf k}-2e_b, l,\mathbbm{l}}(\Sigma,X,L,{\bf d}'+d'',0)\hspace{8cm}\nonumber\\
\cong\widetilde{\mathcal M}^{V,\mathbbm{s}^1}_{{\bf k_1}-e_b,l_1,\mathbbm{l}_1}(\Sigma,X,L,{\bf d}',0)_{evb_{z_0'}}\times_{evb_{z_0''}}\widetilde{\mathcal M}^{V,\mathbbm{s}^2}_{k_2-1,l_2,\mathbbm{l}_2}(D^2,X,L,d'',0).\nonumber
\end{eqnarray}
Then
\begin{eqnarray}
{\mathcal M}^{V,\mathbbm{s}}_{{\bf k},l,\mathbbm{l}}({\Sigma},X,L,{\bf d}'+d'',0)\times\mathrm{Aut}(\Sigma_{{\bf d}'+d''}; x_0,y_1)\hspace{6cm}\nonumber\\
\cong[{\mathcal M}^{V,\mathbbm{s}^1}_{{\bf k}_1+e_b,l_1,\mathbbm{l}_1}({\Sigma},X,L,{\bf d}',0)\times\mathrm{Aut}(\Sigma_{{\bf d}'}; x_0,x_1)]_{evb_{z_0'}}\hspace{4cm}\nonumber\\
\times_{evb_{z_0''}}[{\mathcal M}^{V,\mathbbm{s}^2}_{k_2+1,l_2,\mathbbm{l}_2}(D^2,X,L,d'',0)\times\mathrm{Aut}(D^2_{d''}; y_0,y_1)]\nonumber\\
\cong(-1)^{n+|{\bf k}_1|}\mathrm{Aut}(\Sigma_{{\bf d}'}; x_0,x_1)\times[{\mathcal M}^{V,\mathbbm{s}^1}_{{\bf k}_1+e_b,l_1,\mathbbm{l}_1}({\Sigma},X,L,{\bf d}',0)_{evb_{z_0'}}\hspace{2.3cm}\nonumber\\
\times_{evb_{z_0''}}{\mathcal M}^{V,\mathbbm{s}^2}_{k_2+1,l_2,\mathbbm{l}_2}(D^2,X,L,d'',0)]\times\mathrm{Aut}(D^2_{d''}; y_0,y_1)\nonumber\\
\cong(-1)^{n+|{\bf k}_1|}\R_p\times[{\mathcal M}^{V,\mathbbm{s}^1}_{{\bf k}_1+e_b,l_1,\mathbbm{l}_1}({\Sigma},X,L,{\bf d}',0)_{evb_{z_0'}}\hspace{4.5cm}\nonumber\\
\times_{evb_{z_0''}}{\mathcal M}^{V,\mathbbm{s}^2}_{k_2+1,l_2,\mathbbm{l}_2}(D^2,X,L,d'',0)]\times\mathrm{Aut}(\Sigma_{{\bf d}'+d''}; x_0,y_1),\ \ \ \ \mathrm{by}\ (\ref{p-Aut}).\nonumber
\end{eqnarray}
Therfore, the boundary components of ${\mathcal M}^{V,\mathbbm{s}}_{{\bf k},l,\mathbbm{l}}({\Sigma},X,L,{\bf d},0)$ consist of
$${\mathcal M}^{V,\mathbbm{s}^1}_{{\bf k}_1+e_b,l_1,\mathbbm{l}_1}({\Sigma},X,L,{\bf d}',0)_{evb_{z_0'}}
\times_{evb_{z_0''}}{\mathcal M}^{V,\mathbbm{s}^2}_{k_2+1,l_2,\mathbbm{l}_2}(D^2,X,L,d'',0).$$
From the orientation defined above, the conclusion holds.\qed




\section{Involution and sign}\label{SEC-10}

Here the argument is different from the discussion by Ionel-Parker
\cite{IP1} for closed relative GW-invariant, we must deal with the
codimension 1 frontier of the compactification of the moduli space
of the open $V$-regular maps. We will modify the original method
in \cite{So}. Note in the sequel we still assume $L\cap
V=\emptyset$.

Let us then suppose there exists an anti-symplectic involution
$\phi$ such that $L={\rm Fix}(\phi)$. And suppose $\Sigma$ is
biholomorphic to its conjugation $\bar{\Sigma}$, $i.e.$ there exists
an anti-holomorphic involution $c:\Sigma\rightarrow\Sigma$. Fix
$(J,\nu)\in\mathbb{J}^V_\phi$. Thus, from  the $V$-regular
$(J,\nu)$-holomorphic (resp. $W^{1,p}$) map
$u:(\Sigma,\partial\Sigma)\mapsto(X,L)$ we can define its conjugate
$(J,\nu)$-holomorphic (resp. $W^{1,p}$) map $\tilde{u}=\phi\circ
u\circ c$ representing the homology class $\tilde{d}=[\tilde{u}]$.
So we have an induced map
$$\phi':B^{1,p,\ V}_{{\bf k},l+\mathbbm{l}}(X,L,{\bf
d})\rightarrow B^{1,p,\ V}_{{\bf k},l+\mathbbm{l}}(X,L,\tilde{\bf
d})$$ given by
$${\bf u}=(u,\vec{z},\vec{w},\vec{q})\mapsto\tilde{\bf u}=
(\tilde{u},\ (c|_{\partial\Sigma})^{|{\bf k}|}(\vec{z}),\
c^l(\vec{w}),\ c^\mathbbm{l}(\vec{q})).$$ We denote the relevant Banach
space bundle by $\widetilde{\mathcal E}\rightarrow B^{1,p,\ V}_{{\bf
k},l+\mathbbm{l}}(X,L,\tilde{\bf d})$  with fiber
$$\widetilde{\mathcal E}_{\tilde{\bf u}}:=L^p(\Sigma,\Omega^{0,1}(\tilde{u}^*TX)).$$
And we can similarly get a determinant line bundle of a family of
Fredholm operators $\widetilde{D}:TB^{1,p,\ V}_{{\bf
k},l+\mathbbm{l}}(X,L,\tilde{\bf d})\rightarrow\widetilde{\mathcal E}$
as
$$\widetilde{\mathcal L}:={\rm det}(\widetilde{D})\rightarrow B^{1,p,\ V}_{{\bf
k},l+\mathbbm{l}}(X,L,\tilde{\bf d}).$$ The evaluation maps can also
be similarly defined
$$\widetilde{ev}b_{ai}:B^{1,p,\ V}_{{\bf
k},l+\mathbbm{l}}(X,L,\tilde{\bf d})\rightarrow L,\ \ \
i=1,\cdots,k_a,\ a=1,\cdots,m,$$
$$\widetilde{ev}i_j:B^{1,p,\ V}_{{\bf
k},l+\mathbbm{l}}(X,L,\tilde{\bf d})\rightarrow X,\ \ \
j=1,\cdots,l,$$
$$\widetilde{ev}i^I_\jmath:B^{1,p,\ V}_{{\bf
k},l+\mathbbm{l}}(X,L,\tilde{\bf d})\rightarrow V,\ \ \
\jmath=1,\cdots,\mathbbm{l},$$ by
$$\widetilde{ev}b_{ai}(\tilde{\bf u})=\tilde{u}(c(z_{ai})),\ \ \ \widetilde{ev}i_j(\tilde{\bf u})=\tilde{u}(c(w_j)),
\ \ \
\widetilde{ev}i^I_\jmath(\tilde{\bf u})=\tilde{u}(c(q_\jmath)).$$
The total evaluation map is denoted by
$$\widetilde{\bf ev}:=\prod_{a,i}\widetilde{ev}b_{ai}\times
\prod_j\widetilde{ev}i_j\times\prod_\jmath\widetilde{ev}i^I_\jmath$$
\begin{equation}\label{10-ev-tilde}
\widetilde{\bf ev}: B^{1,p,\ V}_{{\bf
k},l+\mathbbm{l}}(X,L,\tilde{\bf d})\rightarrow L^{|{\bf k}|}\times
X^l\times V^\mathbbm{l}.
\end{equation}
Moreover, we can define a map
$$\Phi:{\mathcal E}\rightarrow\widetilde{\mathcal E}$$
$$\xi\mapsto d\phi\circ\xi\circ dc$$
covering $\phi'$. Also $\Phi$ induces a map $\Psi:{\mathcal
L}\rightarrow\widetilde{\mathcal L}$ covering $\phi'$. Denote by
$$\widetilde{\mathcal L}':=\widetilde{\mathcal L}\otimes\bigotimes_{a,i}\widetilde{ev}b^*_{ai}{\rm
det}(TL)^*.$$ Thus, $\Psi$ also induces a map $$\Psi':{\mathcal
L}'\rightarrow\widetilde{\mathcal L}'$$ covering $\phi'$.

From the Proposition \ref{9-A.2.1} and Definition \ref{9-A-ord-1} we
know that both ${\mathcal L}'$ and $\widetilde{\mathcal L}'$ have canonical
orientation. So the map $\Psi'$ may either preserve the orientation
component or reverse the orientation of some connected components.
We say the sign of $\Psi'$ is 0 if $\Psi'$ preserves the orientation
of ${\mathcal L}'$ to that of $\widetilde{\mathcal L}'$, otherwise, we say
the sign of $\Psi'$ is 1. The following proposition shows an
expression of the sign of $\Psi'$. Recall $g_0$ denotes the genus of
$\Sigma/\partial\Sigma$, and ${\rm dim} L=n$.
\begin{Proposition}\label{1-sign}
If $L$ is relative $Pin^-$, then the map $\Psi'$ has sign
\begin{eqnarray}
\mathfrak{s}^-&\cong&\frac{\mu(d)(\mu(d)+1)}{2}+(1-g_0)n+nm+\deg
\mathbbm{s}+|{\bf
k}|+l+\mathbbm{l}\nonumber\\
  & & +w_2(\V)(\psi(d))+w_1(u^*TL)(\partial
  d)+\sum_{a<b}w_1(u^*TL)(d_a)w_1(u^*TL)(d_b)\nonumber\\
  & & +\sum_aw_1(u^*TL)(d_a)(k_a-1)+(n+1)\sum_a\frac{(k_a-1)(k_a-2)}{2}
  \  \ \mod \ 2.\nonumber
\end{eqnarray}
If $L$ is relative $Pin^+$, then the map $\Psi'$ has sign
\begin{eqnarray}
\mathfrak{s}^-&\cong&\frac{\mu(d)(\mu(d)+1)}{2}+(1-g_0)n+nm+\deg
\mathbbm{s}+|{\bf
k}|+l+\mathbbm{l}\nonumber\\
  & & +w_2(\V)(\psi(d))+\sum_{a<b}w_1(u^*TL)(d_a)w_1(u^*TL)(d_b)\nonumber\\
  & & +\sum_aw_1(u^*TL)(d_a)(k_a-1)+(n+1)\sum_a\frac{(k_a-1)(k_a-2)}{2}
  \  \ \mod \ 2.\nonumber
\end{eqnarray}
Where $\psi:H_*(X,L;\Z/2\Z)\rightarrow H_*(X;\Z/2\Z)$ is a degree 0
homomorphism defined in \cite{So}.
\end{Proposition}


Since the proof of the preceding proposition is similar to the proof
of the next one and we will not actually use it, we do not show the
proof here.

Sometimes, for simplicity, we would denote $\mu=\mu(d)$,
$w_2=w_2(\V)$ and $w_1=w_1(u^*TL)$ if no danger of confusion.
\begin{Corollary}
If ${\rm dim}L\le 3$, then we always have
\begin{eqnarray}\label{dim-<-3}
\mathfrak{s}^-&\cong&\frac{\mu(\mu+1)}{2}+(1-g_0)n+nm+\deg
\mathbbm{s}+|{\bf
k}|+l+\mathbbm{l}\nonumber\\
  & & +\sum_{a<b}w_1(d_a)w_1(d_b)
  +\sum_a w_1(d_a)(k_a-1)\nonumber\\
  & & +(n+1)\sum_a\frac{(k_a-1)(k_a-2)}{2}
  \  \ \mod \ 2.
\end{eqnarray}
In particular, if $\Sigma=D^2$, we have
\begin{eqnarray}\label{D^2-sign}
\mathfrak{s}^-&\cong&\frac{\mu(\mu-1)}{2}+\deg \mathbbm{s}+|{\bf
k}|+l+\mathbbm{l}+\mu(k-1)\nonumber\\
& &+(n+1)\frac{(k-1)(k-2)}{2}
  \  \ \mod \ 2.
\end{eqnarray}
\end{Corollary}
{Proof}.\
When ${\rm dim}L\le 3$, the Wu relations imply that $L$ is
$Pin^-$, we can take the standard $Pin^-$ structure. Note that
$w_2(\V)=0$ and $w_1(\partial d)\cong \mu(d)$ mod 2, thus the
formula is simplified.\qed

Then we define and calculate the sign of a map related to the
boundary of moduli space. Recall that
\begin{eqnarray}
B^{V\#} &=& B^{1,p,\ V}_{{\bf k},\sigma,l,\rho,\mathbbm{l},\varrho}(\hat{\Sigma},X,L,{\bf d}',d'',0)\nonumber\\
&:=& B^{1,p,\ V}_{{\bf k}'+e_b,l',\mathbbm{l}'}(\Sigma,X,L,{\bf d}',0)_{evb_0'}\times_{evb_0''}B^{1,p,\ V}_{k''+1,l'',\mathbbm{l}''}(D^2,X,L,d'',0).\nonumber
\end{eqnarray}
Note that the standard conjugation $c:D^2\rightarrow D^2$ gives a
biholomorphic isomorphism $D^2\simeq \bar{D^2}$. Then from the
involution $\phi$, we just have an induced map
$$\phi_{B''}:B^{1,p,\ V}_{
k''+1,l'',\mathbbm{l}''}(D^2,X,L,d'',0)\rightarrow B^{1,p,\ V}_{
k''+1,l'',\mathbbm{l}''}(D^2,X,L,\widetilde{d''},0)$$ of the second
factor of the fiber product. Then since $L={\rm Fix}(\phi)$,
$\phi_{B''}$ induces an map of the whole fiber product
$$\phi_{B^\#}:B^{V\#}\rightarrow \widetilde{B}^{V\#}$$
$$B^{1,p,\ V}_{{\bf
k},\sigma,l,\rho,\mathbbm{l},\varrho}(\hat{\Sigma},X,L,{\bf
d}',d'',0)\rightarrow B^{1,p,\ V}_{{\bf
k},\sigma,l,\rho,\mathbbm{l},\varrho}(\hat{\Sigma},X,L,{\bf
d}',\widetilde{d''},0).$$

We denote the relevant Banach space bundle by $${\mathcal
E}^\#\rightarrow {B}^{V\#},\ \ \ \ \widetilde{\mathcal E}^\#\rightarrow
\widetilde{B}^{V\#}.$$ And we can similarly get a determinant line
bundle of a family of Fredholm operators
$\widetilde{D}:TB^{V\#}\rightarrow\widetilde{\mathcal E}^\#$ as
$$\widetilde{\mathcal L}^\#:={\rm det}(\widetilde{D})\rightarrow B^{1,p,\ V}_{{\bf
k},\sigma,l,\rho,\mathbbm{l},\varrho}(\hat{\Sigma},X,L,{\bf
d}',\widetilde{d''},0).$$ The obvious evaluation maps are denoted by
$\widetilde{ev}i_{j}$, $\widetilde{ev}b_{ai}$,
$\widetilde{ev}i^{I}_{\jmath}$. Denote by
$$\widetilde{\mathcal L}^{\#'}:=\widetilde{\mathcal L}^\#\otimes\bigotimes_{a,i}\widetilde{ev}b^*_{ai}{\rm
det}(TL)^*.$$

Similarly, we have induced map $\Phi^\#:{\mathcal
E}^\#\rightarrow\widetilde{\mathcal E}^\#$ covering $\phi_{B^\#}$.
Recall that the inhomogeneous term $\nu$ vanishes on bubble
components, so it is $\phi$-invariant. Therefore,
$\bar{\partial}^\#_{(J,\nu)}|_{{\mathcal E}^\#}$ and
$\bar{\partial}^\#_{(J,\nu)}|_{\widetilde{\mathcal E}^\#}$ are two
$\phi_{B^\#}-\Phi^\#$ equivariant sections. Consequently,
$\phi_{B^\#}$ and $\Phi^\#$ induce a map of determinant line bundles
$\Psi^\#:{\mathcal L}^\#\rightarrow\widetilde{\mathcal L}^\#$ covering
$\phi_{B^\#}$. Thus, $\Psi^\#$ also induces a map
$$\Psi^{\#'}:{\mathcal L}^{\#'}\rightarrow\widetilde{\mathcal L}^{\#'}$$
covering $\phi_{B^\#}$. If this map preserves orientation the sign
of it is 0, otherwise the sign is 1.\\

\noindent{\it Remark}. If the homology class $d$ represented by $u$
is $\phi$-anti-invariant, that is, the stable map $u$ is real, then
all the maps above are involutions of their respective objects. In
particular,
we only define and calculate the sign of $\Psi^{\#'}:{\mathcal L}^\#\rightarrow{\mathcal L}^\#$.\\

For stating the formulae for the sign of $\Psi^{\#'}$, we introduce
some new notation. Let
\begin{eqnarray}\label{Ups-1}
\Upsilon^{(1)}(d'',k''):\cong\mu(d'')k'' \cong w_1(\partial d'')k'',
\end{eqnarray}

\begin{eqnarray}\label{Ups-2}
\Upsilon^{(2)}(d_b',d'',k',k''):\cong
\left\{
\begin{array}{ll}
0,& w_1(d_b')=w_1(\partial d'')=0,\\
k' ,& w_1(d_b')=w_1(\partial d'')=1,\\
k''-1 ,& w_1(d_b')=1,\ w_1(\partial d'')=0,\\
k'+k''-1 ,& w_1(d_b')=0,\ w_1(\partial d'')=1.
\end{array}
\right.
\end{eqnarray}

\begin{Proposition}\label{2-sign}
(1) Suppose the marked point $z_{b1}$ does not bubble off, $i.e.$
$1\notin\ \sigma$. Then the map $\Psi^{\#'}$ has sign
\begin{eqnarray}\label{sign-10-1}
\mathfrak{s}^{\#(1)}_{\pm}&\cong&\frac{\mu(d'')(\mu(d'')\pm1)}{2}+\deg
\mathbbm{s}''+
k''+1+l''+\mathbbm{l}''\nonumber\\
  & & +w_2(\V)(\psi(d''))+\Upsilon^{(1)}(d'',k'')\\
  & & +(n+1)\frac{k''(k''-1)}{2}
  \  \ \mod \ 2,\nonumber
\end{eqnarray}
with $+$ in the $Pin^+$ and $-$ in the $Pin^-$ case.

\noindent(2) Suppose  the marked point $z_{b1}$ bubbles off, $i.e.$
$1\in\ \sigma$. Then the map $\Psi^{\#'}$ has sign
\begin{eqnarray}\label{sign-10-2}
\mathfrak{s}^{\#(2)}_{\pm}&\cong&\frac{\mu(d'')(\mu(d'')\pm1)}{2}+\deg
\mathbbm{s}''+
k''+1+l''+\mathbbm{l}''\nonumber\\
  & & +w_2(\V)(\psi(d''))+\Upsilon^{(2)}(d_b',d'',k',k'')
  +w_1(d_b')w_1(\partial d'')\\
  & & +(n+1)[\frac{(k''-1)(k''-2)}{2}+k_b(k''+1)]
  \  \ \mod \ 2,\nonumber
\end{eqnarray}
with $+$ in the $Pin^+$ and $-$ in the $Pin^-$ case.
\end{Proposition}
\noindent{\it Remark}. If $L$ is orientable and ${\rm dim}L=n$ is
odd, note the fact that if $L$ is orientable then $\mu(d'')$ is even
, so $w_1(d_b')=w_1(\partial d'')=0$, then we have
\begin{eqnarray}\label{ori-odd}
\mathfrak{s}^{\#(1)}=\mathfrak{s}^{\#(2)}\hspace{9.5cm}\nonumber\\
\cong\frac{\mu(d'')}{2}+\deg
\mathbbm{s}''+k''+1+l''+\mathbbm{l}''+w_2(\V)(\psi(d''))\  \ \mod \ 2.
\end{eqnarray}

\noindent{\it Proof of Proposition \ref{2-sign}}. The first two
terms in $\mathfrak{s}^{\#(i)}_{\pm}$ come from the formula
(\ref{A4+3}). It is the sign of the conjugation of the determinant
line associated with the restricted $Pin$ boundary problem
$\underline{D_{\bf u}''\oplus D_0''}$ appearing in
(\ref{9-det-iso}), induced from the conjugation on the moduli space
of unmarked discs. The terms $k''+1+l''+\mathbbm{l}''$ account for
conjugation on the configuration space of marked points, satisfying
an extra incidence condition. The term $w_2(V)(\psi(d''))$ accounts
for the change of orientation of the determinant line $\det(D_{\bf
u}''\oplus D_0'')$ arising from the change of $Pin$ structure
induced by the involution $\phi$. Recall that the unique oriented
path from $z\ne z_0$ to $z'$ in the boundary of
$\partial\widehat{\Sigma}$, induced by the complex structure of
$\widehat{\Sigma}$, is very important for determining the canonical
orientation of ${\mathcal L}^{\#'}$. This path will reverse under
conjugation. The terms $\Upsilon^{(1)}$ (resp.
$\Upsilon^{(2)}+w_1(d_b')w_1(\partial d''))$ in
$\mathfrak{s}^{\#(1)}$ (resp. $\mathfrak{s}^{\#(2)}$) reflect this
dependence. Note that  reordering of the boundary marked points will
affect the isomorphism in Definition \ref{9-A-ord-2} when the
dimension of $L$ is even. The last terms account for this
dependence.\qed

\section{Equivalent definition of relatively open invariants}\label{SEC-11}

From the discussion before, we see that we have two moduli spaces of
$V$-regular maps ${\mathcal M}^{V,\mathbbm{s}}_{{\bf k},l,
\mathbbm{l}}(\Sigma,X,L,{\bf d})$ and, corresponding to the anti-symplectic
involution $\phi$, ${\mathcal M}^{V,\mathbbm{s}}_{{\bf k},l,
\mathbbm{l}}(\Sigma, X,L,\tilde{\bf d})$. And we can restrict the two total
evaluation maps (\ref{9-1-ev}) and (\ref{10-ev-tilde}) to have
$${\bf ev}:{\mathcal
M}^{V,\mathbbm{s}}_{{\bf k},l, \mathbbm{l}}(\Sigma, X,L,{\bf d})\rightarrow L^{|{\bf k}|}\times X^l \times V^{\mathbbm{l}},$$
$$\widetilde{\bf ev}:{\mathcal
M}^{V,\mathbbm{s}}_{{\bf k},l, \mathbbm{l}}(\Sigma, X,L,\tilde{\bf d})\rightarrow L^{|{\bf k}|}\times X^l \times V^{\mathbbm{l}}.$$
For
generic choice of points $\vec{x}=(x_{ai})$, $x_{ai}\in L$, and
pairs of points $\vec{y}_{+}=(y_j^+)$, $\vec{y}_{-}=(y_j^-)$ such
that $y_j^+=\phi(y_j^-)$, $j=1,\cdots,l$, and
$\vec{\mathbbm{q}}_{+}=(\mathbbm{q}_\jmath^+)$,
$\vec{\mathbbm{q}}_{-}=(\mathbbm{q}_\jmath^-)$ such that
$\mathbbm{q}_\jmath^+=\phi(\mathbbm{q}_\jmath^-)$,
$\mathbbm{q}_\jmath^\pm\in V$, $\jmath=1,\cdots,\mathbbm{l}$, the
two total evaluation maps will be transverse to
$$\prod_{a,i}x_{ai}\times\prod_jy_j^+
\times\prod_\jmath \mathbbm{q}_\jmath^+\ \ \ {\rm and}\ \ \
\prod_{a,i}x_{ai}\times\prod_jy_j^-\times\prod_\jmath
\mathbbm{q}_\jmath^-$$ in $L^{|{\bf k}|}\times X^l\times
 V^{\mathbbm{l}}$, respectively. For defining invariants,
by index theorem, the following dimension condition would be
satisfied
\begin{equation}\label{11-dim}
(n-1)(|{\bf
k}|+2l)+(n-2)\cdot2\mathbbm{l}=n(1-g)+\mu(d)-2\deg\mathbbm{s}-{\rm dim}\
Aut(\Sigma),
\end{equation}
where $\mu:H_2(X,L)\rightarrow \Z$ denote the Maslov index, $g$
denote the genus of the closed Riemann surface
$\Sigma\cup_{\partial\Sigma}\bar{\Sigma}$ obtained by doubling
$\Sigma$. Note that $\mu(d)=\mu(\tilde{d})$, we can define a number
as
\begin{equation}\label{11-def-1}
M(V,{\bf d},\phi,{\bf k},l,\mathbbm{l})=\#{\bf
ev}^{-1}(\vec{x},\vec{y}_{+},\vec{\mathbbm{q}}_+)+\#\widetilde{\bf
ev}^{-1}(\vec{x},\vec{y}_{-},\vec{\mathbbm{q}}_-),
\end{equation}
where $\#$ denotes the signed count with the sign of a given point,
for example ${\bf u}\in{\bf
ev}^{-1}(\vec{x},\vec{y}_{+},\vec{\mathbbm{q}}_+)$, depending on
whether or not the isomorphism
$$d{\bf ev}_{\bf u}:{\rm det}(T{\mathcal
M}^{V,\mathbbm{s}}_{{\bf k},l, \mathbbm{l}}(\Sigma,X,L,{\bf d}))_{\bf u}
\widetilde{\longrightarrow}\ {\bf ev}^*{\rm det}(T(L^{|{\bf
k}|}\times X^l \times V^\mathbbm{l}))_{\bf u}$$ agrees with the underlying
canonical isomorphism appearing in the Theorem \ref{thm-1.1}
$${\rm det}(T{\mathcal
M}^{V,\mathbbm{s}}_{{\bf k},l, \mathbbm{l}}(\Sigma,X,L,{\bf d}))
\widetilde{\longrightarrow}\bigotimes_{a,i}evb_{ai}^*\det(TL),$$ up
to the action of $\R^+$.

In particular, if $d=\tilde{d}$ we just simply define the number
\begin{equation}\label{11-RN=M}
\mathscr{RN}:=\mathscr{RN}(V,{\bf d},{\bf k},l,\mathbbm{l})
=M(V,{\bf d},\phi,{\bf k},l,\mathbbm{l})=\#{\bf
ev}^{-1}(\vec{x},\vec{y},\vec{\mathbbm{q}}),
\end{equation}
where $(\vec{x},\vec{y},\vec{\mathbbm{q}})$ is a real configuration,
$i.e.$ $l=2c$, $\vec{y}=\{y^+_1,\cdots,y^+_c,y^-_1,\cdots,y^-_c\}$
satisfying $y^+_{j'}=\phi(y^-_{j'})$, $j'=1,\cdots,c$; and
$\mathbbm{l}=2\mathbbm{c}$,
$\vec{\mathbbm{q}}=\{\mathbbm{q}^+_1,\cdots,\mathbbm{q}^+_\mathbbm{c},\mathbbm{q}^-_1,\cdots,\mathbbm{q}^-_\mathbbm{c}\}$
satisfying $\mathbbm{q}^+_{\jmath'}=\phi(\mathbbm{q}^-_{\jmath'})$,
$\jmath'=1,\cdots,\mathbbm{c}$.\\

\noindent {\it Remark}. For such special $d\in H_2(X,L)$,  we can
prove that $\mathscr{RN}$ is an invariant if $L$ is orientable and
${\rm dim}L\le 3$ or if $L$ might not be orientable and ${\rm
dim}L=2$. However, for general homology class $d$, we might not
expect that $M(V,{\bf d},\phi,{\bf k},l, \mathbbm{l})$ is invariant.
In the sequel, we will construct invariants for general homology
class $d$, the proof of invariance of $\mathscr{RN}(V,{\bf d},{\bf
k},l,\mathbbm{l})$ can be considered as a special case of the proof
below
of the invariance of more general invariants.  \\

In order to define relatively open invariants  for general homology
class $d\in H_2(X,L)$, let us introduce more necessary notations. We
denote by $\mathbbm{d}=d_\C$ the doubling of $d$. For any homology
class $\beta\in H_2(X,L)$, denote
$$\bar{\beta}=(\beta,\beta_1,\cdots,\beta_m)\in H_2(X,L)\oplus H_1(L)^{\oplus m},$$  and we
denote the set of $({\bf k},l+\mathbbm{l})$-real configurations by
${\mathcal R}:={\mathcal R}(\vec{x},\vec{y},\vec{\mathbbm{q}})$ which is
\begin{equation}\label{11-real-confg}
\left\{
\begin{array}{l}
\vec{r}=(\vec{x},\vec{\xi},\vec{\lambda}\ )
=(x_{11},\cdots,x_{mk_m},\xi_1,\cdots,\xi_l,
\lambda_1,\cdots,\lambda_\mathbbm{l})\\
\ |\ \xi_j=y^+_j\ {\rm or}\ \xi_j=y^-_j,\ j=1,\cdots,l;
\lambda_\jmath=\mathbbm{q}^+_\jmath\ {\rm or}\
\lambda_\jmath=\mathbbm{q}^-_\jmath,\ \jmath=1,\cdots,\mathbbm{l}.
\end{array}
\right\}.
\end{equation}
Moreover, denote by ${\bf ev}_{(\beta,\vec{r})}$ the total
evaluation map
$${\bf ev}_{(\beta,\vec{r})}:{\mathcal
M}^{V,\mathbbm{s}}_{{\bf
k},l,\mathbbm{l}}(X,L,\bar{\beta})\rightarrow L^{|{\bf k}|}\times
X^l\times V^{\mathbbm{l}},$$
$$(u,\vec{z},\vec{w},\vec{q})\mapsto(\vec{x}=u(\vec{z}),\vec{\xi}=u(\vec{w}),
\vec{\lambda}=u(\vec{q})).$$

Now, we can rewrite the number of (\ref{ROGW-int}) as
\begin{equation}\label{11-inv}
\mathcal {I}:=\mathcal {I}^{\ V,\
\mathbbm{s}}_{X,\phi,g,\mathbbm{d},{\bf
k},l}(\vec{x},\vec{y},\vec{\mathbbm{q}})=\sum_{\forall
\beta:\beta_\C=\mathbbm{d}; \forall\vec{r}\in {\mathcal R}}\#{\bf
ev}^{-1}_{(\beta,\vec{r})}(\vec{x},\vec{\xi},\vec{\lambda}\ )
\end{equation}

To show that the definition (\ref{11-inv}) is independent of the
choices of  $\vec{x}$ and pairs $(\vec{y}_+,\vec{y}_-)$,
$(\vec{\mathbbm{q}}_+,\vec{\mathbbm{q}}_-)$ is equivalent to proving
the expression (\ref{ROGW-int}) is independent of the choices of
${\rm det}(TL)$-valued $n$-forms $\alpha_{ai}$, pairs of $2n$-forms
$(\gamma_j^+,\gamma_j^-)$ and pairs of $(2n-2)$-forms
$(\eta_\jmath^+,\eta_\jmath^-)$, where $\gamma_j^\pm$ represent the
Poincar\'{e} dual of $y^\pm_j$ for $j=1,\cdots,l$, and
$\eta_\jmath^\pm$ represent the Poincar\'{e} dual of
$\mathbbm{q}^\pm_\jmath$ for $\jmath=1,\cdots,\mathbbm{l}$.

\section{Proof of invariance}\label{SEC-12}

In the following, we will
show that the numbers $\mathcal {I}$ (in particular, $\mathscr{RN}$)
are invariants, provided $L$ is orientable and dim$L\le 3$, if $L$
is nonorientable, we suppose that ${\rm dim}L=2$ and the number of
boundary marked points  satisfy some additional conditions.

Suppose that we are given different points of real configuration
$(\vec{x'},\vec{y'}_\pm,\vec{\mathbbm{q}'}_\pm)$ satisfying the same
generic conditions.

Let us denote $${\bf x}:[0,1]\rightarrow L^{|{\bf k}|},\ \ \ {\bf
x}(0)=\vec{x},\ \ {\bf x}(1)=\vec{x'},$$ $${\bf
y}^\pm:[0,1]\rightarrow X^l,\ \ \ {\bf y}^+(t)=\phi({\bf y}^-(t)),$$
$${\bf y}^\pm(0)=\vec{y}_\pm,\ \ \ {\bf y}^\pm(1)=\vec{y'}_\pm,$$
$${\bf
q}^\pm:[0,1]\rightarrow V^\mathbbm{l},\ \ \ {\bf q}^+(t)=\phi({\bf
q}^-(t)),$$
$${\bf q}^\pm(0)=\vec{\mathbbm{q}}_\pm,\ \ \ {\bf q}^\pm(1)=\vec{\mathbbm{q}'}_\pm,$$
$$\Xi:[0,1]\rightarrow X^l,\ \ \ \Xi(0)=\vec{\xi},\ \ \ \Xi(1)=\vec{\xi'},$$
$$\Lambda:[0,1]\rightarrow V^\mathbbm{l},\ \ \ \Lambda(0)=\vec{\lambda},\ \ \ \Lambda(1)=\vec{\lambda'},$$
moreover, we require that $\xi_j=y_j^+$ (resp. $y_j^-$) if and only
if $\xi_j'=y_j'^{+}$ (resp. $y_j'^-$),
$\lambda_\jmath=\mathbbm{q}_\jmath^+$ (resp. $\mathbbm{q}_\jmath^-$)
if and only if $\lambda_\jmath'=\mathbbm{q}_\jmath'^+$ (resp.
$\mathbbm{q}_\jmath'^-$). Denote the set of all paths by
$${\bf R}:={\bf R}({\bf x},{\bf y}^\pm,{\bf q}^\pm)=\{({\bf x},\Xi,\Lambda)\}$$
And denote
\begin{equation}\label{11-W-1}
{\mathcal W}({\bf x},\Xi,\Lambda,\bar{\beta}):={\mathcal M}^{V,\
\mathbbm{s}}_{{\bf k},l,\mathbbm{l}}(X,L,\bar{\beta})_{{\bf
ev}_{(\beta,\vec{r})}}\times_{({\bf
x}\times\Xi\times\Lambda)\circ\triangle}([0,1]),
\end{equation}
\begin{equation}\label{11-W-2}
{\mathcal W}={\mathcal W}({\bf x},{\bf y}^\pm,{\bf
q}^\pm,\mathbbm{d}):=\bigcup_{\begin{small}\begin{array}{c}
\beta:\beta_\C=\mathbbm{d}, \\
({\bf x},\Xi,\Lambda)\in {\bf R}
\end{array}\end{small}}{\mathcal W}({\bf x},\Xi,\Lambda,\bar{\beta}).
\end{equation}
Note that ${\mathcal W}$ gives a smooth oriented cobordism between
$$\bigcup_{\forall
\beta:\beta_\C=\mathbbm{d}; \forall\ \vec{r}\in {\mathcal
R}(\vec{x},\vec{y},\vec{\mathbbm{q}})}{\bf
ev}^{-1}_{(\beta,\vec{r})}(\vec{x},\vec{\xi},\vec{\lambda}\ )$$ and
$$\bigcup_{\forall \beta:\beta_\C=\mathbbm{d}; \forall\ \vec{r'}\in
{\mathcal R}(\vec{x'},\vec{y'},\vec{\mathbbm{q}'})}{\bf
ev}^{-1}_{(\beta,\vec{r'})}(\vec{x'},\vec{\xi'},\vec{\lambda'}\ ).$$

Since in general ${\mathcal W}$ is noncompact, in order to prove the
invariance of $\mathcal {I}:=\mathcal {I}^{\ V,\
\mathbbm{s}}_{X,\phi,g,\mathbbm{d},{\bf k},l}$, we must research the
stable boundary $\partial_G{\mathcal W}$ arising from the Gromov
compactification of ${\mathcal W}$.

Note that for any tuple
$\bar{\beta}=(\beta,\beta_1,\cdots,\beta_m)\in H_2(X,L)\oplus
H_1(L)^{\oplus m}$ the boundary of $V$-stable compactification is
$$\partial\overline{\mathcal
M}^{V,\mathbbm{s}}_{{\bf k},l,
\mathbbm{l}}(\Sigma, X,L,\bar{\beta}):={\mathcal M}^{V,\mathbbm{s}}_{{\bf
k},\sigma,l,\rho,\mathbbm{l},\varrho}(\hat{\Sigma},X,L,\bar{\beta}',\beta'',0),$$
where
$\bar{\beta}'=(\beta',\beta_1,\cdots,\beta_b',\cdots,\beta_m)$,
$\beta'=[u_1(\Sigma)]$, $\beta''=[u_2(D^2)]\in H_2(X,L)$, such that
$\beta=\beta'+\beta''$, $\partial\beta''+\beta_b'=\beta_b$.





For simplicity, denote
$${\mathcal M}(\Sigma)={\mathcal M}^{V,\mathbbm{s}}_{{\bf
k}'+e_b,l',\mathbbm{l}'}(\Sigma,X,L,\bar{\beta}',0),$$
$${\mathcal M}(D^2)={\mathcal M}^{V,\mathbbm{s}}_{
k''+1,l'',\mathbbm{l}''}(D^2,X,L,\beta'',0).$$ Then
\begin{eqnarray}\label{M*M}
\partial\overline{\mathcal M}^{V,\mathbbm{s}}_{{\bf
k},l,\mathbbm{l}}(\Sigma, X,L,\bar{\beta}))&=&{\mathcal M}^{V,\mathbbm{s}}_{{\bf k},\sigma,l,\rho,\mathbbm{l},\varrho}(\hat{\Sigma},X,L,\bar{\beta}',\beta'',0)\nonumber\\
&=&{\mathcal M}(\Sigma)_{evb_{z_0'}}\times_{evb_{z_0''}}{\mathcal M}(D^2).
\end{eqnarray}

We can generically choose $({\bf x},{\bf y}^\pm,{\bf q}^\pm)$ such
that $({\bf x},\Xi,\Lambda)=[{\bf
x},(\Xi',\Xi''),(\Lambda',\Lambda'')]\in{\bf R}$ is transverse to
the total evaluation map
$${\bf ev}_{(\beta'+\beta'',\vec{r})}:{\mathcal M}^{V,\mathbbm{s}}_{{\bf
k},\sigma,l,\rho,\mathbbm{l},\varrho}(\hat{\Sigma},X,L,\bar{\beta}',\beta'',0)\rightarrow
L^{|{\bf k}|}\times X^l\times V^{\mathbbm{l}}.$$ We denote  by
\begin{equation}\label{12-W-bound}
\begin{array}{l}
\partial_G{\mathcal
W}^b_{\sigma,\rho,\varrho,\beta=\beta'+\beta''}(\Xi'',\Lambda'')\\
   \\
\ \ \ \ :={\mathcal M}^{V,\mathbbm{s}}_{{\bf
k},\sigma,l,\rho,\mathbbm{l},\varrho}(\hat{\Sigma},X,L,\bar{\beta}',\beta'',0)_{{\bf
ev}_{(\beta'+\beta'',\vec{r})}} \times_{[{\bf
x}\times(\Xi',\Xi'')\times(\Lambda',\Lambda'')]\circ\triangle}[0,1]
\end{array}\end{equation}
the boundary stratum of the cobordism ${\mathcal W}$ arising from Gromov
compactification. Denote
$$
\begin{array}{l}
\partial_G{\mathcal W}^b_{\sigma, \rho,\varrho,\ \beta',\beta''}\\
  \\
\ \ \ \ :=
\partial_G{\mathcal W}^b_{\sigma, \rho,\varrho,\ \beta'+\beta''}(\Xi'',\Lambda'')
\bigcup \partial_G{\mathcal W}^b_{\sigma, \rho,\varrho,\
\beta'+\widetilde{\beta''}}(\phi(\Xi''),\phi(\Lambda'')),
\end{array}
$$
and
$$
\partial_G{\mathcal W}^b_{\sigma, \rho,\varrho}=\bigcup_{\begin{small}\begin{array}{c}\beta',\beta'':\
(\beta'+\beta'')_\C=\mathbbm{d}\end{array}\end{small}}\partial_G{\mathcal
W}^b_{\sigma, \rho,\varrho,\ \beta',\beta''}\ .
$$ Thus the
total Gromov compactified boundary of ${\mathcal W}$ is
\begin{equation}\label{partial-G-W-1}
\partial_G{\mathcal W}=\bigcup_{
\begin{small}\begin{array}{c} b\in[1,m],\sigma\subset[1,k_b], \\
\rho\subset[1,l],\varrho\subset[1,\mathbbm{l}]\end{array}\end{small}}
\partial_G{\mathcal W}^b_{\sigma, \rho,\varrho}\ .
\end{equation}

Then we have an induced involution map on the boundary of cobordism
$$\phi_{\partial}:\partial_G{\mathcal W}\rightarrow\partial_G{\mathcal W},$$
$$\partial_G{\mathcal W}^b_{\sigma, \rho,\varrho,\
\beta'+\beta''}(\Xi'',\Lambda'') \rightarrow
\partial_G{\mathcal W}^b_{\sigma, \rho,\varrho,\
\beta'+\widetilde{\beta''}}(\phi(\Xi''),\phi(\Lambda'')).$$
We claim
that this map is fixed point free. Indeed, a fixed point of
$\phi_{\partial}$ can only be in the strata satisfying
$\beta''=\widetilde{\beta''}$. And we assume that there is no
nonconstant $\phi$-multiply covered pseudoholomorphic disc. Since a
$\phi$-somewhere injective disc can not be a fixed point of
$\phi_{\partial}$, the fixed point of $\phi_{\partial}$ could only
be the map that a zero energy disc bubbled off. The process of such
bubbling-off corresponds to an interior marked point moving to the
boundary. However, from the definition 
of moduli space (alternatively, we can also define the moduli space as a section of the reparameterization group action, see lemma 4.4 of \cite{So}.) we know that the marked point constrained away from $L$
can not move to the boundary. So the contradiction implies that
$\phi_{\partial}$ has no fixed point.

We below will show that $\phi_\partial$ is orientation reversing if

(i) $L$ is orientable, dim$L=3$; or

(ii) $L$ might not be orientable, dim$L=2$ and if $L$ is
nonorientable we require the
condition (\ref{9-bound}) is satisfied.\\

\noindent $\bullet$ Case (i).  $L$ is orientable and dim$L=3$\\

Since dim$L=3$, by Wu relations, $L$ is $Pin^-$, so $w_2(\V)=0$.
Since $L$ is orientable, $\mu(d'')$ is even. Note the assumption
$L\cap V=\emptyset$. Therefore, by the formula (\ref{ori-odd}), the
sign of $\Psi^{\#'}$ can be simplified as
\begin{equation}\label{12-sign}
\mathfrak{s}^{\#} =\frac{\mu(\beta'')}{2}+\deg
\mathbbm{s}''+k''+1+l''+\mathbbm{l}''\ \ \ \ \ \ {\rm mod}\ 2,
\end{equation}
Moreover, from the definition of the moduli space above we can see
that the map
$$\phi^\#:{\mathcal M}^{V,\mathbbm{s}}_{{\bf
k},\sigma,l,\rho,\mathbbm{l},\varrho}(\hat{\Sigma},X,L,\bar{\beta}',\beta'',0)\rightarrow{\mathcal
M}^{V,\mathbbm{s}}_{{\bf k},\sigma,
l,\rho,\mathbbm{l},\varrho}(\hat{\Sigma},X,L,\bar{\beta}',\widetilde{\beta''},0)$$
has the same sign.

Note that the involution $\phi$ acts on $X$ reversing the
orientation since $\phi^*\omega^{3}=-\omega^{3}$. So $\phi_\partial$
acts non-trivially on $l''$ of the factors of $X^l$. Therefore, the
sign of the map between the two fiber products
$$\phi_\partial:\partial_G{\mathcal W}^b_{\sigma,\rho,\varrho,\
\beta'+\beta''}(\Xi'',\Lambda'') \rightarrow
\partial_G{\mathcal W}^b_{\sigma, \rho,\varrho,\
\beta'+\widetilde{\beta''}}(\phi(\Xi''),\phi(\Lambda''))$$ should be
independent of $l''$.

On the other hand, let us consider general dimension $n={\rm dim}L$
for the moment. Recall (\ref{dim}) the virtual dimension of
$V$-regular moduli space is
\begin{eqnarray*}
{\rm dim}\ {\mathcal M}^{V,\mathbbm{s}}_{{\bf
k},l,\mathbbm{l}}(\Sigma, X,L,\bar{\beta})&=&\mu(\beta)+n(1-g)+|{\bf k}|\nonumber\\
   & &+2(l+\mathbbm{l}-\deg\mathbbm{s})-{\rm dim} Aut(\Sigma).
\end{eqnarray*}
The definition of invariants requires the following dimension
condition
$$\mu(\beta)+n(1-g)+|{\bf k}|+2(l+\mathbbm{l}-\deg\mathbbm{s})-{\rm dim} Aut(\Sigma)=(|{\bf k}|+2l)n+2\mathbbm{l}(n-1).$$
So we have
\begin{equation}\label{12-mud-1}
\mu(\beta)=(|{\bf k}|+2l)(n-1)+2\mathbbm{l}(n-2)+2\deg
\mathbbm{s}-n(1-g)+{\rm dim} Aut(\Sigma).
\end{equation}

We observe that if we restrict the evaluation map to $L^{k''}\times
X^{l''}\times V^{\mathbbm{l} ''}$, the image should be at least codimension
one to have a nontrivial intersection. That is
\begin{eqnarray}\label{12-dim-D}
{\rm dim}{\mathcal M}(D^2)&=&n+\mu(\beta'')+(k''+2l'')+2(\mathbbm{l}''-\deg
\mathbbm{s}'')-3\nonumber\\
&\ge& (k''+2l'')n+2\mathbbm{l}''(n-1)-1,
\end{eqnarray}
therefore
\begin{equation}\label{12-mud-2}
\mu(\beta'')\ge (k''+2l'')(n-1)+2\mathbbm{l}''(n-2)+2\deg
\mathbbm{s}''+2-n.
\end{equation}
Similarly, we have
\begin{eqnarray}\label{12-dim-Sigma}
{\rm dim}{\mathcal M}(\Sigma)&=&n(1-g)+\mu(\beta')+(k'+2l')+2(\mathbbm{l}'-\deg
\mathbbm{s}')-{\rm dim}Aut(\Sigma)\nonumber\\
&\ge& (k'+2l')n+2\mathbbm{l}'(n-1)-1,
\end{eqnarray}
and so
\begin{eqnarray}\label{12-mud-3}
\mu(\beta')&\ge& (k'+2l')(n-1)+2\mathbbm{l}'(n-2)+2\deg
\mathbbm{s}'\nonumber\\
  & &+\ {\rm dim}Aut(\Sigma)-1-n(1-g).
\end{eqnarray}
Then from (\ref{12-mud-1}) and (\ref{12-mud-3}) we have
\begin{eqnarray}\label{12-mud-4}
\mu(\beta'')&=&\mu(\beta)-\mu(\beta')\nonumber\\
         &\le&(k''+2l'')(n-1)+2\mathbbm{l}''(n-2)+2\deg
\mathbbm{s}''+1.
\end{eqnarray}

In particular, when $n=3$, from (\ref{12-mud-2}) and
(\ref{12-mud-4}), and noting that $\mu(\beta'')$ is even, we see
that each set $\partial_G{\mathcal W}^b_{\sigma,\rho,\varrho,\
\beta'+\beta''}(\Xi'',\Lambda'')$ is nonempty if and only if the
following dimension condition is satisfied
\begin{equation}\label{12-mu-beta}
\mu(\beta'')= 2(k''+2l'')+2\mathbbm{l}''+2\deg\mathbbm{s}''.
\end{equation}
Thus we have
$$sign(\phi_\partial)=\frac{\mu(\beta'')}{2}+\deg
\mathbbm{s}''+k''+1+\mathbbm{l}'' \cong 1\ \ \ \ ({\rm mod}\ 2).$$ That is
to say, the map $\phi_{\partial}$ reverses orientation. That means
$\#\partial_G{\mathcal W}=0$. Therefore, we have

\begin{eqnarray*}
0 = \#\partial\overline{\mathcal W} &=& \sum_{\beta; \vec{r'}} \#{\bf
ev}^{-1}_{(\beta,\vec{r'})}(\vec{x'},\vec{\xi'},\vec{\lambda'}\
)-\sum_{\beta;\vec{r}}\#{\bf
ev}^{-1}_{(\beta,\vec{r})}(\vec{x},\vec{\xi},\vec{\lambda}\
)+\#\partial_G{\mathcal
W}\\
&=& \mathcal {I}^{\ V,\mathbbm{s}}_{X,\phi,g,\mathbbm{d},{\bf
k},l}(\vec{x'},\vec{y'},\vec{\mathbbm{q}'})-\mathcal {I}^{\ V,\
\mathbbm{s}}_{X,\phi,g,\mathbbm{d},{\bf k},l}(\vec{x},\vec{y},
\vec{\mathbbm{q}}).
\end{eqnarray*}

So integers $\mathcal {I}^{\ V,\
\mathbbm{s}}_{X,\phi,g,\mathbbm{d},{\bf k},l}$ are independent of
the choice of $(\vec{x},\vec{y}, \vec{\mathbbm{q}})$. Equivalently,
we can say that the integral in (\ref{ROGW-int}) is independent of
the choices of $\alpha_{ai}$,  $\gamma_j$ and $\eta_\jmath$.
Similarly, we can prove that $\mathcal {I}^{\ V,\
\mathbbm{s}}_{X,\phi,g,\mathbbm{d},{\bf k},l}$ are independent of
the generic choice of $J\in{\mathcal J}_{\omega,\phi}$, and the choice
of inhomogeneous perturbation $\nu\in{\mathcal P}_{\phi,c}$. That means
$\mathcal {I}^{\ V,\ \mathbbm{s}}_{X,\phi,g,\mathbbm{d},{\bf k},l}$
are invariants of the tuple $(X,\omega,V,\phi)$. In particular, if
$d=\tilde{d}$, the numbers $\mathscr{RN}=M(V,{\bf d},\phi,{\bf k},l,
\mathbbm{l})$ in (\ref{11-RN=M}) are also
invariants for this case.\\

\noindent $\bullet$ Case (ii).  dim$L=2$, and (\ref{9-bound}) is
satisfied if $L$ is not orientable\\

Also by Wu relations, $L$ is $Pin^-$, so $w_2(\V)=0$. As the
argument above, we conclude that the sign of the map $\phi^\#$ is
given by formulas (\ref{sign-10-1}) and (\ref{sign-10-2}) depending
on whether or not $1\in\sigma$. Since $\phi^*\omega^2=\omega^2$,
$\phi$ preserves the orientation of $X$. It means that
(\ref{sign-10-1}) and (\ref{sign-10-2}) also coincide with the signs
of the map $\phi_\partial$. Then by inequalities (\ref{12-mud-2})
and (\ref{12-mud-4}), note $n=2$, we see that the stratum
$\partial_G{\mathcal W}^b_{\sigma, \rho,\varrho,\ \beta',\beta''}$ will
be empty unless
\begin{equation}\label{12-n=2}
\mu(\beta'')+r=k''+2l''+2\deg\mathbbm{s}''
\end{equation}
for $r=0$ or $-1$. The following calculation shows that the signs
(\ref{sign-10-1}) and (\ref{sign-10-2}) are exactly 1.

First, assume $1\notin\ \sigma$. Using the restriction
(\ref{12-n=2}), we have
\begin{eqnarray}\label{12-k''-1st}
\frac{k''(k''-1)}{2}&=&\frac{(\mu(\beta'')+r-2l''-2\deg\mathbbm{s}'')
(\mu(\beta'')+r-2l''-2\deg\mathbbm{s}''-1)}{2}\nonumber\\
 &\cong&
 \frac{\mu(\beta'')(\mu(\beta'')-1)}{2}+\frac{r(r-1)}{2}+l''+\deg\mathbbm{s}''+r\mu(\beta'')\
 \ ({\rm mod}\ 2).\ \ \ \ \ \ \ \ \ \
\end{eqnarray}
Also from (\ref{12-n=2}), we have
\begin{eqnarray}\label{12--ups-1}
\Upsilon^{(1)}(\beta'',k'')\cong\mu(\beta'')k'' &\cong&
\mu(\beta'')^2+r\mu(\beta'')+2(l''+\deg\mathbbm{s}'')\mu(\beta'')\nonumber\\
 &\cong& \mu(\beta'')^2+r\mu(\beta'')\ \ \ \ ({\rm mod}\ 2).
\end{eqnarray}
Substituting (\ref{12-n=2}), (\ref{12-k''-1st}) and
(\ref{12--ups-1}) into (\ref{sign-10-1}), we have
$$\mathfrak{s}^{\#(1)}_-\cong \frac{r(r+1)}{2}+\mathbbm{l}''+1\ \ \ \ ({\rm mod}\ 2).$$
Also note that when $n=2$, the involution $\phi$ acts on 1
dimensional submanifold $V$ reversing the orientation since
$\phi^*\omega=-\omega$. So $\phi_\partial$ acts non-trivially on
$\mathbbm{l}''$ of the factors of $V^\mathbbm{l}$. Therefore, the sign of the map
between the two fiber products
$$\phi_\partial:\partial_G{\mathcal W}^b_{\sigma,\rho,\varrho,\
\beta'+\beta''}(\Xi'',\Lambda'') \rightarrow
\partial_G{\mathcal W}^b_{\sigma, \rho,\varrho,\
\beta'+\widetilde{\beta''}}(\phi(\Xi''),\phi(\Lambda''))$$ should be
independent of $\mathbbm{l}''$. Thus we have
$$sign(\phi_\partial)=\frac{r(r+1)}{2}+1\cong 1\ \ \ {\rm(mod\ 2)}$$ since $r=0$ or $-1$.\\

Then we consider the case $1\in \sigma$. Using the restriction
(\ref{12-n=2}) again, we have
\begin{eqnarray}\label{12-k''}
\frac{(k''-1)(k''-2)}{2}&=&\frac{(\mu(\beta'')+r-2l''-2\deg\mathbbm{s}''-1)
(\mu(\beta'')+r-2l''-2\deg\mathbbm{s}''-2)}{2}\nonumber\\
 &\cong&
 \frac{\mu(\beta'')(\mu(\beta'')+1)}{2}+\frac{r(r+1)}{2}\nonumber\\
 & & +l''+\deg\mathbbm{s}''+r\mu(\beta'')+1\ \ \ \ ({\rm mod}\ 2).
\end{eqnarray}
Recall the condition (\ref{9-bound}): $w_1(\beta_b)\cong k_b+1$ and
(\ref{12-n=2}), we calculate
\begin{equation}\label{12-kb}
k_b(k''+1)\cong (w_1(\beta_b)+1)(\mu(\beta'')+r+1)\ \ \ \ \ \ ({\rm
mod}\ 2).
\end{equation}
Substituting (\ref{12-n=2}), (\ref{12-k''}) and (\ref{12-kb}) into
(\ref{sign-10-2}), noting that $r=0$ or $-1$, we have
\begin{eqnarray}\label{12-sign-2-substitute}
\mathfrak{s}^{\#(2)}_-&\cong&\frac{r(r+1)}{2}+r\mu(\beta'')+r+\mathbbm{l}''\nonumber\\
  & & +(w_1(\beta_b)+1)(\mu(\beta'')+r+1)+\Upsilon^{(2)}
  +w_1(\beta_b')w_1(\partial \beta'')\nonumber\\
  &\cong& r\mu(\beta'')+r+\mathbbm{l}''
  +(w_1(\beta_b)+1)(\mu(\beta'')+r+1)\nonumber\\
  & &+\Upsilon^{(2)}
  +w_1(\beta_b')w_1(\partial \beta'')
  \hspace{2.5cm} ({\rm mod}\ 2).
\end{eqnarray}
Recall the formula (\ref{Ups-2}) and the fact
$$w_1(\beta_b)=w_1(\beta_b')+w_1(\partial\beta''),$$
we can express
\begin{equation}\label{12-ups-2-eqn}
\Upsilon^{(2)}=w_1(\beta_b)(k''-1)+w_1(\partial\beta'')k'.
\end{equation}
Considering the fact that $\mu(\beta'')\cong w_1(\partial\beta'')\
({\rm mod}\ 2)$, and using (\ref{12-n=2}) and (\ref{9-bound}), we
obtain
\begin{eqnarray}\label{12-w1-k'}
w_1(\partial\beta'')k'&\cong&w_1(\partial\beta'')(k_b-k'')\nonumber\\
 &\cong& w_1(\partial\beta'')(w_1(\beta_b)+1+\mu(\beta''+r))\nonumber\\
 &\cong& w_1(\partial\beta'')w_1(\beta_b')+\mu(\beta'')(1+r)\ \ \ \ \ \ ({\rm
mod}\ 2).
\end{eqnarray}
Substituting (\ref{12-w1-k'}) and (\ref{12-n=2}) in
(\ref{12-ups-2-eqn}), we get
$$\Upsilon^{(2)}\cong w_1(\beta_b)(\mu(\beta'')+r+1)+w_1(\partial\beta'')w_1(\beta_b')
+\mu(\beta'')(1+r)\ \ ({\rm mod}\ 2).$$ Substituting the last
formula in (\ref{12-sign-2-substitute}), we calculate
\begin{eqnarray}\label{12-sign-2-good}
\mathfrak{s}^{\#(2)}_-&\cong& r\mu(\beta'')+r+\mathbbm{l}''
  +(\mu(\beta'')+r+1)+\mu(\beta'')(1+r)\nonumber\\
  &\cong&\mathbbm{l}''+1  \hspace{3.5cm} ({\rm mod}\ 2).
\end{eqnarray}
As the mentioned reason above, the  sign of $\phi_\partial$ is
independent of $\mathbbm{l}''$, so $sign(\phi_\partial)\cong 1 \ \ {\rm
(mod\ 2)}$ , which implies $\mathcal {I}^{\ V,\
\mathbbm{s}}_{X,\phi,g,\mathbbm{d},{\bf k},l}$ are invariants of the
tuple $(X,\omega,V,\phi)$. In particular, the numbers
$\mathscr{RN}=M(V,{\bf d},\phi,{\bf k},l, \mathbbm{l})$ in
(\ref{11-RN=M}) are invariants of the tuple $(X,\omega,V,\phi)$.


\renewcommand{\thesection}{Appendix}
\section{}\label{SEC-APP}
\renewcommand{\theequation}{A.\arabic{equation}}
\renewcommand{\thesubsection}{A.\arabic{subsection}}
\renewcommand{\theDefinition}{A.\arabic{Definition}}
\renewcommand{\theProposition}{A.\arabic{Proposition}}
\renewcommand{\theLemma}{A.\arabic{Lemma}}


For the convenience of the reader, we review some definitions and
important conclusions in \cite{So} about the orientation of
determinant  of real linear Cauchy-Riemann operator. For our
concrete problem of intersection of stable maps with a codimensional
two symplectic submanifold, we state that  parallel conclusions hold
for some kind of restriction of Cauchy-Riemann operator.

We denote by $\Gamma$ an appropriate Banach space completion of the
smooth section of a vector bundle. For a vector bundle $V\rightarrow
B$, we denote by $\mathfrak{F}(V)$ the orthonormal frame bundle of
$V$ which is a principal $O(n)$ bundle. Recall the fact that the Lie
group $Spin(n)$ is the central $\Z/2\Z$ extension of the special
orthogonal group $SO(n)$, similarly, the two groups $Pin^+(n)$ and
$Pin^-(n)$, although are topologically the same, are two different
central extensions of $O(n)$.
\begin{Definition}
A $Pin^\pm(n)$ structure $\mathfrak{P}=(P,p)$ on a vector bundle
$V\rightarrow B$ consists of principal $Pin^\pm(n)$ bundle
$P\rightarrow B$ and a $Pin^\pm(n)-O(n)$ equivariant bundle map
$$p:P\rightarrow\mathfrak{F}(V).$$
A map $\varphi:V\rightarrow V'$ between vector bundles with $Pin$
structure preserves $Pin$ structure if there exists a lifting
$\tilde{\varphi}$ such that the following diagram commutes
\begin{center}
\setlength{\unitlength}{1.0cm}
\begin{picture}(2.5,3)
 \put(0.1,2.0){\makebox(0.5,0.5){$P$}}
 \put(0.6,2.2){\vector(1,0){2}}
 \put(1.4,2.2){\makebox(0.4,0.4){$\tilde{\varphi}$}}
 \put(2.7,2.0){\makebox(0.5,0.5){$P'$}}
 \put(0.3,2.0){\vector(0,-1){1.0}}
 \put(2.9,2.0){\vector(0,-1){1.0}}
 \put(0.3,1.4){\makebox(0.4,0.4){$p$}}
 \put(2.95,1.4){\makebox(0.4,0.4){$p'$}}
 \put(0.1,0.4){\makebox(0.5,0.5){$\mathfrak{F}(V)$}}
 \put(2.8,0.4){\makebox(0.5,0.5){$\mathfrak{F}(V')$}}
 \put(0.9,0.6){\vector(1,0){1.5}}
 \put(1.4,0.6){\makebox(0.4,0.4){$\varphi$}}
\end{picture}
\end{center}
\end{Definition}
The obstruction to putting a $Spin$ structure on a bundle
$V\rightarrow B$ is $w_2(V)\in H^2(B,\Z/2\Z)$. For $Pin^+$ it is
still $w_2(V)$, and for $Pin^-$ it is $w_2(V)+w_1^2(V)$.\\

\begin{Definition}
Let $(\Sigma,\partial\Sigma)$ be a Riemann surface with boundary
$\partial\Sigma=\coprod_{a=1}^m(\partial\Sigma)_a$. A Cauchy-Riemann
(or Riemann-Roch) Pin boundary value problem
$$ \underline{D}=(\Sigma,E,F,\mathfrak{P},D)$$  consists of

$1^\circ$ a complex vector bundle $E\rightarrow\Sigma$,

$2^\circ$ a totally real sub-bundle over the boundary
$F\rightarrow\partial\Sigma$,

$3^\circ$ a $Pin^+$ or $Pin^-$ structure $\mathfrak{P}$ on $F$,

$4^\circ$ an orientation of $F|_{(\partial\Sigma)_a}$ for each $a$
so that $F|_{(\partial\Sigma)_a}$ is orientable,

$5^\circ$ a differential operator
$$D:\Gamma((\Sigma,\partial\Sigma),(E,F))\rightarrow\Gamma(\Sigma,\Omega^{0,1}(E))$$
such that for $f\in C^\infty(\Sigma,R)$,
$\xi\in\Gamma((\Sigma,\partial\Sigma),(E,F))$,
$$D(f\xi)=fD\xi+(\bar{\partial}f)\xi.$$
\end{Definition}
$D$ is called a real linear Cauchy-Riemann operator.\\

\noindent{\it Remark.} For the sake of studying the
$(J,\nu)$-holomorphic maps relative to a codimensional 2 symplectic
submanifold $V\subset X$, say $V$-regular maps, we may consider the
restriction of a real linear Cauchy-Riemann operator $D$ to a
subspace
$\Gamma^{rest}:=\Gamma_{(\mathbbm{s})}((\Sigma,\partial\Sigma),(E,F))$.
A section $\xi\in\Gamma^{rest}$ (or
$\Gamma_{(\mathbbm{s})}$) if and only if it satidfies
some vanishing conditions at each prescribed (say, intersection)
marked point $p_{a\imath}$ (resp. $q_\jmath$). In the present paper,
we only consider the orientability of moduli space under the
condition $L\cap V=\emptyset$, thus we denote the subspace simply by
$\Gamma^{rest}$ or $\Gamma_\mathbbm{s}$ and the restriction of $D$
to this subspace by $D^{rest}$ or $D_\mathbbm{s}$. In particular,
for the concrete case of moduli space, $D_\mathbbm{s}$ is the
restriction by contact condition of
$D_u=D_u\bar{\partial}_{(J,\nu)}$, where $u\in{\mathcal M}^V(X,L,d)$.
When the pair $(J,\nu)\in\mathbb{J}^V$ (or $\mathbb{J}_{\phi}^V$) is
$V$-compatible (see Definition \ref{3.1}), the arguments in Lemma
\ref{dim-M} show that the restriction is transverse. Therefore, in
such concrete situation $D_\mathbbm{s}$ is $D_u$ for $u\in{\mathcal
M}_{\mathbbm{s}}^V(X,L,d)$.

It is not difficult to apply the arguments by Mcduff-Salamon (see
\cite{MS} Appendix C.2 ) to show that $D^{rest}$ is a Fredholm
operator. We call $
\underline{D^{rest}}=(\Sigma,E,F,\mathfrak{P},D^{rest})$ the {\it
restricted Pin boundary value problem}.

\bigskip
For a Fredholm operator $D$, we define its determinant line by
$${\rm det}(D):=\Lambda^{\rm max}({\rm ker}\ D)\otimes\Lambda^{\rm max}({\rm coker}\ D)^*.$$
As explained in \cite{MS}, for a family of Fredholm operators, ${\rm
det}(D)$ denotes a line bundle with the natural topology.

Since we may trivialize $E$ over $\Sigma$ and each component of
boundary $(\partial\Sigma)_a\simeq S^1$, the restriction of $F$ to
each $(\partial\Sigma)_a$ defines a loop of totally real subspaces
of $\C^n$. it is well known such a loop associates a Maslov index
$\mu_a$. we denote by
\begin{equation}\label{eq-A1}
\mu(E,F)=\sum_{a=1}^m\mu_a
\end{equation}
the total Maslov index of the vector bundle pair $(E,F)$. Moreover,
$\mu(E,F)$ doesn't depend on the choice of trivialization of $E$.

We say $\underline{\varphi}:\underline{D}\rightarrow\underline{D}'$
is an {\it isomorphism} of Cauchy-Riemann (or restricted) $Pin$
boundary value problems if

i) there exists a biholomorphic map $f:\Sigma\rightarrow\Sigma'$;

ii) there exists  an isomorphism of bundles $\varphi:E\rightarrow
E'$ covering $f$ such that $\varphi|_{\partial\Sigma}$ maps $F$ to
$F'$ and preserving $Pin$ structure and preserving orientation if
$F,F'$ are orientable;

iii) $\varphi\circ D=D'\circ\varphi$.

J. Solomon showed the following (Proposition 2.8 and Lemma 2.9 in
\cite{So})
\begin{Proposition}\label{A.1.1}
The determinant line of a real-linear Cauchy-Riemann Pin boundary
value problem $\underline{D}$ admits a canonical orientation. If
$\underline{\varphi}:\underline{D}\rightarrow\underline{D}'$ is an
isomorphism, then the induced morphism $$\Psi:{\rm
det}(D)\rightarrow{\rm det}(D')$$ preserves the canonical
orientation. Furthermore, the canonical orientation varies
continuously in a family of Cauchy-Riemann operators. That is, it
defines a single component of the determinant line bundle over that
family. If the boundary condition $F|_{(\partial\Sigma)_a}$  is
orientable, then reversing the orientation on
$F|_{(\partial\Sigma)_a}$ will change the canonical orientation on
${\rm det}(D)$.
\end{Proposition}

\noindent {\it Remark.} It is not difficult to apply the method of
proof by Solomon to generalize the Proposition above to the case for
the restricted $Pin$ boundary value problem $\underline{D^{rest}}$.
For instance, it is important in the proof of Proposition 2.8 in
\cite{So} that, for any chosen Cauchy-Riemann operator $\tilde{D}$
on the restriction $\hat{E}|_{\tilde{\Sigma}}$, $\det(\tilde{D})$
has the canonical complex orientation, where
$\hat{E}\rightarrow\hat{\Sigma}$ is the degenerated vector bundle of
$E\rightarrow\Sigma$,
$\displaystyle\hat{\Sigma}=\tilde{\Sigma}\cup_a\Delta_a$
($\tilde{\Sigma}$ is the closed component, $\Delta_a$ is a disk
corresponding to each $(\partial\Sigma)_a$). In our relative case,
we consider the restriction $\tilde{D}^{rest}$ which is a complex
Fredholm operator, and $\det(\tilde{D}^{rest})$ also can be equipped
with the canonical complex orientation. The remain arguments in the
proof of Proposition 2.8 in \cite{So} can go through with minor
modifications.

\bigskip
We then come to study the sign of conjugation on the canonical
orientation of the determinant line of a Cauchy-Riemann (or
restricted) $Pin$ boundary value problem. More precisely, given a
Riemann surface $\Sigma$, let $\bar{\Sigma}$ denote the same
topological surface with conjugate complex structure, and let
$$c:\Sigma\rightarrow\bar{\Sigma}$$ denote the tautological
anti-holomorphic map. Similarly, let $(\bar{E}, \bar{F})$ denote the
same real bundle pair with the opposite complex structure on $E$,
and denote by
$$C:E\rightarrow\bar{E}$$
the tautological anti-complex-linear bundle map. Furthermore, a
Cauchy-Riemann (resp. restricted) operator $D$ (resp.
$D_\mathbbm{s}$) on the bundle $E\rightarrow\Sigma$ is the same as a
Cauchy-Riemann (resp. restricted) operator $\bar{D}$ (resp.
$\bar{D}_\mathbbm{s}$) on the bundle
$\bar{E}\rightarrow\bar{\Sigma}$. So, given any Cauchy-Riemann
(resp. restricted) $Pin$ boundary problem $\underline{D}$ (resp.
$\underline{D_\mathbbm{s}}$), we may construct its conjugate
$\underline{\bar{D}}$ (resp. $\underline{\bar{D}_\mathbbm{s}}$).
Clearly, we have a tautological map of Cauchy-Riemann (resp.
restricted) $Pin$ boundary value problems
$$\underline{C}:\underline{D}\rightarrow\underline{\bar{D}}\ \ ({\rm resp}.
\
\underline{D_\mathbbm{s}}\rightarrow\underline{\bar{D}_\mathbbm{s}}).$$

Suppose the genus of $\Sigma/\partial\Sigma$ is $g_0$, the number of
boundary components of $\Sigma$ is $m$ and rank$(F)=n$. Denote by
$\mu=\mu(E,F)$ and by $w_1$ the first Stiefel-Whitney class. Solomon
calculated the sign of the induced isomorphism $\Psi:$
$\det(D)\rightarrow$ $\det(\bar{D})$ (Proposition 2.12 in
\cite{So}).
\begin{Proposition}\label{Pro-A2}
The sign of the induced isomorphism $\Psi$ relative to the
respective canonical orientation is given by
\begin{eqnarray}\label{A2}
sign^+(\underline{D})& =&\frac{\mu(\mu+1)}{2}+(1-g_0)n+mn\nonumber\\
                      & & +\sum_{a<b}w_1(F)((\partial\Sigma)_a)w_1(F)((\partial\Sigma)_b)\
\ \ {\rm mod}\ 2,
\end{eqnarray}
for a $Pin^+$ structure and

\begin{eqnarray}\label{A3}
sign^-(\underline{D})&=&\frac{\mu(\mu+1)}{2}+(1-g_0)n+mn\nonumber\\
                     & &
+\sum_{a<b}w_1(F)((\partial\Sigma)_a)w_1(F)((\partial\Sigma)_b)\nonumber\\
                     &  &+\sum_{a=1}^mw_1(F)((\partial\Sigma)_a)\
\ \ \ \ \ \ \ \ \ \ \ \ \ \ \ \ \ \  \ {\rm mod}\ 2,
\end{eqnarray}
for a $Pin^-$ structure. In particular, when $\Sigma=D^2$, $g_0=0$
and $m=1$ we have
\begin{eqnarray}\label{A4}
sign^\pm(\underline{D})=\frac{\mu(\mu\pm 1)}{2}\ \ \ \ \ \ \ {\rm
mod}\ 2.
\end{eqnarray}
\end{Proposition}

As to our concrete case of regular maps, for the restricted $Pin$
boundary value problem $\underline{D_\mathbbm{s}}$, we study the
sign of the induced isomorphism $\Psi:$
$\det(D_\mathbbm{s})\rightarrow$ $\det(\bar{D}_\mathbbm{s})$. Recall
that when the pair $(J,\nu)\in\mathbb{J}^V$ (or
$\mathbb{J}_{\phi}^V$), we can consider $D_\mathbbm{s}=D_u$ for
$u\in{\mathcal M}_{\mathbbm{s}}^V(X,L,d)$, where
$\mathbbm{s}=(s_1,\cdots,s_\mathbbm{l})$ is the list of
multiplicities of interior intersection points of a $V$-regular map,
and ${\rm deg}\mathbbm{s}=\sum_{\jmath=1}^{\mathbbm{l}} s_\jmath$.
The argument is a modification of the one of Proposition 2.12 in
\cite{So}, involving the interior contact conditions. Basically, one
need modify the formula (9) in \cite{So} respecting to our concrete
case of $V$-regular $(J,\nu)$-maps for generic $(J,\nu)\in
\mathbb{J}^V$, the transversality of contact conditions (see Lemma
\ref{dim-M}) implies the relation
$${\rm index}_{\C}(\tilde{D}_\mathbbm{s})={\rm index}_{\C}(\tilde{D})-{\rm deg}\mathbbm{s},$$
where $\tilde{D}_\mathbbm{s}$ is the restriction of the chosen
complex Cauchy-Riemann operator $\tilde{D}$ on the restriction
$\hat{E}|_{\tilde{\Sigma}}$. Thus, we have
\begin{Proposition}\label{Pro-A3}
The  sign of the induced isomorphism $\Psi:\
\det(D_\mathbbm{s})\rightarrow \det(\bar{D}_\mathbbm{s})$ relative
to the respective canonical orientation is given by
\begin{eqnarray}\label{A2+3}
sign^+(\underline{D_\mathbbm{s}})& =&\frac{\mu(\mu+1)}{2}+{\rm deg}\mathbbm{s}+(1-g_0)n+mn\nonumber\\
                      & & +\sum_{a<b}w_1(F)((\partial\Sigma)_a)w_1(F)((\partial\Sigma)_b)\
\ \ {\rm mod}\ 2,
\end{eqnarray}
for a $Pin^+$ structure and

\begin{eqnarray}\label{A3+3}
sign^-(\underline{D_\mathbbm{s}})&=&\frac{\mu(\mu+1)}{2}+{\rm deg}\mathbbm{s}+(1-g_0)n+mn\nonumber\\
                     & &
+\sum_{a<b}w_1(F)((\partial\Sigma)_a)w_1(F)((\partial\Sigma)_b)\nonumber\\
                     &  &+\sum_{a=1}^mw_1(F)((\partial\Sigma)_a)\
\ \ \ \ \ \ \ \ \ \ \ \ \ \ \ \ \ \  \ {\rm mod}\ 2,
\end{eqnarray}
for a $Pin^-$ structure. In particular, when $\Sigma=D^2$, $g_0=0$
and $m=1$ we have
\begin{eqnarray}\label{A4+3}
sign^\pm(\underline{D_\mathbbm{s}})=\frac{\mu(\mu\pm 1)}{2}+{\rm
deg}\mathbbm{s}\ \ \ \ \ \ \ {\rm mod}\ 2.
\end{eqnarray}
\end{Proposition}
{\it Sketch of proof}.\ From Proposition \ref{A.1.1} and the remark after it, we  see that the determinant line of
the restricted real-linear Cauchy-Riemann $Pin$ boundary value problem $\underline{D_\mathbbm{s}}$ ($i.e.$ $\underline{D^{rest}}$) admits a canonical orientation. Note that the formula (7) in \cite{So} can be modified as
\begin{equation}\label{modfy-7}
\det(\#_aD_a\#\tilde{D}_\mathbbm{s})\simeq\bigotimes_a\det(D_a)\otimes\det(\tilde{D}_\mathbbm{s})
\otimes\bigotimes_a\det(E_{\hat{\gamma}_a})^*.
\end{equation}
Since conjugation on a complex  vector space will make a sign change equal to its complex dimension (mod 2), that changes the orientation of $\det(\tilde{D}_\mathbbm{s})$  according to the index $\mathrm{index}_{\C}(\tilde{D}_\mathbbm{s})={\rm index}_{\C}(\tilde{D})-\deg\mathbbm{s}$. Combining the formula (\ref{modfy-7}) and the formula (9) and arguments in Proposition 2.12 of \cite{So}, one can verify the equalities (\ref{A2+3}), (\ref{A3+3}) and (\ref{A4+3}).\qed

\begin{Definition}\label{short-exact}
A short exact sequence of families of Fredholm operators
$$0\rightarrow D'\rightarrow D\rightarrow D''\rightarrow 0$$
consists of a base parameter space $B$, two short exact sequences of
Banach space bundles over $B$
$$0\rightarrow X'\rightarrow X\rightarrow X''\rightarrow 0,\ \ \ \ \ \
0\rightarrow Y'\rightarrow Y\rightarrow Y''\rightarrow 0.$$ and
three Fredholm morphisms of Banach bundles
$$D:X\rightarrow Y,\ \ \ \ \ D':X'\rightarrow Y'\ \ \ \ \ D'':X''\rightarrow Y'',$$
such that the following diagram commutes
$$
\begin{array}{ccccccccc}
0&\longrightarrow &Y'&\longrightarrow &Y&\longrightarrow
&Y''&\longrightarrow &0 \\
 & &D'\uparrow\ \ \ \ \ \ & &D\uparrow\ \ \ \ \ & &D''\uparrow\ \ \ \ \ \ & & \\
0&\longrightarrow &X'&\longrightarrow &X&\longrightarrow
&X''&\longrightarrow &0
\end{array}
$$
\end{Definition}

\bigskip

The following lemma will be used in section 8.
\begin{Lemma}\label{A-iso}
A short exact sequence of families of Fredholm operators
$$0\rightarrow D'\rightarrow D\rightarrow D''\rightarrow 0$$induces
an isomorphism $${\rm det}(D')\otimes{\rm det}(D'')\simeq{\rm
det}(D).$$
\end{Lemma}

\bigskip


\bigskip

\small

\noindent School of Mathematical Sciences\\
Nanjing Normal University\\
Nanjing 210023, P. R. China\\
Email:\ hailongher@126.com

\end{document}